\documentclass[preprint,review,3p,12pt,sort&compress]{elsarticle}

\usepackage{amsmath}
\usepackage{amsfonts}
\usepackage{amssymb}
\usepackage{amsthm}
\usepackage{siunitx}
\usepackage{graphicx} % for includegraphics
\usepackage{xspace} % needed for \eg, \ie, \etc
\usepackage{bm} % for bold math
\usepackage{threeparttable} % for tables with footnotes
\usepackage[ruled,lined,linesnumbered]{algorithm2e}
\usepackage[english]{babel}
\usepackage{amsmath}
\usepackage{graphicx}
\usepackage[justification=centering]{caption}
\usepackage{booktabs}
\usepackage{placeins}
\usepackage{soul}
\usepackage[table]{xcolor}

\usepackage{hyperref}
\hypersetup{
 colorlinks=true,
 citecolor=blue,
 linkcolor=blue,
 urlcolor=blue}
\usepackage{tikz}
\usetikzlibrary{matrix,backgrounds,shapes}
\usetikzlibrary{positioning}
\usetikzlibrary{fit}
\usetikzlibrary{plotmarks}
\usetikzlibrary{decorations.pathreplacing}
\usepackage{pgfplots}
\usetikzlibrary{shapes}
\usetikzlibrary{positioning,arrows}
\usetikzlibrary{fadings,shapes.arrows,shadows}
\usetikzlibrary{arrows.meta}
\definecolor{darkorange}{rgb}{1.0, 0.55, 0.0}
\definecolor{Royalblue}{rgb}{0.254,0.41,0.88}
\definecolor{royalblue}{rgb}{0.254,0.41,0.88}
\usepackage{xcolor}

\newtheorem{definition}{Definition}
\newtheorem{remark}{Remark}
\usepackage{xspace}
\usepackage{amsmath}
\usepackage{amsfonts}
\usepackage{amssymb}
\usepackage{bm}
\usepackage{xspace}
\usepackage{afterpage}
\usepackage{color}
\usepackage{soul}
\usepackage{cancel}
\usepackage{gensymb}
\usepackage{tikz}
\usetikzlibrary{matrix,backgrounds,shapes}
\usetikzlibrary{positioning}
\usetikzlibrary{fit}
\usetikzlibrary{plotmarks}
\usetikzlibrary{decorations.pathreplacing}
\usepackage{pgfplots}
\usetikzlibrary{shapes}
\usetikzlibrary{positioning,arrows}
\usetikzlibrary{fadings,shapes.arrows,shadows}  
\usetikzlibrary{arrows.meta}
\definecolor{darkorange}{rgb}{1.0, 0.55, 0.0}
\definecolor{do}{rgb}{1.0, 0.55, 0.0}
\definecolor{Royalblue}{rgb}{0.254,0.41,0.88}
\usepackage{xcolor}

\definecolor{darkorange}{rgb}{1.0, 0.55, 0.0}
\definecolor{royalblue}{rgb}{0.254,0.41,0.88}

\newcommand{\Tr}{\ensuremath{^{\mr{T}}}}

\newcommand{\mr}[1]{\ensuremath{\mathrm{#1}}}
\newcommand{\fnc}[1]{\ensuremath{\mathcal{#1}}}
\newcommand{\bfnc}[1]{\ensuremath{\bm{\mathcal{#1}}}}
\newcommand{\mat}[1]{\ensuremath{\mathsf{#1}}}

\newcommand{\ie}[0]{{i.e.\@}\xspace}

\newcommand{\Th}[0]{\ensuremath{^{\mathrm{th}}}}
\newtheorem{assume}{Assumption}

\DeclareMathOperator{\diag}{diag}

\newcommand{\xm}[1]{\ensuremath{x_{#1}}}
\newcommand{\xil}[1]{\ensuremath{\xi_{#1}}}
\newcommand{\alphal}[1]{\ensuremath{\alpha_{#1}}}
\newcommand{\betal}[1]{\ensuremath{\beta_{#1}}}
\newcommand{\Nl}[1]{\ensuremath{N_{#1}}}
\newcommand{\bxil}[1]{\ensuremath{\bm{\xi}_{#1}}}
\newcommand{\bxili}[2]{\ensuremath{\bm{\xi}_{#1}^{(#2)}}}

\newcommand{\Q}[0]{\ensuremath{\bm{\fnc{Q}}}}
\newcommand{\Jk}[0]{\ensuremath{\fnc{J}}}

\newcommand{\Jdxildxm}[2]{\ensuremath{\Jk\frac{\partial\xil{#1}}{\partial\xm{#2}}}}
\newcommand{\Fxm}[1]{\ensuremath{\bm{\fnc{F}}_{\xm{#1}}}}
\newcommand{\FxmI}[1]{\ensuremath{\bm{\fnc{F}}_{\xm{#1}}^{(\mathrm{I})}}}
\newcommand{\FxmV}[1]{\ensuremath{\bm{\fnc{F}}_{\xm{#1}}^{(\mathrm{V})}}}

\newcommand{\Um}[1]{\ensuremath{\fnc{U}_{#1}}}
\newcommand{\E}[0]{\ensuremath{\fnc{E}}}

\newcommand{\GB}[0]{\ensuremath{\bm{\fnc{G}}^{(\mathrm{B})}}}
\newcommand{\Gzero}[0]{\ensuremath{\bm{\fnc{G}}^{(0)}}}

\newcommand{\DxiloneD}[1]{\ensuremath{\mat{D}_{\xil{#1}}^{(\mathrm{1D})}}}
\newcommand{\PxiloneD}[1]{\ensuremath{\mat{P}_{\xil{#1}}^{(\mathrm{1D})}}}
\newcommand{\QxiloneD}[1]{\ensuremath{\mat{Q}_{\xil{#1}}^{(\mathrm{1D})}}}
\newcommand{\SxiloneD}[1]{\ensuremath{\mat{S}_{\xil{#1}}^{(\mathrm{1D})}}}
\newcommand{\ExiloneD}[1]{\ensuremath{\mat{E}_{\xil{#1}}^{(\mathrm{1D})}}}
\newcommand{\txilalpha}[1]{\ensuremath{\bm{t}_{\alphal}}}
\newcommand{\txilbeta}[1]{\ensuremath{\bm{t}_{\betal}}}

\newcommand{\Imat}[1]{\ensuremath{\mat{I}_{#1}}}
\newcommand{\M}[0]{\ensuremath{\mat{P}}}

\newcommand{\Dxil}[1]{\ensuremath{\mat{D}_{\xil{#1}}}}

\newcommand{\matFxm}[3]{\ensuremath{\mat{F}_{\xm{#1}}\left(#2,#3\right)}}

\newcommand{\qk}[1]{\ensuremath{\bm{q}_{#1}}}

\newcommand{\ones}[1]{\ensuremath{\bm{1}_{#1}}}

\newcommand{\wk}[1]{\ensuremath{\bm{w}_{#1}}}

\newcommand{\matJk}[1]{\ensuremath{{\color{black}\mat{J}}_{#1}}}

\newcommand{\matAlmk}[3]{\ensuremath{\left[{\color{black}\fnc{J}\frac{\partial\xil{#1}}{\partial\xm{#2}}}\right]_{#3}}}

\newcommand{\Ok}[0]{\ensuremath{\Omega_{\kappa}}}
\newcommand{\pOk}[0]{\ensuremath{\partial\Omega_{\kappa}}}

\newcommand{\Ohatk}[0]{\ensuremath{\widehat{\Omega}_{\kappa}}}

\newcommand{\pOhatk}[0]{\ensuremath{\partial\widehat{\Omega}_{\kappa}}}

\newcommand{\uk}[0]{\ensuremath{\bm{u}_{\kappa}}}

\newcommand{\Cij}[2]{\ensuremath{\mat{C}_{#1,#2}}}
\newcommand{\Chatij}[2]{\ensuremath{\widehat{\mat{C}}_{#1,#2}}}

\newcommand{\eonel}[1]{\ensuremath{\bm{e}_{1}}}
\newcommand{\eNl}[1]{\ensuremath{\bm{e}_{N}}}
\newcommand{\bmxi}[1]{\ensuremath{\bm{\xi}^{(#1)}}}
\newcommand{\U}[0]{\ensuremath{\fnc{U}}}

\newcommand{\GBlin}[0]{\ensuremath{{\fnc{G}}^{(\mathrm{B})}}}
\newcommand{\Gzerolin}[0]{\ensuremath{{\fnc{G}}^{(0)}}}
\newcommand{\am}[1]{\ensuremath{a_{#1}}}
\newcommand{\Bm}[1]{\ensuremath{b_{#1}}}
\newcommand{\Chatla}[2]{\ensuremath{\widehat{\fnc{C}}_{#1,#2}}}
\newcommand{\Thetaa}[1]{\ensuremath{\Theta_{#1}}}
\newcommand{\matChatla}[2]{\ensuremath{\left[{\color{black}\Chatla{#1}{#2}}\right]}}

\newcommand{\thetaa}[1]{\ensuremath{\bm{\theta}_{#1}}}
\newcommand{\IP}[0]{\ensuremath{\bm{I}_{P}}}

\renewcommand{\equiv}{=}
\newcommand{\T}{\ensuremath{\mathcal{T}}}

\newcommand{\Rey}{\mathrm{Re}}
\newcommand{\Ma}{\mathrm{Ma}}

\usepackage{cleveref}

\newlength{\mywidthgraph}

\newlength{\mywidthsubfig}
\usepackage{tikz,pgfplots,pgfplotstable}
\pgfplotsset{compat=1.13}
\definecolor{mycolor11}{rgb}{1.0, 0.55, 0.0}
\definecolor{mycolor2}{rgb}{0.254,0.41,0.88}
\definecolor{mycolor1}{rgb}{0.8, 0.0, 0.0}
\definecolor{mycolor3}{rgb}{0.825,0.11,0.50}
\definecolor{mycolor4}{rgb}{0.0, 0.7, 0.0}
\definecolor{mycolor5}{rgb}{0.54,0.81,0.21}
\usepackage{bbding}
\usepackage{tabu}
\usepackage{booktabs}% for better rules in the table
\usepackage{graphicx,subcaption}
\newcommand{\sublabelsty}[1]{\textbf{\small{(#1)}}}
\usepackage[switch, modulo]{lineno}

\begin{document}

%\linenumbers
\modulolinenumbers[1]
\begin{frontmatter}

%%%%%%%%%%%%%%%%%%%%%%%%%%%%%%%%%%%%%%%%%%%%%%%%%%%%%%%%%%%%%%%%%%%%%%%%%%%%%%%%
\title{On the robustness and performance of entropy stable discontinuous collocation methods
  for the compressible Navier--Stokes equations}

%%%%%%%%%%%%%%%%%%%%%%%%%%%%%%%%%%%%%%%%%%%%%%%%%%%%%%%%%%%%%%%%%%%%%%%%%%%%%%%%
\author[KAUST1]{Diego Rojas\fnref{fn1}}
\ead{diego.rojasblanco@kaust.edu.sa}
\author[KAUST]{Radouan Boukharfane\fnref{fn2}}
\ead{radouan.boukharfane@kaust.edu.sa}
\author[KAUST]{Lisandro Dalcin\fnref{fn3}}
\ead{dalcinl@gmail.com}
\author[nia,nasa]{David C.~Del Rey Fern\'andez\fnref{fn3}}
\ead{dcdelrey@gmail.com}
\author[KAUST]{Hendrik Ranocha\fnref{fn2}}
\ead{hendrik.ranocha@kaust.edu.sa}
\author[KAUST]{David E. Keyes\fnref{fn4}}
\ead{david.keyes@kaust.edu.sa}
\author[KAUST]{Matteo Parsani\fnref{fn5}\corref{cor1}}
\ead{matteo.parsani@kaust.edu.sa}

\cortext[cor1]{Corresponding author}

\fntext[fn1]{Ph.D. student}
\fntext[fn2]{Postdoctoral Fellow}
\fntext[fn3]{Research Scientist}
\fntext[fn4]{Professor and ECRC Director}
\fntext[fn5]{Assistant Professor}

\address[KAUST1]{King Abdullah University of Science and Technology (KAUST), Physical Science and Engineering Division (PSE), Extreme Computing Research Center (ECRC), Thuwal, Saudi Arabia}
\address[KAUST]{King Abdullah University of Science and Technology (KAUST), Computer Electrical and Mathematical Science and Engineering Division (CEMSE), Extreme Computing Research Center (ECRC), Thuwal, Saudi Arabia}
\address[nia]{National Institute of Aerospace, Hampton, Virginia, United States}
\address[nasa]{Computational AeroSciences Branch, NASA Langley Research Center,
Hampton, Virginia, United States}

%%%%%%%%%%%%%%%%%%%%%%%%%%%%%%%%%%%%%%%%%%%%%%%%%%%%%%%%%%%%%%%%%%%%%%%%%%%%%%%%
\begin{abstract}
In computational fluid dynamics, the demand for increasingly multidisciplinary reliable simulations, for both analysis and design optimization purposes, requires transformational advances in individual components of future solvers.
At the algorithmic level, hardware compatibility and efficiency are of paramount importance in determining viability at exascale and beyond.
However, equally important (if not more so) is algorithmic robustness with minimal user intervention, which 
becomes progressively more challenging to achieve as problem size and physics complexity increase.
We numerically show that low and high order entropy stable discontinuous spatial discretizations based on summation-by-part operators and simultaneous-approximation-terms technique provide an essential step toward a truly enabling technology in terms of reliability and robustness for both under-resolved turbulent flow simulations and flows with discontinuities.
\end{abstract}

%%%%%%%%%%%%%%%%%%%%%%%%%%%%%%%%%%%%%%%%%%%%%%%%%%%%%%%%%%%%%%%%%%%%%%%%%%%%%%%%

\end{frontmatter}

%%%%%%%%%%%%%%%%%%%%%%%%%%%%%%%%%%%%%%%%%%%%%%%%%%%%%%%%%%%%%%%%%%%%%%%%%%%%%%%%
\section{Introduction}
Efficient numerical algorithms are sought that exploit $\mathcal{O}(10^{9})$ flops $10^{9}$ times per second, or exaflop/s on next-generation data-centric hardware.
Hardware compatibility and efficiency are of paramount importance in determining an algorithm's viability at exascale.
Equally important (if not more so) is algorithmic robustness, which becomes progressively more challenging to achieve as problem size and physics complexity increase. 
The requirement is that every step of the solution chain executes high level of reliability/robustness to minimize user intervention. 
In computational fluid dynamics (CFD), very compact high-order accurate methods are natural candidates for next-generation hardware because they are accurate, and their ratio of communications to local computations is reduced relative to low order methods of the same accuracy (see, for instance, \cite{hadri_ccpe_2019} for sustained petascale production turbulent flow runs).

Among modern, unstructured high order methods we can mention discontinuous collocation (DC), discontinuous Galerkin (DG), spectral difference (SD), and flux reconstruction (FR) methods, which can produce highly accurate solutions with minimum numerical dispersion and dissipation.
Although DC, DG, SD, and FR methods are well suited for smooth solutions, numerical instabilities may occur if the flow contains \textit{under-resolved physical features} (e.g., under-resolved turbulent flows) or \textit{discontinuities} (e.g., shocks).
A variety of mathematical stabilization strategies are commonly used to alleviate this problem, e.g., filtering \cite{hesthaven_2008_nodal_dg}, artificial dissipation, polynomial de-aliasing through over-integration \cite{gassner_underresolved_turbulence_2013,mengaldo_dealiasing_2015}, and weighted essentially non-oscillatory limiters \cite{zhu_dg_weno_2013}, to cite a few.
However, such stabilization techniques possess several drawbacks since i) they reduce accuracy \cite{wang_high_order_workshop_2013}, ii) they usually require tuning parameters for each problem configuration, and iii) they do not yet possess rigorous stability proofs.
Thus, the use of high-order accurate methods for complex flow applications is still problematic, and most commercial and industrial software rely on robust nominally second-order accurate discretizations.

Over the past few years, there have been rapid developments in entropy stable high order methods, which can be proven rigorously to be nonlinearly stable (entropy stable).
These discretizations are expected to be an essential component in future CFD solvers for complex practical flow simulations \cite{nasa_2030_vision}.
High order entropy stable schemes are often based on the well known matrix-vector nodal formulation collocated at quadrature points; see, for instance, \cite{hesthaven_2008_nodal_dg}.
Because of the approximation error induced by quadrature, we no longer have the integration by parts property and the chain rule in all cases \cite{ranocha2019mimetic}.
However, since integration by parts is an essential ingredient for stability proofs at the continuous level, it is necessary to transfer this property to the discrete level. 
This is precisely the design goal of summation-by-parts (SBP) methods.
SBP operators were developed at first in the context of finite difference schemes \cite{kreiss1974finite} and later transferred to other frameworks such as finite volume \cite{nordstrom2001finite}, flux reconstruction \cite{ranocha2016summation}, and discontinuous Galerkin schemes \cite{gassner2013skew,carpenter_ssdc_2014}. The methods used in this article belong to the latter class.
In particular, they are based on collocated Legendre--Gauss--Lobatto (LGL) quadrature rules in one space dimension and the corresponding discrete SBP operators on curvilinear, unstructured tensor product elements \cite{carpenter_ssdc_2014,parsani_entropy_stability_solid_wall_2015,parsani_ssdc_staggered_2016}.
The discrete integral and derivative operators associated with this quadrature were shown to satisfy the SBP property (see \cite{Fernandez2014,Svard2014} for a review of SBP operators).

While SBP operators can be viewed as matrix difference operators that are mimetic of integration by parts, additional techniques are necessary to compensate for the lack of the chain rule.
For example, ad hoc split form methods have been provided for some PDEs such as Burgers' equation and the compressible Euler equations \cite{carpenter_ssdc_2014,carpenter_entropy_stable_staggered_2015,gassner_split_form_sbp_2016}.
In \cite{FisherCarpenter2013JCPb,carpenter_ssdc_2014,parsani_entropy_stability_solid_wall_2015}, Carpenter and co-authors demonstrated the generic logic behind the splitting procedure by showing the flux differencing technique with the telescoping property, i.e., a telescoping flux form at the element level (see also \cite{carpenter_entropy_stability_ssdc_2016}).
Flux differencing is essentially a high order difference operation on Tadmor's entropy conservative fluxes \cite{Tadmor2003}, and applies to any system with any given entropy function.
Discrete stability over the whole domain is achieved by combining the SBP operator with suitable inter-element coupling procedures and boundary conditions, e.g., the simultaneous-approximation-terms (SATs); see, for instance, \citet{parsani_entropy_stable_interfaces_2015,parsani_entropy_stability_solid_wall_2015,svard_entropy_stable_euler_wall_2014}.

Tadmor's basic idea has led to the construction of several high order and low order entropy stable schemes (see, for instance, \cite{Fjordholm2012,Ray2016}).
An alternative approach, developed by \citet{Olsson1994}, \citet{Margot1996} and \citet{Yee2000} (see also \cite{Sandham2002,Bjorn2018}), relies on choosing entropy functions that result in a homogeneity property on the compressible Euler fluxes.
By using this property, splitting of the compressible Euler fluxes are constructed such that when contracted with the entropy variables result in stability estimates analogous in form to energy estimates obtained for linear PDEs.
Thus, discretizing the resulting split form using SBP operators, the nonlinear stability analysis performed at the continuous level is mimicked at the semi-discrete level.
However, because of the choice of the entropy functions, these approaches cannot be used for the compressible Navier--Stokes equations.

A complementary and general extension of Tadmor's ideas to finite domains was initiated by Fisher and co-workers \cite{Fisher2013,FisherCarpenter2013JCPb} who combined the SBP framework, using classical finite difference SBP operators, with Tadmor's two-point entropy conservative flux.
A key feature of this approach is that it extends directly to the compressible Navier--Stokes equations.
The resulting schemes follow the continuous entropy analysis, can be shown to be entropy conservative, and can be made entropy stable by adding appropriate interface dissipation in multi-block domains.
This nonlinearly stable approach inherits all of the mechanics of SBP-SAT schemes for the imposition of boundary conditions and inter-element coupling.
Therefore, it gives a systematic methodology for discretizing problems on complex geometries \cite{parsani_entropy_stable_interfaces_2015,parsani_entropy_stability_solid_wall_2015}.
Moreover, by constructing schemes that are discretely mimetic of the continuous stability analysis, the need to assume exact integration in the stability proofs is eliminated (see, for example, the work of \citet{Hughes1986}).
These ideas have been extended to include discontinuous collocated spectral elements \cite{carpenter_ssdc_2014,parsani_entropy_stability_solid_wall_2015}, fully- and semi-staggered discontinuous collocated spectral elements \cite{parsani_ssdc_staggered_2016,carpenter_entropy_stability_ssdc_2016}, Cartesian, semi-staggered, discontinuous collocated spectral elements with $p$-refinement \cite{carpenter_entropy_stability_ssdc_2016}, WENO spectral collocation \cite{Yamaleev2017}, multidimensional SBP operators~\cite{crean_entropy_stable_sbp_curvilinear_euler,Chen2017}, multidimensional staggered SBP operators~\cite{Fernandez2019_staggered}, modal decoupled SBP operators \cite{Chan2018},
$p$- and $hp$-adaptive discontinuous collocated spectral elements \cite{fernandez2019entropy,fernandez_entropy_stable_p_ns_2019,fernandez_entropy_stable_hp_ref_snpdea_2019}, 
and fully discrete entropy stable schemes \cite{Friedrich2019,ranocha2019relaxation}, as well as to a number of PDEs besides the compressible Euler and Navier--Stokes equations (for example, the
magneto-hydrodynamics~\cite{Winters2017} and the shallow water equations~\cite{Winters2015}).

This paper aims to shed some light on the robustness and cost of high-order accurate entropy stable discontinuous collocation discretizations for compressible viscous flows.
The literature that reports robustness studies for standard and entropy stable DC/DG formulations in CFD is scarce and focuses exclusively on inviscid flows \cite{gassner_split_form_sbp_2016,winters_split_form_2018,pazner_es_line_dg_2019} or low-speed viscous flows \cite{flad2017use, klose2019robustness}.
Hence, a detailed analysis of the numerical simulation of compressible viscous flows with under-resolved physical features or discontinuities is needed, given that the additional dissipation introduced by the viscous terms, despite obvious expectations, may not help in resolving robustness issues.

The paper is organized as follows.  
In \Cref{sec:notation}, we introduce the notation that is extensively used in the article. 
\Cref{sec:adv_diff_eq} presents the coordinate transformation from physical to computational space and key elements of the general spatial discretization framework in the context of the linear convection-diffusion equation.
In \Cref{sec:compressible_nse}, we briefly introduce the compressible Navier--Stokes equations and give an overview of three algorithms that we use to solve them numerically. 
These schemes cover entropy stable, split form, and conventional discretization choices.
\Cref{sec:numerical_results} presents and discusses the numerical results including accuracy, cost, and robustness of these three discretizations.
The test cases used for the study are the propagation of a three-dimensional (3D) isentropic vortex at $\Ma_{\infty}=0.5$, a 3D supersonic viscous problem constructed with the method of manufactured solutions at $\Rey=4\cdot 10^6$ and $\Ma\approx 2.14$, the Taylor--Green vortex at a Reynolds number of $\Rey=1{,}600$ and $\Ma=0.05$, a 3D simulation of homogeneous isotropic turbulence at $\Rey_\lambda=192$ and $\Ma_t=0.62$ with the formation of shocklets, and a 3D supersonic flow past a rod with square cross-section at $\Rey_{\infty}=10{,}000$ and $\Ma_{\infty}=1.5$.
Conclusions are drawn in \Cref{sec:conclusion}.

%%%%%%%%%%%%%%%%%%%%%%%%%%%%%%%%%%%%%%%%%%%%%%%%%%%%%%%%%%%%%%%%%%%%%%%%%%%%%%%%
\section{Notation}\label{sec:notation}

Partial differential equations (PDEs) are discretized on tensor-product cells having Cartesian computational coordinates denoted by the triple $(\xil{1},\xil{2},\xil{3})$, where the physical coordinates are denoted by the triple $(\xm{1},\xm{2},\xm{3})$. 
Vectors are represented by lowercase bold font, for example $\bm{u}$, while matrices are represented using sans-serif font, for example, $\mat{B}$. 
Continuous functions on a space-time domain are denoted by capital letters in script font.  
For example, 
\begin{equation*}
\fnc{U}\left(\xil{1},\xil{2},\xil{3},t\right)\in L^{2}\left(\left[\alphal{1},\betal{1}\right]\times
\left[\alphal{2},\betal{2}\right]\times\left[\alphal{3},\betal{3}\right]\times\left[0,T\right]\right)
\end{equation*}
represents a square integrable function, where $t$ is the temporal coordinate. 
The restriction of such function onto a set of mesh nodes is denoted by lower case bold font. 
For example, the restriction of $\fnc{U}$ onto a grid of $\Nl{1}\times\Nl{2}\times\Nl{3}$ nodes is given by the vector
\begin{equation*}
\bm{u} = \left[\fnc{U}\left(\bxili{}{1},t\right),\dots,\fnc{U}\left(\bxili{}{N},t\right)\right]\Tr,
\end{equation*}
where $N$ is the total number of nodes ($N\equiv\Nl{1}\Nl{2}\Nl{3}$), and the square brackets are used 
to delineate vectors and matrices, as well as ranges for variables (the context will make clear which meaning is being used). 
Moreover, $\bm{\xi}$ is a vector of vectors constructed from the three vectors $\bxil{1}$, $\bxil{2}$, and $\bxil{3}$, which are vectors of size $\Nl{1}$, $\Nl{2}$, and $\Nl{3}$ and contain the coordinates of the mesh in the three computational directions, respectively. 
Finally, $\bxil{}$ is constructed as 
\begin{equation*}
\bxil{}(3(i-1)+1:3i)\equiv  \bxili{}{i}
\equiv\left[\bxil{1}(i),\bxil{2}(i),\bxil{3}(i)\right]\Tr,
\end{equation*}
where the notation $\bm{u}(i)$ means the $i\Th$ entry of the vector $\bm{u}$ and $\bm{u}(i:j)$ is the subvector constructed from $\bm{u}$ using the $i\Th$ through $j\Th$ entries (\ie, Matlab notation is used).

Oftentimes, monomials are discussed and the following notation is used:
\begin{equation*}
\bxil{l}^{j} \equiv \left[\left(\bxil{l}(1)\right)^{j},\dots,\left(\bxil{l}(\Nl{l})\right)^{j}\right]\Tr,
\end{equation*}
with the  convention that $\bxil{l}^{j}=\bm{0}$ for $j<0$.

Herein, one-dimensional SBP operators are used to discretize derivatives.
The definition of a one-dimensional SBP operator in the $\xil{l}$ direction, $l=1,2,3$, is \cite{DCDRF2014,Fernandez2014,Svard2014}
\begin{definition}\label{SBP}
\textbf{Summation-by-parts operator for the first derivative}: A matrix operator, $\DxiloneD{l}\in\mathbb{R}^{\Nl{l}\times\Nl{l}}$, is an SBP operator of degree $p$ approximating the derivative $\frac{\partial}{\partial \xil{l}}$ on the domain $\xil{l}\in\left[\alphal{l},\betal{l}\right]$ with nodal distribution $\bxil{l}$ having $\Nl{l}$ nodes, if
\begin{enumerate}
\item $\DxiloneD{l}\bxil{l}^{j}=j\bxil{l}^{j-1}$, $j=0,1,\dots,p$;
\item $\DxiloneD{l}\equiv\left(\PxiloneD{l}\right)^{-1}\QxiloneD{l}$, where the norm matrix, $\PxiloneD{l}$, is symmetric positive definite;
\item $\QxiloneD{l}\equiv\left(\SxiloneD{l}+\frac{1}{2}\ExiloneD{l}\right)$, $\SxiloneD{l}=-\left(\SxiloneD{l}\right)\Tr$, $\ExiloneD{l}=\left(\ExiloneD{l}\right)\Tr$, \\ 
$\ExiloneD{l} = \diag\left(-1,0,\dots,0,1\right)=\eNl{l}\eNl{l}\Tr-\eonel{l}\eonel{l}\Tr$, 
$\eonel{l}\equiv\left[1,0,\dots,0\right]\Tr$, and $\eNl{l}\equiv\left[0,0,\dots,1\right]\Tr$. 
\end{enumerate}
Thus, a degree $p$ SBP operator is one that differentiates exactly monomials up to degree $p$.
\end{definition}

In this work, one-dimensional SBP operators are extended to multiple dimensions using tensor products ($\otimes$).  
The tensor product between the matrices $\mat{A}$ and $\mat{B}$ is given as $\mat{A}\otimes\mat{B}$. 
When referencing individual entries in a matrix the notation $\mat{A}(i,j)$ is used, which means 
the $i\Th$ $j\Th$ entry in the matrix $\mat{A}$.

The focus of this paper is exclusively on diagonal-norm SBP operators. 
Moreover, the same one-dimensional SBP operator is used in each direction, each operating on $N_l$ nodes. 
Specifically, diagonal-norm SBP operators constructed on the Legendre--Gauss--Lobatto (LGL) nodes are used, \ie, a discontinuous collocated spectral element approach is utilized (see, for instance, \cite{carpenter_ssdc_2014,parsani_entropy_stability_solid_wall_2015,parsani_ssdc_staggered_2016,gassner_entropy_shallow_water_2016,gassner_split_form_sbp_2016}).

When solving PDEs numerically, the physical domain $\Omega\subset\mathbb{R}^{3}$, with boundary $\Gamma\equiv\partial\Omega$, with Cartesian coordinates $\left(\xm{1},\xm{2},\xm{3}\right)\subset\mathbb{R}^{3}$, 
is partitioned into $K$ non-overlapping elements. 
The domain of the $\kappa^{\text{th}}$ element is denoted by $\Ok$ and has boundary $\pOk$. Numerically, we solve PDEs in computational coordinates, $\left(\xil{1},\xil{2},\xil{3}\right)\subset\mathbb{R}^{3}$, 
where each $\Ok$ is locally transformed to the reference element $\Ohatk$, with boundary $\pOhatk$,
using a pull-back curvilinear coordinate transformation which satisfies the following assumption:
\begin{assume}\label{assume:curv}
Each element in physical space is transformed using a local and invertible curvilinear coordinate transformation that is compatible at shared interfaces, meaning that the push-forward element-wise mappings are continuous across physical element interfaces.
Note that this is the standard assumption requiring that the curvilinear coordinate transformation is water-tight.
\end{assume}

Precisely, one maps from the reference coordinates $\left(\xil{1},\xil{2},\xil{3}\right) \in [-1,1]^{3}$
to the physical element by the push-forward transformation 
\begin{equation}\label{eq:3d_mapping}
\left(\xm{1},\xm{2},\xm{3}\right) = X \left(\xil{1},\xil{2},\xil{3}\right),
\end{equation}
which, in the presence of curved elements, is usually a high-order degree polynomial.

%%%%%%%%%%%%%%%%%%%%%%%%%%%%%%%%%%%%%%%%%%%%%%%%%%%%%%%%%%%%%%%%%%%%%%%%%%%%%%%%
\section{The linear convection-diffusion equation}\label{sec:adv_diff_eq}

Many of the technical details for constructing conservative and stable discretizations for the compressible Navier--Stokes can be presented in the simple context of the linear convection-diffusion equation. 
The linear convection-diffusion equation in Cartesian physical coordinates is given as  
\begin{equation}\label{eq:cartconvectiondiffusion}
\begin{aligned}
& \frac{\partial\U}{\partial t}+\sum\limits_{m=1}^{3}\frac{\partial\left(\am{m}\U\right)}{\partial\xm{m}}=
\sum\limits_{m=1}^{3}\frac{\partial^{2}(\Bm{m}\U)}{\partial\xm{m}^{2}},
&& \forall \left(\xm{1},\xm{2},\xm{3}\right)\in\Omega,\quad t\ge 0,\\
&  \U\left(\xm{1},\xm{2},\xm{3},t\right)=\GBlin\left(\xm{1},\xm{2},\xm{3},t\right), 
&& \forall \left(\xm{1},\xm{2},\xm{3}\right)\in\Gamma,\quad t\ge 0,\\
&  \U\left(\xm{1},\xm{2},\xm{3},0\right)=\Gzerolin\left(\xm{1},\xm{2},\xm{3}\right), 
&& \forall \left(\xm{1},\xm{2},\xm{3}\right)\in\Omega,
\end{aligned}
\end{equation}
where $(\am{m}\U)$ are the inviscid fluxes, $\am{m}$ are the (constant) components of the convection speed, 
$\frac{\partial(\Bm{m}\U)}{\partial\xm{m}}$ are the viscous fluxes, and $\Bm{m}$ are the (constant and positive) diffusion coefficients.  
The boundary data, $\GBlin$, and the initial condition, $\Gzerolin$, are assumed to be in $L^{2}(\Omega)$,
with the further assumption that $\GBlin$ is prescribed so that either energy conservation or energy stability is achieved.

Since derivatives are approximated with differentiation operators defined in computational space, we use the Jacobian of the push-forward mapping and the chain rule
\begin{equation*}
\frac{\partial}{\partial\xm{m}}
=\sum\limits_{l=1}^{3}\frac{\partial\xil{l}}{\partial\xm{m}}\frac{\partial}{\partial\xil{l}},\quad
\frac{\partial^{2}}{\partial\xm{m}^{2}}
=\sum\limits_{l,a=1}^{3}\frac{\partial\xil{l}}{\partial\xm{m}}
\frac{\partial}{\partial\xil{l}}\left(
\frac{\partial\xil{a}}{\partial \xm{m}}\frac{\partial}{\partial\xil{a}}  
\right),
\end{equation*} 
to transform \cref{eq:cartconvectiondiffusion} from physical to computational space as
\begin{equation}\label{eq:convectiondiffusionchain}
\Jk\frac{\partial\U}{\partial t}
+\sum\limits_{l,m=1}^{3}\Jk\frac{\partial\xil{l}}{\partial\xm{m}}
\frac{\partial \left(a_{m}\U\right)}{\partial\xil{l}}
=\sum\limits_{l,a,m=1}^{3}
\Jk\frac{\partial\xil{l}}{\partial\xm{m}}\frac{\partial}{\partial\xil{l}}
\left(\frac{\partial\xil{a}}{\partial \xm{m}}
\frac{\partial(\Bm{m}\U)}{\partial\xil{a}}
\right),
\end{equation}
where $\Jk$ is the determinant of the metric Jacobian.
Bringing the metric terms $\Jdxildxm{l}{m}$ inside the derivative, and using the product rule, gives
\begin{equation}\label{eq:convectiondiffusionstrong1}
\begin{split}
\Jk\frac{\partial\U}{\partial t}+\sum\limits_{l,m=1}^{3}
\frac{\partial}{\partial\xil{l}}\left(\Jdxildxm{l}{m}a_{m}\U\right)
-&\sum\limits_{l,m=1}^{3}a_{m}\U\frac{\partial}{\partial\xil{l}}\left(\Jdxildxm{l}{m}\right)
=\\
\sum\limits_{l,a,m=1}^{3}
\frac{\partial}{\partial\xil{l}}
\left(\Jk\frac{\partial\xil{l}}{\partial\xm{m}}\frac{\partial\xil{a}}{\partial \xm{m}}
\frac{\partial(\Bm{m}\U)}{\partial\xil{a}}\right)
-&\sum\limits_{l,a,m=1}^{3}
\frac{\partial\xil{a}}{\partial \xm{m}}
\frac{\partial(\Bm{m}\U)}{\partial\xil{a}}
\frac{\partial}{\partial\xil{l}}\left( \Jk\frac{\partial\xil{l}}{\partial\xm{m}}\right).
\end{split}
\end{equation}
The last terms on the left- and right-hand sides of \cref{eq:convectiondiffusionstrong1} are zero via the GCL relations
\begin{equation}\label{eq:GCL}
\sum\limits_{l=1}^{3}\frac{\partial}{\partial\xil{l}}\left(\Jdxildxm{l}{m}\right)=0,\quad m=1,2,3,
\end{equation}
leading to the strong conservation form of the convection-diffusion equation in computational space
\begin{equation}\label{eq:convectiondiffusionstrong}
\Jk\frac{\partial\U}{\partial t}+\sum\limits_{l,m=1}^{3}
\frac{\partial}{\partial\xil{l}}\left(\Jdxildxm{l}{m}a_{m}\U\right)= \sum\limits_{l,a,m=1}^{3}
\frac{\partial}{\partial\xil{l}}
\left(\Jk\frac{\partial\xil{l}}{\partial\xm{m}}\frac{\partial\xil{a}}{\partial \xm{m}}
\frac{\partial(\Bm{m}\U)}{\partial\xil{a}}\right).
\end{equation}

Now, consider discretizing \cref{eq:convectiondiffusionstrong} by using the following differentiation matrices
\begin{equation*}
\Dxil{1}\equiv\DxiloneD{1}\otimes\Imat{N_2}\otimes\Imat{N_3},\;
\Dxil{2}\equiv\Imat{N_1}\otimes\DxiloneD{2}\otimes\Imat{N_3},\;
\Dxil{3}\equiv\Imat{N_1}\otimes\Imat{N_2}\otimes\DxiloneD{3},
\end{equation*} 
where $\Imat{N_l}$ is an $N_l\times N_l$ identity matrix and $N_l$ is the number of LGL points per direction in a given element. 
The diagonal matrix containing the metric Jacobian is defined as
\begin{equation*}
\matJk{\kappa}\equiv\diag\left(\Jk(\bmxi{1}),\dots,\Jk(\bmxi{\Nl{\kappa}})\right),
\end{equation*}
while the diagonal matrix of the metric terms, $\matAlmk{l}{m}{\kappa}$, has to be chosen to be a discretization of
\begin{equation*}
\diag\left(\Jdxildxm{l}{m}(\bmxi{1}),\dots, \Jdxildxm{l}{m}(\bmxi{\Nl{\kappa}})\right),
\end{equation*}
where $\Nl{\kappa}\equiv N_1 N_2 N_3$ is the total number of nodes in element $\kappa$.
Using this nomenclature, the discretization of \cref{eq:convectiondiffusionstrong} on the $\kappa\Th$ element reads
\begin{equation}\label{eq:convectionstrongdisc}
\begin{split}
&\matJk{\kappa}\frac{\mr{d}\uk}{\mr{d}t}+\sum\limits_{l,m=1}^{3}a_{m}\Dxil{l}\matAlmk{l}{m}{\kappa}\uk=
\sum\limits_{l,m,a=1}^{3}\Bm{m}\Dxil{l}\matJk{\kappa}^{-1}\matAlmk{l}{m}{\kappa}\matAlmk{a}{m}{\kappa}\Dxil{a}\uk
+\bm{\mathrm{SAT}}_{\kappa},
\end{split}
\end{equation}
where $\bm{\mathrm{SAT}}_{\kappa}$ is the vectors of the SATs used to impose boundary conditions and inter-element connectivity \cite{parsani_entropy_stable_interfaces_2015,carpenter_entropy_stable_staggered_2015}. 
The $\bm{\mathrm{SAT}}_{\kappa}$ vector is in general composed from inviscid and viscous contributions, \ie $\bm{\mathrm{SAT}}_{\kappa} = \bm{\mathrm{SAT}}^{(\mathrm{I})}_{\kappa} + \bm{\mathrm{SAT}}^{(\mathrm{V})}_{\kappa}$.

Unfortunately, the scheme \eqref{eq:convectionstrongdisc} is not guaranteed to be stable. 
However, a well-known remedy is to canonically split the inviscid terms into one half of the inviscid terms in~\eqref{eq:convectiondiffusionchain} and one half of the inviscid terms in~\eqref{eq:convectiondiffusionstrong1} (see, for instance, \cite{carpenter_entropy_stable_staggered_2015}), while the viscous terms are treated in strong conservation form.  
In the continuum, this process leads to
\begin{equation}\label{eq:convectiondiffusionsplit}
\begin{split}
&\Jk\frac{\partial\U}{\partial t}+\frac{1}{2}\sum\limits_{l,m=1}^{3}\left\{
\frac{\partial}{\partial\xil{l}}\left(\Jdxildxm{l}{m}a_{m}\U\right)+
\Jdxildxm{l}{m}\frac{\partial}{\partial\xil{l}}\left(a_{m}\U\right)
\right\}\\
&-\frac{1}{2}\sum\limits_{l,m=1}^{3}\left\{
a_{m}\U\frac{\partial}{\partial\xil{l}}\left(\Jdxildxm{l}{m}\right)\right\}=
\sum\limits_{l,a,m=1}^{3}
\frac{\partial}{\partial\xil{l}}
\left(\Jk\frac{\partial\xil{l}}{\partial\xm{m}}\frac{\partial\xil{a}}{\partial \xm{m}}
\frac{\partial(\Bm{m}\U)}{\partial\xil{a}}\right), 
\end{split}
\end{equation}
where the last set of terms on the left-hand side are zero by the GCL conditions \eqref{eq:GCL}. 
Then, a stable semi-discrete form is constructed in the same manner as the split form \eqref{eq:convectiondiffusionsplit} by discretizing the inviscid portion of \eqref{eq:convectiondiffusionchain} and \eqref{eq:convectiondiffusionstrong} using $\Dxil{l}$, $\matJk{\kappa}$, and $\matAlmk{l}{m}{\kappa}$, and by averaging the results. 
The viscous terms result from the discretization of the viscous portion of \cref{eq:convectiondiffusionstrong}. 
This procedure yields
\begin{equation}\label{eq:convectionsplitdisc}
\begin{split}
&\matJk{\kappa}\frac{\mr{d}\uk}{\mr{d} t}+\frac{1}{2}\sum\limits_{l,m=1}^{3}
a_{m}\left\{\Dxil{l}\matAlmk{l}{m}{\kappa}+\matAlmk{l}{m}{\kappa}\Dxil{l}\right\}\uk
\\&-\frac{1}{2}\sum\limits_{l,m=1}^{3}\left\{
a_{m}\diag\left(\uk\right)\Dxil{l}\matAlmk{l}{m}{\kappa}\ones{\kappa}\right\}=\\
&\sum\limits_{l,m,a=1}^{3}\Bm{m}\Dxil{l}\matJk{\kappa}^{-1}\matAlmk{l}{m}{\kappa}\matAlmk{a}{m}{\kappa}\Dxil{a}\uk
+\bm{\mathrm{SAT}}_{\kappa},
\end{split}
\end{equation}
where $\ones{\kappa}$ is a vector of ones of size $\Nl{\kappa}$.

As in the continuous case, the semi-discrete form has a set of discrete GCL conditions
\begin{equation}\label{eq:discGCLconvection}
\sum\limits_{l=1}^{3} \Dxil{l}\matAlmk{l}{m}{\kappa}\ones{\kappa}=\bm{0}, \quad m = 1,2,3,
\end{equation}
that, when satisfied, lead to the following telescoping, provably stable, semi-discrete form
\begin{equation}\label{eq:convectionsplitdisctele}
\begin{split}
&\matJk{\kappa}\frac{\mr{d}\uk}{\mr{d} t}+\frac{1}{2}\sum\limits_{l,m=1}^{3}
a_{m}\left\{\Dxil{l}\matAlmk{l}{m}{\kappa}+\matAlmk{l}{m}{\kappa}\Dxil{l}\right\}\uk=\\
&\sum\limits_{l,m,a=1}^{3}\Bm{m}\Dxil{l}{\matJk{\kappa}}^{-1}\matAlmk{l}{m}{\kappa}\matAlmk{a}{m}{\kappa}\Dxil{a}\uk
+\bm{\mathrm{SAT}}_{\kappa}.
\end{split}
\end{equation}

Herein, we consider only conforming interfaces and optimize the metric terms, $\matAlmk{l}{m}{\kappa}$, as presented in \cite{nolasco_optim_metrics_2019}:
\begin{itemize}
\item The surface metric terms are specified using analytic metrics,
\item Each discrete GCL system \eqref{eq:discGCLconvection} is highly undertermined and is solved using an optimization approach that minimizes the difference between the numerical and analytic volume metrics.
\end{itemize}

In contrast to the metrics for the inviscid terms, the metrics used for the viscous terms need only be, at worst, consistent and design order approximations. 
Herein, we use the analytic metrics for the viscous terms calculation. 

To make the presentation easier and to introduce the general discretization that will later be used for the viscous portion of the compressible Navier--Stokes equations, the inviscid term is lumped into $\fnc{I}^{(\mathrm{I})}$, while the viscous terms are simplified. 
Thus, \cref{eq:convectiondiffusionsplit} reduces to 
\begin{equation}\label{eq:diffusionsplit}
\begin{split}
&\Jk\frac{\partial\U}{\partial t}+\fnc{I}^{(\mathrm{I})}= \sum\limits_{l,a=1}^{3}
\frac{\partial}{\partial\xil{l}} \left(\Chatla{l}{a} \Thetaa{a}\right), \\
&\Chatla{l}{a}\equiv\sum\limits_{m=1}^{3}\Jk\frac{\partial\xil{l}}{\partial\xm{m}}\frac{\partial\xil{a}}{\partial\xm{m}}\Bm{m},\quad 
\Thetaa{a}\equiv\frac{\partial\U}{\partial\xil{a}}.
\end{split}
\end{equation}
A local discontinuous Galerkin (LDG) and interior penalty approach (IP) approach are used (see references~\cite{parsani_entropy_stable_interfaces_2015,carpenter_entropy_stable_staggered_2015,parsani_ssdc_staggered_2016}). 
In the LDG approach, the discretization of the viscous terms in \cref{eq:diffusionsplit} proceeds in two steps. 
First, the gradients $\Thetaa{a}$ are discretized, then the derivatives of the viscous fluxes are discretized. 
Notice that all the metric terms are contained in $\Chatla{l}{a}$, and therefore the critical requirement for stability is to use an SBP operator \cite{parsani_entropy_stable_interfaces_2015,carpenter_entropy_stable_staggered_2015,fernandez_entropy_stable_p_ref_nasa_2019}. 
Plugging everything together, the final discretization reads
\begin{equation}\label{eq:diffusionsplitdisc}
\matJk{\kappa}\frac{\mr{d}\uk}{\mr{d}t}+\bm{I}^{(\mathrm{I})}_{\kappa}
= \sum\limits_{l,a=1}^{3}\Dxil{l}\matChatla{l}{a}_{\kappa}\thetaa{a}^{\kappa}\:+\: \bm{\mathrm{SAT}}^{(\mathrm{I})}_{\kappa} 
+ \bm{\mathrm{SAT}}^{(\mathrm{V})}_{\kappa} ,\quad
\thetaa{a}^{\kappa}=\Dxil{a}\uk 
+ \bm{\mathrm{SAT}}^{\theta}_{\kappa},
\end{equation}
where the inviscid contributions are contained in $\bm{I}^{(\mathrm{I})}_{\kappa}$, while $\bm{\mathrm{SAT}}^{\theta}_{\kappa}$ contains the LDG penalty on the gradient of the variables \cite{parsani_entropy_stable_interfaces_2015}.
The proposed discretization of the viscous terms telescopes the viscous fluxes to the boundary and adds a dissipative term \cite{parsani_entropy_stable_interfaces_2015}. 
Thus, it mimics the continuous energy analysis, and leads to a provably energy stable discretization, provided appropriate boundary SATs are used.

%%%%%%%%%%%%%%%%%%%%%%%%%%%%%%%%%%%%%%%%%%%%%%%%%%%%%%%%%%%%%%%%%%%%%%%%%%%%%%%%
\section{Discretization of the compressible Navier--Stokes equations}\label{sec:compressible_nse}

In this section, the algorithm for the convection-diffusion equation presented in the previous section is applied to the compressible Navier--Stokes.
These equations in Cartesian coordinates read 
\begin{equation}\label{eq:compressible_ns}
\begin{aligned}
&  \frac{\partial\Q}{\partial t}+\sum\limits_{m=1}^{3}\frac{\partial \FxmI{m}}{\partial \xm{m}} 
=  \sum\limits_{m=1}^{3}\frac{\partial \FxmV{m}}{\partial\xm{m}}, 
&& \forall \left(\xm{1},\xm{2},\xm{3}\right)\in\Omega,\quad t\ge 0,\\
&  \Q\left(\xm{1},\xm{2},\xm{3},t\right)=\GB\left(\xm{1},\xm{2},\xm{3},t\right), 
&& \forall \left(\xm{1},\xm{2},\xm{3}\right)\in\Gamma,\quad t\ge 0,\\
&  \Q\left(\xm{1},\xm{2},\xm{3},0\right)=\Gzero\left(\xm{1},\xm{2},\xm{3}\right), 
&& \forall \left(\xm{1},\xm{2},\xm{3}\right)\in\Omega,
\end{aligned}
\end{equation}
where the vectors $\Q$, $\FxmI{m}$ and $\FxmV{m}$ denote the conserved variables, the inviscid fluxes, and the
viscous fluxes, respectively. 
The boundary data, $\GB$, and the initial condition, $\Gzero$, are assumed to be in $L^{2}(\Omega)$,
with the further assumption that $\GB$ will be set to coincide with linear, well-posed boundary conditions, prescribed in such a way that either entropy conservation or entropy stability is achieved.

The vector of conserved variables is given by 
\begin{equation*}
\Q = \left[\rho,\rho\Um{1},\rho\Um{2},\rho\Um{3},\rho\E\right]\Tr,
\end{equation*}
where $\rho$ denotes the density, $\bm{\fnc{U}} = \left[\Um{1},\Um{2},\Um{3}\right]\Tr$ is the velocity 
vector, and $\E$ is the specific total energy. 
The inviscid fluxes are given as
\begin{equation}\label{eq:ns_iflux}
\begin{split}
\FxmI{m} = 
& \left[\rho\Um{m},\rho\Um{m}\Um{1}+\delta_{m,1}\fnc{P},\rho\Um{m}\Um{2}+\delta_{m,2}\fnc{P},\right.
\left.\rho\Um{m}\Um{3}+\delta_{m,3}\fnc{P},\rho\Um{m}\fnc{H}\right]\Tr,
\end{split}
\end{equation}
where $\fnc{P}$ is the pressure, $\fnc{H}$ is the specific total enthalpy and $\delta_{i,j}$ is the 
Kronecker delta.

The required constituent relations are
\begin{equation*}
\fnc{H} = c_{\fnc{P}}\fnc{T}+\frac{1}{2}\bm{\fnc{U}}\Tr\bm{\fnc{U}},\quad \fnc{P} = \rho R \fnc{T},\quad R = \frac{R_{u}}{M_{w}},
\end{equation*}
where $\fnc{T}$ is the temperature, $R_{u}$ is the universal gas constant, $M_{w}$ is the molecular weight of the gas, and $c_{\fnc{P}}$ is the specific heat capacity at constant pressure. 
Finally, the specific thermodynamic entropy is given as 
\begin{equation*}
s=\frac{R}{\gamma-1}\log\left(\frac{\fnc{T}}{\fnc{T}_{\infty}}\right)-R\log\left(\frac{\rho}{\rho_{\infty}}\right),\quad \gamma=\frac{c_{p}}{c_{p}-R},
\end{equation*}
where $\fnc{T}_{\infty}$ and $\rho_{\infty}$ are the reference temperature and density, respectively
(the stipulated convention has been broken here and $s$ has been used rather than $\fnc{S}$ for reasons that will be clear next). 

The viscous fluxes $\FxmV{m}$ are given by
\begin{equation}\label{eq:Fv}
\FxmV{m}=\left[0,\tau_{1,m},\tau_{2,m},\tau_{3,m}, \sum\limits_{i=1}^{3}\tau_{i,m}\fnc{U}_{i}-\kappa\frac{\partial \fnc{T}}{\partial\xm{m}}\right]\Tr,
\end{equation} 
while the viscous stresses are defined as
\begin{equation}\label{eq:tau}
\tau_{i,j} = \mu\left(\frac{\partial\fnc{U}_{i}}{\partial x_{j}}+\frac{\partial\fnc{U}_{j}}{\partial x_{i}}
-\delta_{i,j}\frac{2}{3}\sum\limits_{n=1}^{3}\frac{\partial\fnc{U}_{n}}{\partial x_{n}}\right),
\end{equation}
where $\mu(\fnc{T})$ is the dynamic viscosity and $\kappa(\fnc{T})$ is the thermal conductivity.

The compressible Navier--Stokes equations given in \cref{eq:NSCCS1} have a convex extension, that when integrated over the physical domain, $\Omega$, depends only on the boundary data and negative semi-definite dissipation terms.
This convex extension depends on an entropy function, $\fnc{S}$, that is constructed from the thermodynamic entropy as
\begin{equation*}
\fnc{S}=-\rho s,
\end{equation*}
and provides a mechanism for proving stability in the $L^{2}$ norm. 
The entropy variables $\bfnc{W}$ are an alternative variable set related to the conservative variables via a one-to-one mapping.  
They are defined in terms of the entropy function $\fnc{S}$ by the relation $\bfnc{W}\Tr=\partial\fnc{S}/\partial\bfnc{Q}$ and they are extensively used in the entropy stability proofs of the algorithms used herein; see for instance \cite{carpenter_ssdc_2014,parsani_ssdc_staggered_2016,friedrich_hp_entropy_stability_2018,fernandez_entropy_stable_hp_ref_snpdea_2019}.
In addition, they simultaneously symmetrize the inviscid and the viscous flux Jacobians in all three spatial directions.
Further details on continuous entropy analysis are available elsewhere \cite{dafermos_book_2010,parsani_entropy_stability_solid_wall_2015,carpenter_entropy_stable_staggered_2015}.

The entropy stability for the viscous terms in the compressible Navier--Stokes equations \eqref{eq:NSCCS1} is readily demonstrated by exploiting the symmetrizing properties of the entropy variables.
Thus, we recast the viscous fluxes in terms of the entropy variables
\begin{equation}\label{eq:Fxment}
\FxmV{m} =\sum\limits_{j=1}^{3}\Cij{m}{j}\frac{\partial\bfnc{W}}{\partial x_{j}},
\end{equation}
with the flux Jacobian matrices satisfying $\Cij{m}{j} \:=\: {(\Cij{j}{m})}\Tr$.

Furthermore, in order to apply the algorithm outlined for the convection-diffusion case to the compressible Navier--Stokes equations, we have to recast system \eqref{eq:compressible_ns} in a skew-symmetric form with respect to the metric terms. 
This procedure results in
\begin{equation}\label{eq:NSCCS1}
\begin{split}
&\Jk\frac{\partial\bfnc{Q}}{\partial t}+\sum\limits_{l,m=1}^{3}\frac{1}{2}
\frac{\partial }{\partial \xil{l}}\left(\Jdxildxm{l}{m}\FxmI{m}\right)
+\frac{1}{2}\Jdxildxm{l}{m}\frac{\partial \FxmI{m}}{\partial \xil{l}}
=\sum\limits_{l,m=1}^{3}\frac{\partial}{\partial\xil{l}}\left(\Jdxildxm{l}{m}\FxmV{m}\right)
\end{split},
\end{equation}
where the GCL relations given in \cref{eq:GCL} are used to obtain \cref{eq:NSCCS} from the divergence form \eqref{eq:compressible_ns}.
Substituting \cref{eq:Fxment} into \cref{eq:NSCCS1}, we arrive at the system of equations
\begin{equation}\label{eq:NSCCS}
\begin{split}
&\Jk\frac{\partial\bfnc{Q}}{\partial t}+\sum\limits_{l,m=1}^{3}
\frac{1}{2}\frac{\partial }{\partial \xil{l}}\left(\Jdxildxm{l}{m}\Fxm{l}\right)
+\frac{1}{2}\Jdxildxm{l}{m}\frac{\partial \Fxm{m}}{\partial \xil{l}}
=\sum\limits_{l,a=1}^{3}\frac{\partial}{\partial\xil{l}}\left(\Chatij{l}{a}\frac{\partial\bfnc{W}}{\partial \xil{a}}\right),
\end{split}
\end{equation}
where
\begin{equation}\label{eq:Chatij}
\Chatij{l}{a}=\Jdxildxm{l}{m}\sum\limits_{m,j=1}^{3}\Cij{m}{j}\frac{\partial\xil{a}}{\partial x_{j}}.
\end{equation}
The symmetric properties of the viscous flux Jacobians are preserved by the rotation into curvilinear
coordinates, \ie $\Chatij{l}{a} \:=\: {(\Chatij{a}{l})}\Tr$.
We remark that this form of the equations, \ie skew-symmetric form plus the quadratic form of the viscous terms, is necessary for the construction of the entropy stable schemes used in this work.
For further details on the derivation of these viscous coefficient matrices see \cite{fisher_phd_2012,parsani_entropy_stability_solid_wall_2015}.

The discretization of the compressible Euler equations, i.e., the inviscid part of \eqref{eq:NSCCS}, is given by 
\begin{equation}\label{eq:discEL}
\begin{split}
&\matJk{\kappa}\frac{\mr{d}\qk{\kappa}}{\mr{d}t}+\frac{1}{2}\sum\limits_{l,m=1}^{3}\left(\Dxil{l}\matAlmk{l}{m}{\kappa}+\matAlmk{m}{l}{\kappa}\Dxil{l}\right)
\circ\matFxm{m}{\qk{\kappa}}{\qk{\kappa}}\ones{\kappa}=
\bm{\mathrm{SAT}}^{(\mathrm{I})}_{\kappa}, 
\end{split}
\end{equation}
where the symbol $\circ$ indicates the Hadamard product, and $\matFxm{m}{\cdot}{\cdot}$ is a two argument matrix flux function which, for an entropy conservative schemes, is constructed from a two-point entropy
conservative flux function, $\FxmI{m}(\bfnc{Q}_i,\bfnc{Q}_{j})$, (see, for instance, \cite{fernandez_entropy_stable_p_ref_nasa_2019}).

Herein, to construct the entropy conservative discretization, we use the two-point entropy conservative flux by \citet{chandrashekar2013kinetic} which reads:
\begin{equation}
\FxmI{m}(\bfnc{Q}_i,\bfnc{Q}_{j}) =
\begin{bmatrix}
\widehat{\rho}\overline{\mathcal{U}}_m \\
\widehat{\rho}\overline{\mathcal{U}}_m \overline{\mathcal{U}}_1 + \delta_{m,1}\widetilde{\mathcal{P}} \\
\widehat{\rho}\overline{\mathcal{U}}_m \overline{\mathcal{U}}_2 + \delta_{m,2}\widetilde{\mathcal{P}} \\
\widehat{\rho}\overline{\mathcal{U}}_m \overline{\mathcal{U}}_3 + \delta_{m,3}\widetilde{\mathcal{P}} \\
\widehat{\rho}\overline{\mathcal{U}}_m \left( \frac{1}{2(\gamma-1)\widehat{\beta}} -\frac{1}{2} \overline{ \left| \bm{\mathcal{U}} \right|^2 } \right) + \overline{ \bm{ \mathcal{U} } } \cdot \bm{\mathcal{M}}
\end{bmatrix},
\end{equation}
where, for the generic quantity $\phi$,
\[
\begin{aligned}
\widehat{\phi} = \frac{\phi_i-\phi_{j}}{\log{\phi_i}-\log{\phi_{j}}} && \mathrm{,} &&
\overline{\phi} = \frac{\phi_i+\phi_{j}}{2}
\end{aligned}
\]
are the logarithmic average and the arithmetic average, respectively,
\[
\begin{aligned}
\widetilde{\mathcal{P}}=R\overline{\rho}\frac{\mathcal{T}_i \mathcal{T}_{j}}{\overline{\mathcal{T}}} && \mathrm{,} && \beta=\frac{1}{2R\mathcal{T}},
\end{aligned}
\]
and the term $\bm{\mathcal{M}}$ corresponds to a vector of the three momentum components of this two-point flux. We will refer to our entropy stable discretization as ES-C.

The Hadamard formalism is capable of compactly representing various split forms, and more importantly, extends to nonlinear equations for which a canonical split form is inappropriate.
In this work, we also use this formalism to compute the divergence of the inviscid fluxes using the cheaper two-point flux of \citet{kennedy2008reduced}:
\begin{equation}
\FxmI{m}(\bfnc{Q}_i,\bfnc{Q}_{j}) =
\begin{bmatrix}
\overline{\rho}\overline{\mathcal{U}}_m \\
\overline{\rho}\overline{\mathcal{U}}_m \overline{\mathcal{U}}_1 + \delta_{m,1}\overline{\mathcal{P}} \\
\overline{\rho}\overline{\mathcal{U}}_m \overline{\mathcal{U}}_2 + \delta_{m,2}\overline{\mathcal{P}} \\
\overline{\rho}\overline{\mathcal{U}}_m \overline{\mathcal{U}}_3 + \delta_{m,3}\overline{\mathcal{P}} \\
\overline{\rho}\overline{\mathcal{U}}_m \overline{\mathcal{E}} + \overline{\mathcal{P}} \  \overline{\mathcal{U}}_m
\end{bmatrix}.
\end{equation}
We will refer to the discretization that uses the \citet{kennedy2008reduced} flux as SF-KG.

Finally, we also compare the above two discretizations with a standard discontinuous collocation type discretization (DC).
In this case the divergence of the inviscid fluxes is computed by applying the SBP differentiation matrix to the inviscid flux vector whose components are the direct evaluation of \eqref{eq:ns_iflux} at each node:
\begin{equation}
\begin{split}
\FxmI{m}(\bfnc{Q}_i) = &\left[\rho \mathcal{U}_m,\rho \mathcal{U}_m \mathcal{U}_1 + \delta_{m,1} \mathcal{P}, \rho u_m \mathcal{U}_2 +\delta_{m,2} p, \rho \mathcal{U}_m \mathcal{U}_3 + \delta_{m,3} \mathcal{P} , \rho \mathcal{P}_m \mathcal{H} \right]\Tr.
\end{split}
\end{equation}
The above flux is the cheapest flux among the fluxes considered herein.

\begin{remark}
It is worth noting the differences in the cost of these discretizations. 
The computing of logarithms contributes significantly to the cost of the inviscid flux evaluations of the entropy stable discretization.
Furthermore, the standard DC methods require fewer inviscid flux evaluations than the other two discretizations with $\mathcal{O}(N^3)$ evaluations instead of the $\mathcal{O}(N^4)$ evaluations in the three-dimensional case.
\Cref{sec:numerical_results} presents results exploring the differences in the cost of these discretizations.
\end{remark}

Next, recasting the viscous fluxes in terms of entropy variables as shown in \cref{eq:Fxment} yields the following form for the discretization of the divergence of the viscous fluxes
\begin{equation}\label{eq:MacroForm}
\sum\limits_{l,a=1}^{3}\frac{\partial}{\partial\xil{l}}\left(\Chatla{l}{a}\frac{\partial\bfnc{W}}{\partial \xil{a}}\right)
\approx  \sum\limits_{l,a=1}^{3} \Dxil{l}\matChatla{l}{a}\thetaa{a}^{\kappa}, \qquad \thetaa{a}^{\kappa}=\Dxil{a}\wk{\kappa}. 
\end{equation}
Note that \cref{eq:MacroForm} is precisely the symmetric generalization of the convection-diffusion operator to a viscous system.  

The discretization on the $\kappa\Th$ element reads
\begin{equation}\label{eq:NSL}
\begin{split}
&\matJk{\kappa}\frac{\mr{d}\qk{\kappa}}{\mr{d}t}+\bm{I}^{(\mathrm{E})}_{\kappa}=\sum\limits_{l,a=1}^{3}
\Dxil{a}\matChatla{l}{a}\thetaa{a}^{\kappa} + \bm{\mathrm{SAT}}^{(\mathrm{I})}_{\kappa} + \bm{\mathrm{SAT}}^{(\mathrm{V})}_{\kappa} + \IP^{\kappa},\quad
\thetaa{a}^{\kappa}=\Dxil{a}\wk{\kappa} + \bm{\mathrm{SAT}}^{\theta}_{\kappa},
\end{split}
\end{equation}
where $\bm{I}^{(\mathrm{E})}_{\kappa}$ represents the discretization of the divergence of the inviscid fluxes and the interior penalty term, $\IP^{\kappa}$, adds interface dissipation \cite{parsani_entropy_stable_interfaces_2015}.
This term is a design-order zero interface dissipation term that is constructed to damp neutrally stable ``odd-even'' eigenmodes that arise from the LDG viscous operator. 
Scheme \eqref{eq:NSL} telescopes to the boundaries where appropriate SATs need to be imposed to obtain a stability statement \cite{parsani_wall_bc_entropy_2015,parsani_ssdc_staggered_2016,svard_entropy_stable_solid_wall_2018}.

%%%%%%%%%%%%%%%%%%%%%%%%%%%%%%%%%%%%%%%%%%%%%%%%%%%%%%%%%%%%%%%%%%%%%%%%%%%%%%%%
\section{Numerical results}\label{sec:numerical_results}

The curvilinear, unstructured grid, CFD framework used in this article has been developed at the Extreme Computing Research Center (ECRC) at KAUST.
The conforming numerical solver is based on the algorithms proposed in \cite{carpenter_ssdc_2014,parsani_entropy_stability_solid_wall_2015,parsani_entropy_stable_interfaces_2015}.
It is built on top of the Portable and Extensible Toolkit for Scientific computing (PETSc) \cite{petsc-user-ref}, its mesh topology abstraction (DMPLEX) \cite{KnepleyKarpeev09} and scalable ordinary differential equation (ODE)/differential algebraic equations (DAE) solver library \cite{abhyankar2018petsc}.
The fifth order explicit Runge--Kutta scheme by \citet{DORMAND198019} with an adaptive time-stepping based on signal processing \cite{Soderlind2003,Soderlind2006}) is used to integrate the numerical solution in time. 
The tolerances of the time integrator used in all test cases are small enough to render the time error negligible.

%%%%%%%%%%%%%%%%%%%%%%%%%%%%%%%%%%%%%%%%%%%%%%%%%%%%%%%%%%%%%%%%%%%%%%%%%%%%%%%%
\subsection{Error vs. cost}

Convergence studies for 3D flows are conducted, and the serial time to solution is measured.
The convergence studies are done on structured cubic domains, using sequences of nested grids.
Finally, the $L^2$ norm of error is calculated using the matrix norm ($\M$-norm) associated with the SBP-SAT scheme.

\subsubsection{Inviscid flow: 3D isentropic vortex}
Although the focus of this paper is on the compressible Navier--Stokes equations, in this section, we study the propagation of an inviscid 3D isentropic vortex.
This first test allows us to verify and characterize the cost and accuracy of the inviscid component of our algorithm, a key element in any compressible or incompressible Navier--Stokes solver \cite{carpenter_ssdc_2014,parsani_ssdc_staggered_2016,fernandez_entropy_stable_p_euler_2019}.

The simulation of an isentropic vortex is a widely used benchmark problem \cite{shu1998essentially} because it has an analytical solution. 
In particular, for the stationary case, the exact solution is given by
\begin{equation} \label{eq:iv_sol}
\begin{cases}
\rho              &= \T^\frac{1}{\gamma-1},
\\ \U_{\mathrm{t}}&= \frac{r\beta}{2\pi}\exp{\left(\frac{1-r^2}{2}\right)}, 
\\ \T             &= \T_{\infty}-\frac{(\gamma-1)\Ma_\infty^2\beta^2}{8\pi^2}\exp{\left(1-r^2\right)},
\end{cases}
\end{equation}
where $r$ is the distance from the axis of the vortex, and $\U_{\mathrm{t}}$ is the tangential velocity.
The moving vortex solution is obtained by applying to \eqref{eq:iv_sol} a uniform translation in the direction of the velocity vector field.

Herein, the simulation domain is a cube $\Omega=[0,10]^3$ where the vortex rotates around the axis $(1, 1, 1)\Tr$, a direction not aligned with the grid.
A constant velocity field $\Um{m}=\Um{m}^{\infty}$ is imposed and the vortex is simulated for a short amount of time.
The parameters for this test are $\gamma=1.4$, $\Ma_\infty=0.5$, $\beta=5$ and $\T_\infty=1$.

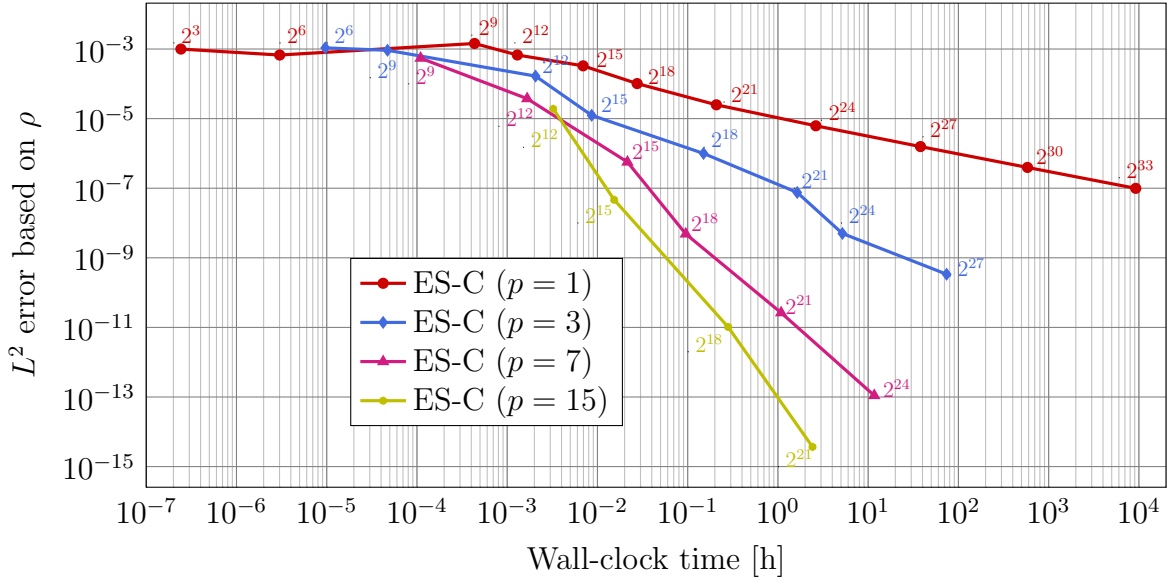
\begin{figure}[!ht]
\centering
\begin{tikzpicture}
\definecolor{color0}{rgb}{0.75,0.75,0}
\definecolor{color1}{rgb}{0.75,0,0.75}
\begin{axis}[height=8cm,width=15.0cm,legend cell align={left},legend entries={ES-C ($p=1$),ES-C ($p=3$),ES-C ($p=7$),ES-C ($p=15$)},legend style={at={(0.20,0.3)},anchor=west},legend columns=4,log basis x={10},log basis y={10},tick align=outside,tick pos=left,xlabel={Wall-clock time [h]},ylabel={$L^2$ error based on $\rho$},xmajorgrids,xminorgrids,xmode=log,ymode=log,grid=both,major grid style={black!50},major tick length=0pt,minor tick length=0pt,legend columns=1,xmin=0.0000001,xmax=20000]
\addplot [very thick, mycolor1, mark=*, mark size=1.5, mark options={solid}]table[x index=0,y index=1]{data/IV3D/order1.dat};
\addplot [very thick, mycolor2, mark=diamond*, mark size=1.5, mark options={solid}]table[x index=0,y index=1]{data/IV3D/order3.dat};
\addplot [very thick, mycolor3, mark=triangle*, mark size=1.5, mark options={solid}]table[x index=0,y index=1]{data/IV3D/order7.dat};
\addplot [very thick, color0, mark=asterisk, mark size=1.5, mark options={solid}]table[x index=0,y index=1]{data/IV3D/order15.dat};

\draw[] (axis cs:2.0e-07,1.5e-03) -- (axis cs:2.0e-07,1.5e-03); \node at (axis cs:2.0e-07,1.5e-03)[scale=0.70,anchor=base west,text=mycolor1,rotate=0.0]{$2^{3}$};
\draw[] (axis cs:3.0e-06,1.5e-03) -- (axis cs:3.0e-06,1.5e-03); \node at (axis cs:3.0e-06,1.5e-03)[scale=0.70,anchor=base west,text=mycolor1,rotate=0.0]{$2^{6}$};
\draw[] (axis cs:4.0e-04,2.5e-03) -- (axis cs:4.0e-04,2.5e-03); \node at (axis cs:4.0e-04,2.5e-03)[scale=0.70,anchor=base west,text=mycolor1,rotate=0.0]{$2^{9}$};
\draw[] (axis cs:1.2e-03,1.5e-03) -- (axis cs:1.2e-03,1.5e-03); \node at (axis cs:1.2e-03,1.5e-03)[scale=0.70,anchor=base west,text=mycolor1,rotate=0.0]{$2^{12}$};
\draw[] (axis cs:8.0e-03,3.0e-04) -- (axis cs:8.0e-03,3.0e-04); \node at (axis cs:8.0e-03,3.0e-04)[scale=0.70,anchor=base west,text=mycolor1,rotate=0.0]{$2^{15}$};
\draw[] (axis cs:3.0e-02,1.0e-04) -- (axis cs:3.0e-02,1.0e-04); \node at (axis cs:3.0e-02,1.0e-04)[scale=0.70,anchor=base west,text=mycolor1,rotate=0.0]{$2^{18}$};
\draw[] (axis cs:2.5e-01,3.0e-05) -- (axis cs:2.5e-01,3.0e-05); \node at (axis cs:2.5e-01,3.0e-05)[scale=0.70,anchor=base west,text=mycolor1,rotate=0.0]{$2^{21}$};
\draw[] (axis cs:3.0e-00,9.0e-06) -- (axis cs:3.0e-00,9.0e-06); \node at (axis cs:3.0e-00,9.0e-06)[scale=0.70,anchor=base west,text=mycolor1,rotate=0.0]{$2^{24}$};
\draw[] (axis cs:4.0e+01,3.0e-06) -- (axis cs:4.0e+01,3.0e-06); \node at (axis cs:4.0e+01,3.0e-06)[scale=0.70,anchor=base west,text=mycolor1,rotate=0.0]{$2^{27}$};
\draw[] (axis cs:6.0e+02,5.0e-07) -- (axis cs:6.0e+02,5.0e-07); \node at (axis cs:6.0e+02,5.0e-07)[scale=0.70,anchor=base west,text=mycolor1,rotate=0.0]{$2^{30}$};
\draw[] (axis cs:6.0e+03,2.0e-07) -- (axis cs:6.0e+03,2.0e-07); \node at (axis cs:6.0e+03,2.0e-07)[scale=0.70,anchor=base west,text=mycolor1,rotate=0.0]{$2^{33}$};

\draw[] (axis cs:1.0e-05,1.5e-03) -- (axis cs:1.0e-05,1.5e-03); \node at (axis cs:1.0e-05,1.5e-03)[scale=0.70,anchor=base west,text=mycolor2,rotate=0.0]{$2^{6}$};
\draw[] (axis cs:3.0e-05,1.5e-04) -- (axis cs:3.0e-05,1.5e-04); \node at (axis cs:3.0e-05,1.5e-04)[scale=0.70,anchor=base west,text=mycolor2,rotate=0.0]{$2^{9}$};
\draw[] (axis cs:2.0e-03,2.0e-04) -- (axis cs:2.0e-03,2.0e-04); \node at (axis cs:2.0e-03,2.0e-04)[scale=0.70,anchor=base west,text=mycolor2,rotate=0.0]{$2^{12}$};
\draw[] (axis cs:9.0e-03,1.5e-05) -- (axis cs:9.0e-03,1.5e-05); \node at (axis cs:9.0e-03,1.5e-05)[scale=0.70,anchor=base west,text=mycolor2,rotate=0.0]{$2^{15}$};
\draw[] (axis cs:1.5e-01,1.5e-06) -- (axis cs:1.5e-01,1.5e-06); \node at (axis cs:1.5e-01,1.5e-06)[scale=0.70,anchor=base west,text=mycolor2,rotate=0.0]{$2^{18}$};
\draw[] (axis cs:1.5e-00,1.0e-07) -- (axis cs:1.5e-00,1.0e-07); \node at (axis cs:1.5e-00,1.0e-07)[scale=0.70,anchor=base west,text=mycolor2,rotate=0.0]{$2^{21}$};
\draw[] (axis cs:5.0e-00,1.0e-08) -- (axis cs:5.0e-00,1.0e-08); \node at (axis cs:5.0e-00,1.0e-08)[scale=0.70,anchor=base west,text=mycolor2,rotate=0.0]{$2^{24}$};
\draw[] (axis cs:8.0e+01,3.0e-10) -- (axis cs:8.0e+01,3.0e-10); \node at (axis cs:8.0e+01,3.0e-10)[scale=0.70,anchor=base west,text=mycolor2,rotate=0.0]{$2^{27}$};

\draw[] (axis cs:8.0e-05,1.0e-04) -- (axis cs:8.0e-05,1.0e-04); \node at (axis cs:8.0e-05,1.0e-04)[scale=0.70,anchor=base west,text=mycolor3,rotate=0.0]{$2^{9}$};
\draw[] (axis cs:8.0e-04,6.0e-06) -- (axis cs:8.0e-04,6.0e-06); \node at (axis cs:8.0e-04,6.0e-06)[scale=0.70,anchor=base west,text=mycolor3,rotate=0.0]{$2^{12}$};
\draw[] (axis cs:2.0e-02,7.0e-07) -- (axis cs:2.0e-02,7.0e-07); \node at (axis cs:2.0e-02,7.0e-07)[scale=0.70,anchor=base west,text=mycolor3,rotate=0.0]{$2^{15}$};
\draw[] (axis cs:9.0e-02,6.0e-09) -- (axis cs:9.0e-02,6.0e-09); \node at (axis cs:9.0e-02,6.0e-09)[scale=0.70,anchor=base west,text=mycolor3,rotate=0.0]{$2^{18}$};
\draw[] (axis cs:1.0e-00,2.5e-11) -- (axis cs:1.0e-00,2.5e-11); \node at (axis cs:1.0e-00,2.5e-11)[scale=0.70,anchor=base west,text=mycolor3,rotate=0.0]{$2^{21}$};
\draw[] (axis cs:1.2e+01,1.0e-13) -- (axis cs:1.2e+01,1.0e-13); \node at (axis cs:1.2e+01,1.0e-13)[scale=0.70,anchor=base west,text=mycolor3,rotate=0.0]{$2^{24}$};

\draw[] (axis cs:1.5e-03,1.5e-06) -- (axis cs:1.5e-03,1.5e-06); \node at (axis cs:1.5e-03,1.5e-06)[scale=0.70,anchor=base west,text=color0,rotate=0.0]{$2^{12}$};
\draw[] (axis cs:6.0e-03,1.0e-08) -- (axis cs:6.0e-03,1.0e-08); \node at (axis cs:6.0e-03,1.0e-08)[scale=0.70,anchor=base west,text=color0,rotate=0.0]{$2^{15}$};
\draw[] (axis cs:1.0e-01,2.0e-12) -- (axis cs:1.0e-01,2.0e-12); \node at (axis cs:1.0e-01,2.0e-12)[scale=0.70,anchor=base west,text=color0,rotate=0.0]{$2^{18}$};
\draw[] (axis cs:1.0e-00,1.0e-15) -- (axis cs:1.0e-00,1.0e-15); \node at (axis cs:1.0e-00,1.0e-15)[scale=0.70,anchor=base west,text=color0,rotate=0.0]{$2^{21}$};

\end{axis}
\end{tikzpicture}
\caption{Error vs. cost for the 3D isentropic vortex. The numbers represent the number of degrees of freedom of each simulation.}
\label{fig_iv3D_costVSerror}
\end{figure}

\Cref{fig_iv3D_costVSerror} shows the $L^2$ norm of the error of the density field computed at the final time $t_f = 2.5$ against the serial wall-clock time for the ES-C discretization.
We remark that even for large errors, e.g. $10^{-4}$, the smallest time to solution is always achieved with high order discretizations (namely $p=7$ or even $p=15$).

\subsubsection{Viscous flow: 3D manufactured solution}

\begin{figure}[!ht]
\centering
\setlength{\mywidthsubfig}{0.32\textwidth}
\setlength{\mywidthgraph}{0.45\textwidth}
\begin{subfigure}{\mywidthsubfig}\centering
\includegraphics[width=1.0\columnwidth,trim = 0 0 -5 0,clip]{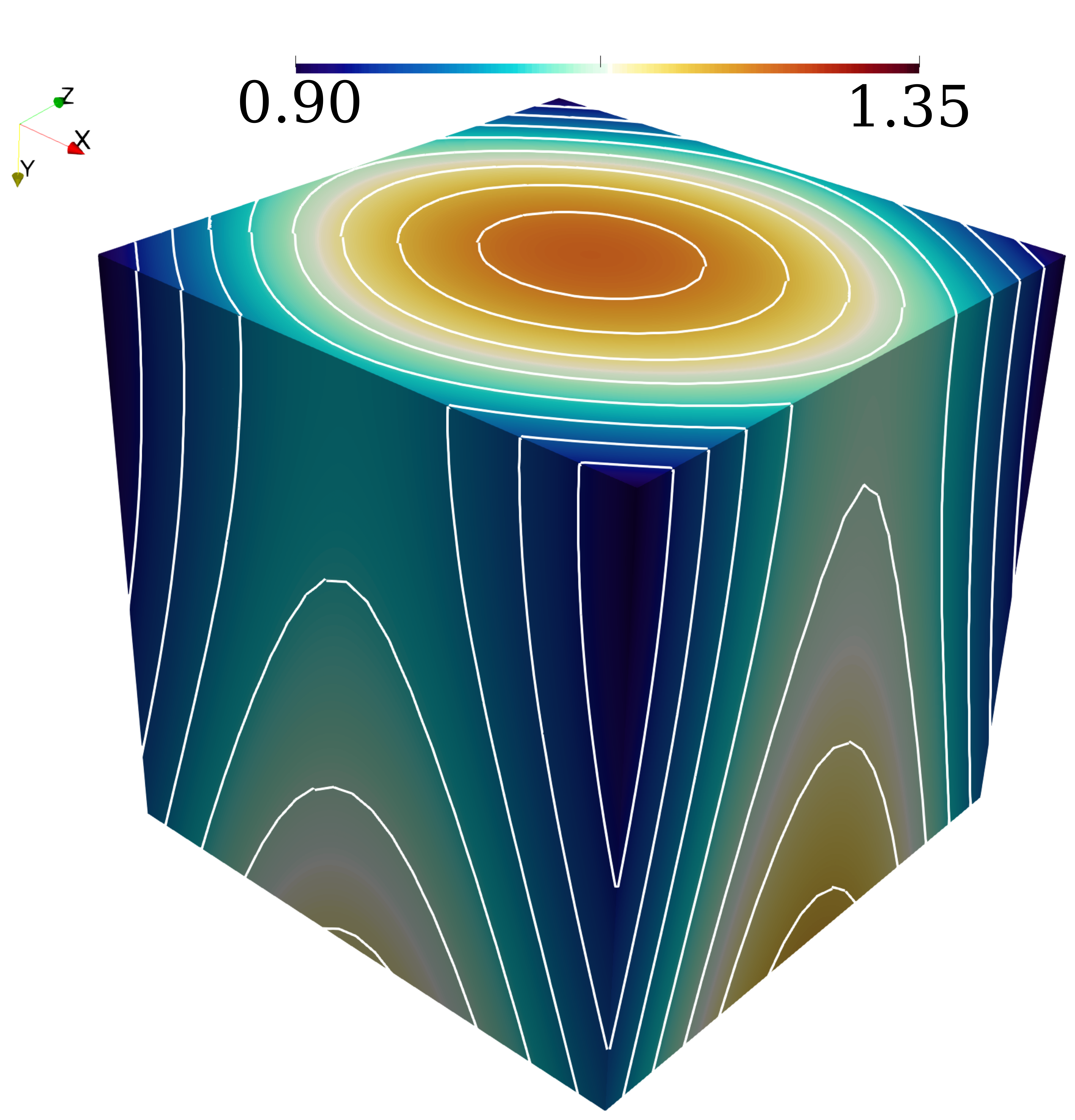}
\sublabelsty{a}\phantomsubcaption
\label{subfig-rho}
\end{subfigure}
\begin{subfigure}{\mywidthsubfig}\centering
\includegraphics[width=1.0\columnwidth,trim = 0 0 -5 0,clip]{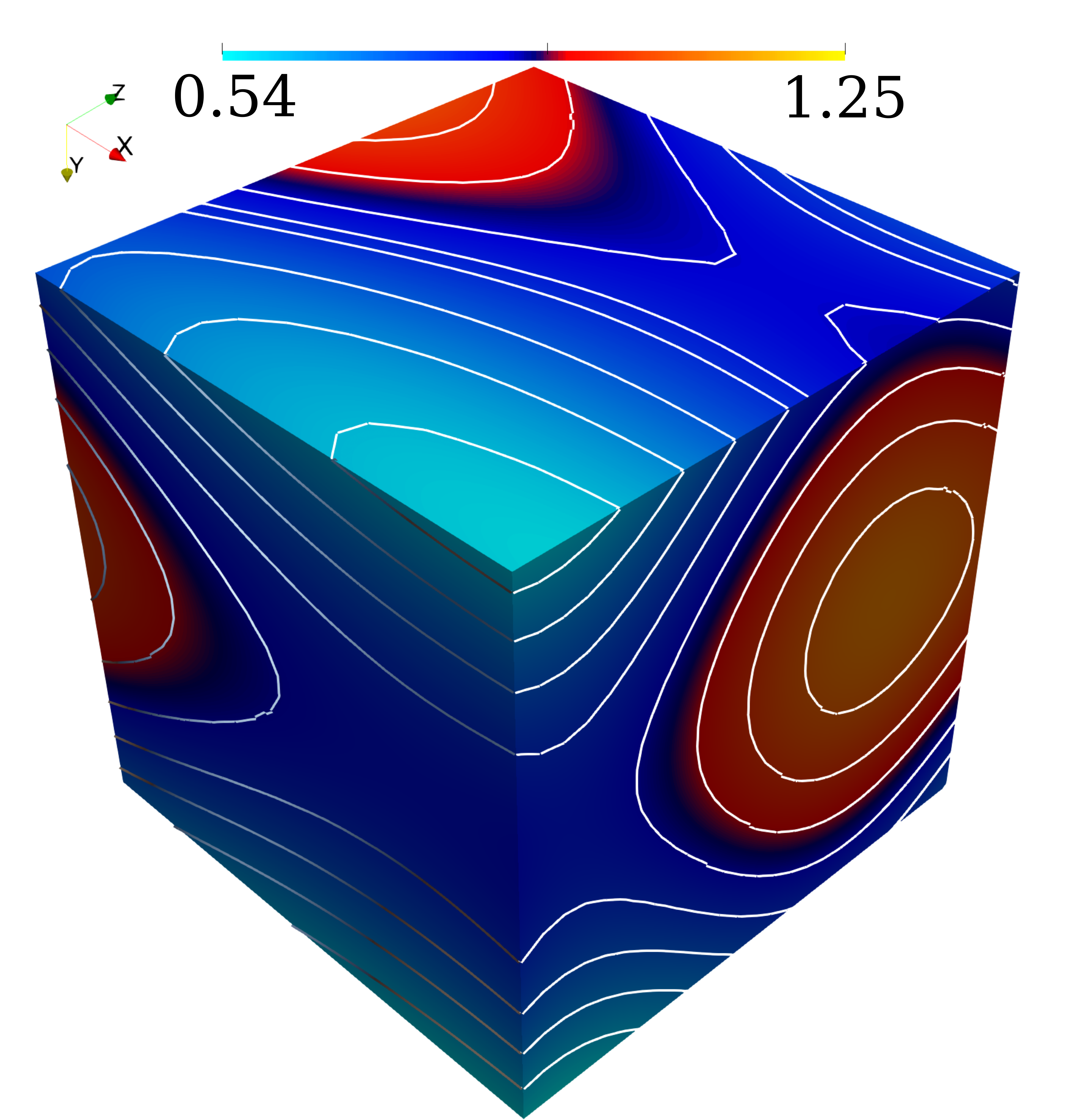}
\sublabelsty{b}\phantomsubcaption
\label{subfig-T}
\end{subfigure}
\begin{subfigure}{\mywidthsubfig}\centering
\includegraphics[width=1.0\columnwidth,trim = 0 0 -5 0,clip]{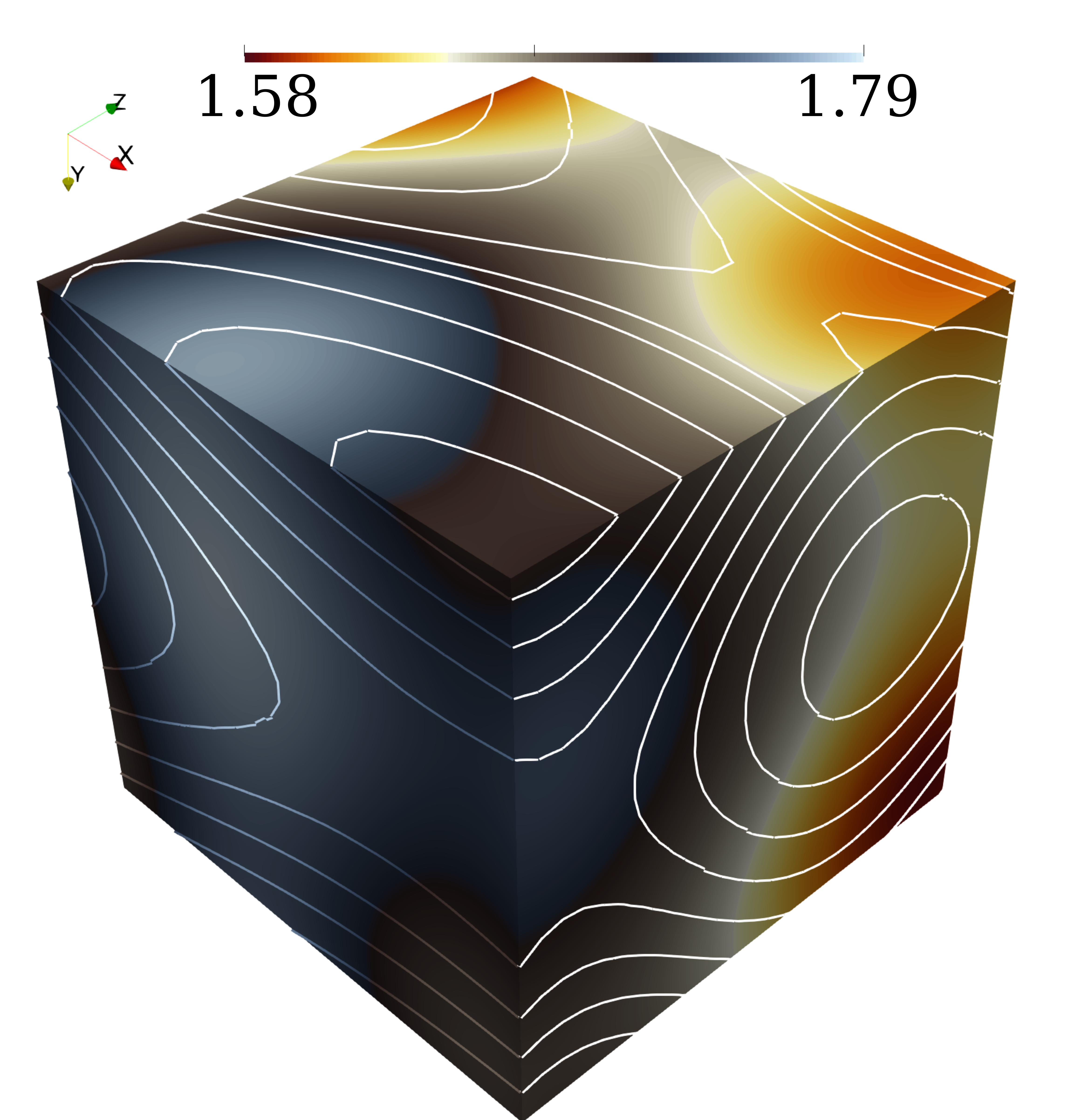}
\sublabelsty{c}\phantomsubcaption
\label{subfig-vel}
\end{subfigure}
\caption{3D supersonic manufactured solution: (\subref{subfig-rho}) density, (\subref{subfig-T}) temperature and (\subref{subfig-vel}) velocity magnitude.}
\label{fig-mms-fields}
\end{figure}

To extend the error vs. cost discussion to the realm of viscous flows, we use the method of manufactured solutions (MMS).
The MMS is a powerful technique for code verification widely used in the scientific community \cite{10.1115/1.1436090,roy2005} where, in principle, the manufactured solution doesn't need to satisfy the PDE under investigation.
Here, a proposed smooth solution is inserted into the original PDE so that all derivatives can be calculated analytically.
The result is then simplified, and the residual obtained is used as a source term that, once added to the original equation, creates a modified problem for which the analytical solution is known.
This solution allows for the comparison between the exact and numerical solutions.
In this section, this technique is used to evaluate the accuracy versus wall-clock time for different polynomial degrees using both the conventional DC and ES-C discretizations for a supersonic compressible viscous flow.
To do so, we consider a cubic domain $\Omega=[0,1]^3$ with the following manufactured solution
\begin{equation}\label{eq:mms}
\begin{cases}
\rho&=\rho^0+\rho^{x_1}\sin\left(\frac{\alpha_{\rho}^{x_1}\pi x_1}{L}\right)+\rho^{x_2}\cos\left(\frac{\alpha_{\rho}^{x_2}\pi x_2}{L}\right)+\rho^{x_3}\sin\left(\frac{\alpha_{\rho}^{x_3}\pi x_3}{L}\right),
\\
\mathcal{U}_1&=\mathcal{U}_1^0+\mathcal{U}_1^{x_1}\sin\left(\frac{\alpha_{1}^{x_1}\pi x_1}{L}\right)+\mathcal{U}_1^{x_2}\cos\left(\frac{\alpha_{1}^{x_2}\pi x_2}{L}\right)+\mathcal{U}_3^{x_3}\cos\left(\frac{\alpha_{1}^{x_3}\pi x_3}{L}\right),
\\
\mathcal{U}_2&=\mathcal{U}_2^0+\mathcal{U}_2^{x_1}\cos\left(\frac{\alpha_{2}^{x_1}\pi x_1}{L}\right)+\mathcal{U}_2^{x_2}\sin\left(\frac{\alpha_{2}^{x_2}\pi x_2}{L}\right)+\mathcal{U}_2^{x_3}\sin\left(\frac{\alpha_{2}^{x_3}\pi x_3}{L}\right),
\\
\mathcal{U}_3&=\mathcal{U}_3^0+\mathcal{U}_3^{x_1}\sin\left(\frac{\alpha_{3}^{x_1}\pi x_1}{L}\right)+\mathcal{U}_3^{x_2}\sin\left(\frac{\alpha_{3}^{x_2}\pi x_2}{L}\right)+\mathcal{U}_3^{x_3}\cos\left(\frac{\alpha_{3}^{x_3}\pi x_3}{L}\right),
\\
\mathcal{P}&=\mathcal{P}^0+\mathcal{P}^{x_1}\sin\left(\frac{\alpha_{\mathcal{P}}^{x_1}\pi x_1}{L}\right)+\mathcal{P}^{x_2}\sin\left(\frac{\alpha_{\mathcal{P}}^{x_2}\pi x_2}{L}\right)+\mathcal{P}^{x_3}\cos\left(\frac{\alpha_{p}^{x_3}\pi x_3}{L}\right).
\end{cases}
\end{equation}
Expressions \eqref{eq:mms} are similar to those proposed by \citet{roy2005} and \citet{katz2011}.  
The constants appearing in this manufactured solution are chosen to give supersonic flow in the three spatial directions (cf. \Cref{tab_3dmms_ns}).
Therefore, the exact Dirichlet values for all primitive variables are specified on both inflow and outflow boundaries.
The value of the parameters are: $L=1$, $\gamma=1.4$, $\Rey=4\times10^6$, and $\Ma\approx 2.14$.
\Cref{fig-mms-fields} shows a visualization of this solution.

\begin{table}[htbp]
\centering
\begin{tabular*}{\linewidth}{@{\extracolsep{\fill}}*8c@{}}
\toprule
Equation, $\varphi$ & $\varphi^0$ & $\varphi^{x_1}$ & $\varphi^{x_2}$ & $\varphi^{x_3}$ & $\alpha_{\varphi}^{x_1}$ & $\alpha_{\varphi}^{x_2}$ & $\alpha_{\varphi}^{x_3}$ \\
\midrule
$\rho$              & 1           & 3/20           & -1/10         & 1/5   & 1   & 1/2 & 1   \\
$\mathcal{U}_1$     & 800         & 50             & -30           & 50    & 3/2 & 3/5 & 1   \\
$\mathcal{U}_2$     & 800         & -75            & 40            & 10    & 1/2 & 3/2 & 3/2 \\
$\mathcal{U}_3$     & 800         & 15             & -25           & 20    & 3/2 & 1/2 & 5/4 \\
$\mathcal{P}$       & $10^5$      & $2\times 10^4$ &$5\times 10^4$ & $2/3$ & 1   & 3/2 & 1   \\
\bottomrule
\end{tabular*}
\caption{Constants for the 3D compressible Navier--Stokes supersonic manufactured solution.}
\label{tab_3dmms_ns}
\end{table}

\newcommand\annotate[1]{\drawannotate(#1)}%Command to produce annotations based on path
\def\drawannotate(#1,#2){
\path (axis cs:#1,#2)-- +(5pt,5pt) node[rotate=60,scale=.3,pos=.25] {\textless} node[font=\tiny] {#2};
}

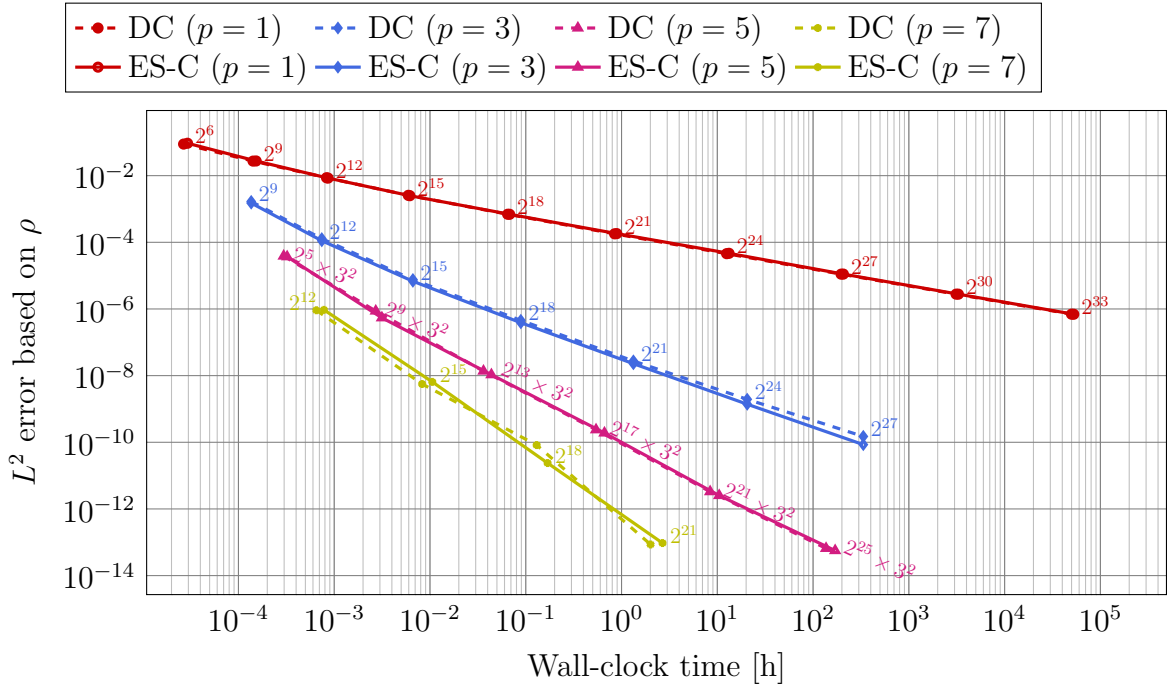
\begin{figure}[!ht]
\centering
\begin{tikzpicture}
\definecolor{color0}{rgb}{0.75,0.75,0}
\definecolor{color1}{rgb}{0.75,0,0.75}
\begin{axis}[height=8cm,width=15.0cm,legend cell align={left},legend entries={DC ($p=1$),DC ($p=3$),DC ($p=5$),DC ($p=7$),ES-C ($p=1$),ES-C ($p=3$),ES-C ($p=5$),ES-C ($p=7$)},legend style={at={(-0.08,1.13)},anchor=west},legend columns=4,log basis x={10},log basis y={10},tick align=outside,tick pos=left,xlabel={Wall-clock time [h]},ylabel={$L^2$ error based on $\rho$},xmajorgrids,xminorgrids,xmode=log,ymode=log,grid=both,major grid style={black!50},major tick length=0pt,minor tick length=0pt,ymax=0.9,xmin=0.000011]
\addplot [very thick, mycolor1, dashed, mark=*, mark size=1.5, mark options={solid}]table[x index=0,y index=1]{data/MMS/order1.dat};
\addplot [very thick, mycolor2, dashed, mark=diamond*, mark size=1.5, mark options={solid}]table[x index=0,y index=1]{data/MMS/order3.dat};
\addplot [very thick, mycolor3, dashed, mark=triangle*, mark size=1.5, mark options={solid}]table[x index=0,y index=1]{data/MMS/order5.dat};
\addplot [very thick, color0, dashed, mark=asterisk, mark size=1.5, mark options={solid}]table[x index=0,y index=1]{data/MMS/order7.dat};
\addplot [very thick, mycolor1, mark=o, mark size=1.5, mark options={fill=white},text mark as node=true]table[x index=2,y index=3]{data/MMS/order1.dat};
\addplot [very thick, mycolor2, mark=diamond, mark size=1.5, mark options={fill=white}]table[x index=2,y index=3]{data/MMS/order3.dat};
\addplot [very thick, mycolor3, mark=triangle, mark size=1.5, mark options={fill=white}]table[x index=2,y index=3]{data/MMS/order5.dat};
\addplot [very thick, color0, mark=asterisk, mark size=1.5, mark options={fill=white}]table[x index=2,y index=3]{data/MMS/order7.dat};

\draw[] (axis cs:2.93822222222222e-05,0.0937401)        -- (axis cs:2.93822222222222e-05,0.0937401); \node at (axis cs:2.93822222222222e-05,0.0937401)[scale=0.70,anchor=base west,text=mycolor1,rotate=0.0]{$2^{6}$};
\draw[] (axis cs:0.000152893333333333,0.0275705)        -- (axis cs:0.000152893333333333,0.0275705); \node at (axis cs:0.000152893333333333,0.0275705)[scale=0.70,anchor=base west,text=mycolor1,rotate=0.0]{$2^{9}$};
\draw[] (axis cs:0.000862311111111111,0.00844906)       -- (axis cs:0.000862311111111111,0.00844906);\node at (axis cs:0.000862311111111111,0.00844906)[scale=0.70,anchor=base west,text=mycolor1,rotate=0.0]{$2^{12}$};
\draw[] (axis cs:0.00616106666666667,0.00250349)        -- (axis cs:0.00616106666666667,0.00250349);\node at (axis cs:0.00616106666666667,0.00250349)[scale=0.70,anchor=base west,text=mycolor1,rotate=0.0]{$2^{15}$};
\draw[] (axis cs:0.0678364444444444,0.000692466)        -- (axis cs:0.0678364444444444,0.000692466);\node at (axis cs:0.0678364444444444,0.000692466)[scale=0.70,anchor=base west,text=mycolor1,rotate=0.0]{$2^{18}$};
\draw[] (axis cs:0.892124444444444,0.000182258)         -- (axis cs:0.892124444444444,0.000182258); \node at (axis cs:0.892124444444444,0.000182258)[scale=0.70,anchor=base west,text=mycolor1,rotate=0.0]{$2^{21}$};
\draw[] (axis cs:13.2472888888889,4.67408e-05)          -- (axis cs:13.2472888888889,4.67408e-05);\node at (axis cs:13.2472888888889,4.67408e-05)[scale=0.70,anchor=base west,text=mycolor1,rotate=0.0]{$2^{24}$};
\draw[] (axis cs:207.018666666667,1.11188e-05)          -- (axis cs:207.018666666667,1.11188e-05);\node at (axis cs:207.018666666667,1.11188e-05)[scale=0.70,anchor=base west,text=mycolor1,rotate=0.0]{$2^{27}$};
\draw[] (axis cs:3287.89333333333,2.77907720185137e-06) -- (axis cs:3287.89333333333,2.77907720185137e-06);\node at (axis cs:3287.89333333333,2.77907720185137e-06)[scale=0.70,anchor=base west,text=mycolor1,rotate=0.0]{$2^{30}$};
\draw[] (axis cs:52689.5497299883,6.94492126918759e-07) -- (axis cs:52689.5497299883,6.94492126918759e-07);\node at (axis cs:52689.5497299883,6.94492126918759e-07)[scale=0.70,anchor=base west,text=mycolor1,rotate=0.0]{$2^{33}$};
\draw[] (axis cs:0.000135877777777778,0.0016768)        -- (axis cs:0.000135877777777778,0.0016768);\node at (axis cs:0.000135877777777778,0.0016768)[scale=0.70,anchor=base west,text=mycolor2,rotate=0.0]{$2^{9}$};
\draw[] (axis cs:0.000741688888888889,0.000128571)      -- (axis cs:0.000741688888888889,0.000128571);\node at (axis cs:0.000741688888888889,0.000128571)[scale=0.70,anchor=base west,text=mycolor2,rotate=0.0]{$2^{12}$};
\draw[] (axis cs:0.00663102222222222,7.61852e-06)       -- (axis cs:0.00663102222222222,7.61852e-06);\node at (axis cs:0.00663102222222222,7.61852e-06)[scale=0.70,anchor=base west,text=mycolor2,rotate=0.0]{$2^{15}$};
\draw[] (axis cs:0.0888124444444444,4.54158e-07)        -- (axis cs:0.0888124444444444,4.54158e-07);\node at (axis cs:0.0888124444444444,4.54158e-07)[scale=0.70,anchor=base west,text=mycolor2,rotate=0.0]{$2^{18}$};
\draw[] (axis cs:1.33009777777778,2.75305e-08)          -- (axis cs:1.33009777777778,2.75305e-08);\node at (axis cs:1.33009777777778,2.75305e-08)[scale=0.70,anchor=base west,text=mycolor2,rotate=0.0]{$2^{21}$};
\draw[] (axis cs:20.5219555555556,1.94706e-09)          -- (axis cs:20.5219555555556,1.94706e-09);\node at (axis cs:20.5219555555556,1.94706e-09)[scale=0.70,anchor=base west,text=mycolor2,rotate=0.0]{$2^{24}$};
\draw[] (axis cs:334.520888888889,1.49844e-10)          -- (axis cs:334.520888888889,1.49844e-10);\node at (axis cs:334.520888888889,1.49844e-10)[scale=0.70,anchor=base west,text=mycolor2,rotate=0.0]{$2^{27}$};
\draw[] (axis cs:0.000296636666666667,3.86568e-05)      -- (axis cs:0.000296636666666667,3.86568e-05);\node at (axis cs:0.000296636666666667,3.86568e-05)[scale=0.70,anchor=base west,text=mycolor3,rotate=-25.0]{$2^{5}\times 3^2$};
\draw[] (axis cs:0.00271473333333333,8.67154e-07)       -- (axis cs:0.00271473333333333,8.67154e-07);\node at (axis cs:0.00271473333333333,8.67154e-07)[scale=0.70,anchor=base west,text=mycolor3,rotate=-25]{$2^{9}\times 3^2$};
\draw[] (axis cs:0.036188,1.3674e-08)                   -- (axis cs:0.036188,1.3674e-08);\node at (axis cs:0.045188,1.3674e-08)[scale=0.70,anchor=base west,text=mycolor3,rotate=-25]{$2^{13}\times 3^2$};
\draw[] (axis cs:0.541706666666667,2.37275e-10)         -- (axis cs:0.541706666666667,2.37275e-10);\node at (axis cs:0.611706666666667,2.37275e-10)[scale=0.70,anchor=base west,text=mycolor3,rotate=-25]{$2^{17}\times 3^2$};
\draw[] (axis cs:8.41706666666667,3.29433e-12)          -- (axis cs:8.41706666666667,3.29433e-12);\node at (axis cs:9.39706666666667,3.29433e-12)[scale=0.70,anchor=base west,text=mycolor3,rotate=-25]{$2^{21}\times 3^2$};
\draw[] (axis cs:136.650666666667,6.54874e-14)          -- (axis cs:136.650666666667,6.54874e-14);\node at (axis cs:170.950666666667,6.54874e-14)[scale=0.70,anchor=base west,text=mycolor3,rotate=-25]{$2^{25}\times 3^2$};
\draw[] (axis cs:0.000781466666666667,9.39503e-07)      -- (axis cs:0.000781466666666667,9.39503e-07);\node at (axis cs:0.000281466666666667,9.39503e-07)[scale=0.70,anchor=base west,text=color0,rotate=0.0]{$2^{12}$};
\draw[] (axis cs:0.0105662222222222,6.49191e-09)        -- (axis cs:0.0105662222222222,6.49191e-09);\node at (axis cs:0.0105662222222222,6.49191e-09)[scale=0.70,anchor=base west,text=color0,rotate=0.0]{$2^{15}$};
\draw[] (axis cs:0.169134222222222,2.38087e-11)         -- (axis cs:0.169134222222222,2.38087e-11);\node at (axis cs:0.169134222222222,2.38087e-11)[scale=0.70,anchor=base west,text=color0,rotate=0.0]{$2^{18}$};
\draw[] (axis cs:2.67818666666667,9.46888e-14)          -- (axis cs:2.67818666666667,9.46888e-14);\node at (axis cs:2.67818666666667,9.46888e-14)[scale=0.70,anchor=base west,text=color0,rotate=0.0]{$2^{21}$};
\end{axis}
\end{tikzpicture}
\caption{Error vs. cost for 3D manufactured solution.}
\label{fig_mms_costVSerror}
\end{figure}

\Cref{fig_mms_costVSerror} shows the results for the $L^2$ error of the density field against the time to solution for several solution polynomial degrees.
The numbers superimposed to the data series represent the total number of degrees of freedom (DOFs) of each mesh, showing that the wall-clock time is almost constant given a total of DOFs, regardless of the polynomial degree.
This observation added to the accuracy of high order solutions can result in higher order discretizations being more economical, given some acceptable error threshold.
This observation should serve as a motivation to further explore the behavior of high order DC-type discretizations in more applied settings.

%%%%%%%%%%%%%%%%%%%%%%%%%%%%%%%%%%%%%%%%%%%%%%%%%%%%%%%%%%%%%%%%%%%%%%%%%%%%%%%%
\subsection{The cost of entropy stability}

\begin{figure}[!ht]
\centering
\pgfplotstableread[row sep=\\,col sep=&]{
degree & ES         & KG   & SDG       \\
1      & 3.1798E+00 & 2.8634e+00 & 2.7771E+00 \\
3      & 2.5059E+00 & 2.01675    & 1.9593E+00 \\
7      & 2.6099E+00 & 1.73415    & 1.5335E+00 \\
15     & 3.5589E+00 & 1.9914e+00 & 1.4830E+00 \\
%31     & 6.6816E+00 & 3.2718e+00 & 4.7928E+00 \\
}\TGVtimedata
\begin{tikzpicture}
\begin{axis}[ybar,bar width = 0.7cm,width = 0.8\textwidth,height = 0.5\textwidth,legend style={at={(0.5,1)},anchor=north,legend columns=-1},symbolic x coords={1,3,7,15},xtick=data,nodes near coords,nodes near coords align={vertical},ylabel={Wall-clock time per TS [s]},xlabel={Degree ($p$)},enlarge x limits=0.18,legend columns=1,legend cell align={left},ymin=0 ]
\addplot[black,fill=mycolor11] table[x=degree,y=ES]{\TGVtimedata};
\addplot[black,fill=mycolor2] table[x=degree,y=KG]{\TGVtimedata};
\addplot[black,fill=mycolor3] table[x=degree,y=SDG]{\TGVtimedata};
\legend{ES-C,SF-KG,DC}
\end{axis}
\end{tikzpicture}
\caption{Wall-clock time per time step for different solution polynomial degrees and discretizations.}
\label{fig_scaling}
\end{figure}
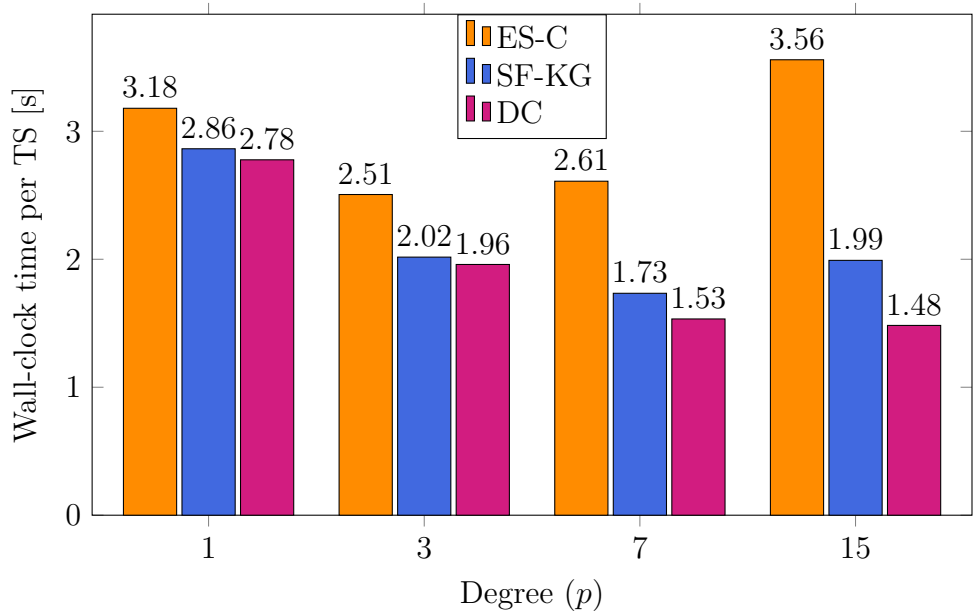

We devised a test to compare the cost of the entropy stable discretization with the non-entropy stable ones.
In this test, the Taylor--Green vortex (TGV) at $\Rey=1{,}600$ and $\Ma=0.05$ is simulated using a fixed number of DOFs for four different solution polynomial degrees. The initial condition reads
\begin{equation}
\begin{cases}
\rho&= 1 + \frac{\gamma \Ma^2}{16} \left(\cos(2x_1)+\cos(2x_2)\right)\left(\cos(2x_3)+2\right), \\
\mathcal{U}_1&= \sin(x_1)\cos(x_2)\cos(x_3), \\
\mathcal{U}_2&= -\cos(x_1)\sin(x_2)\cos(x_3),\\
\mathcal{U}_3&= 0,\\
\T &= 1.
\end{cases}
\end{equation}
This problem is generally used to study the transition to turbulence, the energy decay in turbulent flow, the simulation of under-resolved turbulent features and the robustness of solvers \cite{debonis2013solutions,moura2015dg,Parsani2016,Gassner2016}.
Here we use it as a nontrivial flow pattern where we measure the cost of our discretizations.
The mesh is a structured cube domain $\Omega=[-\pi,\pi]^3$ with periodic boundary conditions, discretized with enough elements to guarantee a total of 256 DOFs per direction (for a total of $256^3$ DOFs).
We use two nodes of Shaheen XC 40 \cite{shaheen_2} with a perfectly balanced partition\footnote{Each Shaheen XC 40 node has 32 Intel Haswell cores.}. In fact, in our solver, the parallel overhead for this number of nodes is completely negligible.
During the execution, a fixed small time step is used and all I/O operations are turned off. The time per time step (TS) is averaged until no noticeable change is observed.
The final results are shown in Figure \ref{fig_scaling}.

All the results for wall-clock time per time step fall within the same order of magnitude.
This similarity is consistent with the results shown in the \Cref{sec:numerical_results} since the computational cost of the discretization mostly depends on the number of DOFs.
The differences within this range are caused by two main factors, namely the number of flux evaluations in the volume terms (with its complexity), and the ratio of the volume terms to the surface terms (SATs).
Here the former causes discretizations to be more expensive with increasing polynomial degree, while the latter has the opposite effect.
The results for both ES-C and SF-KG discretizations show this with a minimum time per time step at $p=3$ and $p=7$, respectively, with the SF-KG flux being cheaper.
A minimum time per time step for the conventional DC discretization is not found within the studied range of solution polynomial degrees.
Obviously, the conventional DC discretization is cheaper because it is characterized by fewer inviscid flux evaluations for all the solution polynomial degree, as discussed in \Cref{sec:compressible_nse}.

%%%%%%%%%%%%%%%%%%%%%%%%%%%%%%%%%%%%%%%%%%%%%%%%%%%%%%%%%%%%%%%%%%%%%
\subsection{Robustness}
Although high-order accurate methods are well suited for smooth solutions, numerical instabilities may occur if the flow contains \textit{under-resolved physical features} (e.g., under-resolved turbulent flows) 
or \textit{discontinuities} (e.g., shocks). In this section we investigate the robusteness of the three discretizations for under-resolved turbulent and non-smooth flows.

\definecolor{g}{rgb}{0.0, 0.68, 0.15}
\definecolor{r}{rgb}{0.73, 0.0, 0.0}

\subsubsection{Taylor--Green vortex}

The Taylor--Green vortex (TGV) problem at $\Rey=1{,}600$ and $\Ma=0.05$ is a flow that degenerates to turbulence over time; therefore its solution is representative of the behavior 
of the algorithm for the solution of turbulent flows.
By simply using coarse grids, this problem can be made a challenging benchmark for the simulation of under-resolved turbulent features, a situation that easily leads to stability issues with high order operators.
The study is performed using several nested meshes (see the first column in \Cref{tab_3D_TGV}) and solution polynomial degrees from $p=1$ to $p=15$. 
The final time is set to $t_f=20$.
\Cref{tab_3D_TGV} summarizes the results.

\begin{table}[!ht]
\centering
\scriptsize
\begin{tabu} to \columnwidth {lc|X[c]|X[c]|X[c]|X[c]|X[c]|X[c]|X[c]|X[l]|X[l]|X[l]|X[l]|X[l]|X[l]|X[l]|X[l]}
\hline\hline
& {\scriptsize Degree}& 1 & 2 & 3 & 4 & 5 & 6 & 7 & 8 & 9 & 10 & 11 & 12 & 13 & 14 & 15 \\
\hline\hline
{\scriptsize ES-C} & $3^3$ &
\textcolor{mycolor2}{\Checkmark}&
\textcolor{mycolor2}{\Checkmark}&
\textcolor{mycolor2}{\Checkmark}&
\textcolor{mycolor2}{\Checkmark}&
\textcolor{mycolor2}{\Checkmark}&
\textcolor{mycolor2}{\Checkmark}&
\textcolor{mycolor2}{\Checkmark}&
\textcolor{mycolor2}{\Checkmark}&
\textcolor{mycolor2}{\Checkmark}&
\textcolor{mycolor2}{\Checkmark}&
\textcolor{mycolor2}{\Checkmark}&
\textcolor{mycolor2}{\Checkmark}&
\textcolor{mycolor2}{\Checkmark}&
\textcolor{mycolor2}{\Checkmark}&
\textcolor{mycolor2}{\Checkmark}\\
&    $6^3$ &
\textcolor{mycolor2}{\Checkmark}&
\textcolor{mycolor2}{\Checkmark}&
\textcolor{mycolor2}{\Checkmark}&
\textcolor{mycolor2}{\Checkmark}&
\textcolor{mycolor2}{\Checkmark}&
\textcolor{mycolor2}{\Checkmark}&
\textcolor{mycolor2}{\Checkmark}&
\textcolor{mycolor2}{\Checkmark}&
\textcolor{mycolor2}{\Checkmark}&
\textcolor{mycolor2}{\Checkmark}&
\textcolor{mycolor2}{\Checkmark}&
\textcolor{mycolor2}{\Checkmark}&
\textcolor{mycolor2}{\Checkmark}&
\textcolor{mycolor2}{\Checkmark}&
\textcolor{mycolor2}{\Checkmark}\\
&    $12^3$ &
\textcolor{mycolor2}{\Checkmark}&
\textcolor{mycolor2}{\Checkmark}&
\textcolor{mycolor2}{\Checkmark}&
\textcolor{mycolor2}{\Checkmark}&
\textcolor{mycolor2}{\Checkmark}&
\textcolor{mycolor2}{\Checkmark}&
\textcolor{mycolor2}{\Checkmark}&
\textcolor{mycolor2}{\Checkmark}&
\textcolor{mycolor2}{\Checkmark}&
\textcolor{mycolor2}{\Checkmark}&
\textcolor{mycolor2}{\Checkmark}&
\textcolor{mycolor2}{\Checkmark}&
\textcolor{mycolor2}{\Checkmark}&
\textcolor{mycolor2}{\Checkmark}&
\textcolor{mycolor2}{\Checkmark}\\
&    $24^3$ &
\textcolor{mycolor2}{\Checkmark}&
\textcolor{mycolor2}{\Checkmark}&
\textcolor{mycolor2}{\Checkmark}&
\textcolor{mycolor2}{\Checkmark}&
\textcolor{mycolor2}{\Checkmark}&
\textcolor{mycolor2}{\Checkmark}&
\textcolor{mycolor2}{\Checkmark}&
\textcolor{mycolor2}{\Checkmark}&
\textcolor{mycolor2}{\Checkmark}&
\textcolor{mycolor2}{\Checkmark}&
\textcolor{mycolor2}{\Checkmark}&
\textcolor{mycolor2}{\Checkmark}&
\textcolor{mycolor2}{\Checkmark}&
\textcolor{mycolor2}{\Checkmark}&
\textcolor{mycolor2}{\Checkmark}\\
&    $48^3$ &
\textcolor{mycolor2}{\Checkmark}&
\textcolor{mycolor2}{\Checkmark}&
\textcolor{mycolor2}{\Checkmark}&
\textcolor{mycolor2}{\Checkmark}&
\textcolor{mycolor2}{\Checkmark}&
\textcolor{mycolor2}{\Checkmark}&
\textcolor{mycolor2}{\Checkmark}&
\textcolor{mycolor2}{\Checkmark}&
\textcolor{mycolor2}{\Checkmark}&
\textcolor{mycolor2}{\Checkmark}&
\textcolor{mycolor2}{\Checkmark}&
\textcolor{mycolor2}{\Checkmark}&
\textcolor{mycolor2}{\Checkmark}&
\textcolor{mycolor2}{\Checkmark}&
\textcolor{mycolor2}{\Checkmark}\\
\hline\hline
{\scriptsize SF-KG} & $3^3$ &
\textcolor{mycolor2}{\Checkmark}&
\textcolor{mycolor2}{\Checkmark}&
\textcolor{mycolor2}{\Checkmark}&
\textcolor{mycolor2}{\Checkmark}&
\textcolor{mycolor2}{\Checkmark}&
\textcolor{mycolor2}{\Checkmark}&
\textcolor{mycolor2}{\Checkmark}&
\textcolor{mycolor2}{\Checkmark}&
\textcolor{mycolor2}{\Checkmark}&
\textcolor{mycolor2}{\Checkmark}&
\textcolor{mycolor2}{\Checkmark}&
\textcolor{mycolor2}{\Checkmark}&
\textcolor{mycolor2}{\Checkmark}&
\textcolor{mycolor2}{\Checkmark}&
\textcolor{mycolor2}{\Checkmark}\\
& $6^3$ &
\textcolor{mycolor2}{\Checkmark}&
\textcolor{mycolor2}{\Checkmark}&
\textcolor{mycolor2}{\Checkmark}&
\textcolor{mycolor2}{\Checkmark}&
\textcolor{mycolor2}{\Checkmark}&
\textcolor{mycolor2}{\Checkmark}&
\textcolor{mycolor2}{\Checkmark}&
\textcolor{mycolor2}{\Checkmark}&
\textcolor{mycolor2}{\Checkmark}&
\textcolor{mycolor2}{\Checkmark}&
\textcolor{mycolor2}{\Checkmark}&
\textcolor{mycolor2}{\Checkmark}&
\textcolor{mycolor2}{\Checkmark}&
\textcolor{mycolor2}{\Checkmark}&
\textcolor{mycolor2}{\Checkmark}\\
& $12^3$ &
\textcolor{mycolor2}{\Checkmark}&
\textcolor{mycolor2}{\Checkmark}&
\textcolor{mycolor2}{\Checkmark}&
\textcolor{mycolor2}{\Checkmark}&
\textcolor{mycolor2}{\Checkmark}&
\textcolor{mycolor2}{\Checkmark}&
\textcolor{mycolor2}{\Checkmark}&
\textcolor{mycolor2}{\Checkmark}&
\textcolor{mycolor2}{\Checkmark}&
\textcolor{mycolor2}{\Checkmark}&
\textcolor{mycolor2}{\Checkmark}&
\textcolor{mycolor2}{\Checkmark}&
\textcolor{mycolor2}{\Checkmark}&
\textcolor{mycolor2}{\Checkmark}&
\textcolor{mycolor2}{\Checkmark}\\
& $24^3$ &
\textcolor{mycolor2}{\Checkmark}&
\textcolor{mycolor2}{\Checkmark}&
\textcolor{mycolor2}{\Checkmark}&
\textcolor{mycolor2}{\Checkmark}&
\textcolor{mycolor2}{\Checkmark}&
\textcolor{mycolor2}{\Checkmark}&
\textcolor{mycolor2}{\Checkmark}&
\textcolor{mycolor2}{\Checkmark}&
\textcolor{mycolor2}{\Checkmark}&
\textcolor{mycolor2}{\Checkmark}&
\textcolor{mycolor2}{\Checkmark}&
\textcolor{mycolor2}{\Checkmark}&
\textcolor{mycolor2}{\Checkmark}&
\textcolor{mycolor2}{\Checkmark}&
\textcolor{mycolor2}{\Checkmark}\\
& $48^3$ &
\textcolor{mycolor2}{\Checkmark}&
\textcolor{mycolor2}{\Checkmark}&
\textcolor{mycolor2}{\Checkmark}&
\textcolor{mycolor2}{\Checkmark}&
\textcolor{mycolor2}{\Checkmark}&
\textcolor{mycolor2}{\Checkmark}&
\textcolor{mycolor2}{\Checkmark}&
\textcolor{mycolor2}{\Checkmark}&
\textcolor{mycolor2}{\Checkmark}&
\textcolor{mycolor2}{\Checkmark}&
\textcolor{mycolor2}{\Checkmark}&
\textcolor{mycolor2}{\Checkmark}&
\textcolor{mycolor2}{\Checkmark}&
\textcolor{mycolor2}{\Checkmark}&
\textcolor{mycolor2}{\Checkmark}\\
\hline\hline
{\scriptsize DC} & $3^3$ &
\textcolor{mycolor2}{\Checkmark}&
\textcolor{mycolor2}{\Checkmark}&
\textcolor{mycolor1}{$\bm{\times}$}&
\textcolor{mycolor1}{$\bm{\times}$}&
\textcolor{mycolor1}{$\bm{\times}$}&
\textcolor{mycolor1}{$\bm{\times}$}&
\textcolor{mycolor1}{$\bm{\times}$}&
\textcolor{mycolor1}{$\bm{\times}$}&
\textcolor{mycolor1}{$\bm{\times}$}&
\textcolor{mycolor1}{$\bm{\times}$}&
\textcolor{mycolor1}{$\bm{\times}$}&
\textcolor{mycolor1}{$\bm{\times}$}&
\textcolor{mycolor1}{$\bm{\times}$}&
\textcolor{mycolor1}{$\bm{\times}$}&
\textcolor{mycolor1}{$\bm{\times}$}\\
& $6^3$ &
\textcolor{mycolor2}{\Checkmark}&
\textcolor{mycolor2}{\Checkmark}&
\textcolor{mycolor1}{$\bm{\times}$}&
\textcolor{mycolor1}{$\bm{\times}$}&
\textcolor{mycolor1}{$\bm{\times}$}&
\textcolor{mycolor1}{$\bm{\times}$}&
\textcolor{mycolor1}{$\bm{\times}$}&
\textcolor{mycolor1}{$\bm{\times}$}&
\textcolor{mycolor1}{$\bm{\times}$}&
\textcolor{mycolor1}{$\bm{\times}$}&
\textcolor{mycolor1}{$\bm{\times}$}&
\textcolor{mycolor1}{$\bm{\times}$}&
\textcolor{mycolor1}{$\bm{\times}$}&
\textcolor{mycolor1}{$\bm{\times}$}&
\textcolor{mycolor1}{$\bm{\times}$}\\
& $12^3$ &
\textcolor{mycolor2}{\Checkmark}&
\textcolor{mycolor2}{\Checkmark}&
\textcolor{mycolor1}{$\bm{\times}$}&
\textcolor{mycolor1}{$\bm{\times}$}&
\textcolor{mycolor1}{$\bm{\times}$}&
\textcolor{mycolor1}{$\bm{\times}$}&
\textcolor{mycolor1}{$\bm{\times}$}&
\textcolor{mycolor1}{$\bm{\times}$}&
\textcolor{mycolor1}{$\bm{\times}$}&
\textcolor{mycolor1}{$\bm{\times}$}&
\textcolor{mycolor1}{$\bm{\times}$}&
\textcolor{mycolor1}{$\bm{\times}$}&
\textcolor{mycolor1}{$\bm{\times}$}&
\textcolor{mycolor1}{$\bm{\times}$}&
\textcolor{mycolor1}{$\bm{\times}$}\\
& $24^3$ &
\textcolor{mycolor2}{\Checkmark}&
\textcolor{mycolor2}{\Checkmark}&
\textcolor{mycolor2}{\Checkmark}&
\textcolor{mycolor1}{$\bm{\times}$}&
\textcolor{mycolor1}{$\bm{\times}$}&
\textcolor{mycolor1}{$\bm{\times}$}&
\textcolor{mycolor1}{$\bm{\times}$}&
\textcolor{mycolor1}{$\bm{\times}$}&
\textcolor{mycolor1}{$\bm{\times}$}&
\textcolor{mycolor2}{\Checkmark}&
\textcolor{mycolor1}{$\bm{\times}$}&
\textcolor{mycolor2}{\Checkmark}&
\textcolor{mycolor2}{\Checkmark}&
\textcolor{mycolor2}{\Checkmark}&
\textcolor{mycolor2}{\Checkmark}\\
& $48^3$ &
\textcolor{mycolor2}{\Checkmark}&
\textcolor{mycolor2}{\Checkmark}&
\textcolor{mycolor2}{\Checkmark}&
\textcolor{mycolor2}{\Checkmark}&
\textcolor{mycolor1}{$\bm{\times}$}&
\textcolor{mycolor1}{$\bm{\times}$}&
\textcolor{mycolor2}{\Checkmark}&
\textcolor{mycolor2}{\Checkmark}&
\textcolor{mycolor2}{\Checkmark}&
\textcolor{mycolor2}{\Checkmark}&
\textcolor{mycolor2}{\Checkmark}&
\textcolor{mycolor2}{\Checkmark}&
\textcolor{mycolor2}{\Checkmark}&
\textcolor{mycolor2}{\Checkmark}&
\textcolor{mycolor2}{\Checkmark}\\
& $96^3$ &
\textcolor{mycolor2}{\Checkmark}&
\textcolor{mycolor2}{\Checkmark}&
\textcolor{mycolor2}{\Checkmark}&
\textcolor{mycolor2}{\Checkmark}&
\textcolor{mycolor2}{\Checkmark}&
\textcolor{mycolor2}{\Checkmark}&
\textcolor{mycolor2}{\Checkmark}&--&--&--&--&--&--&--&--\\
\end{tabu}
\caption{Numerical stability for the 3D Taylor--Green vortex. \textcolor{mycolor2}{\Checkmark}= success, \textcolor{mycolor1}{$\bm{\times}$} = failure.}
\label{tab_3D_TGV}
\end{table}

As expected, the DC discretization is frequently unstable and leads to a crash of the solver when coarse meshes with high order solution polynomial degrees are used. 
On the contrary, the other two methods are always stable.
We emphasize that in our test, only the ES-C discretization has a rigorous proof of stability and the particular conditions of the test may be the cause for the perfect score of SF-KG.
So it is natural to expect there is a set of problems where even SF-KG is bound to fail.
The next section explores further in this difference.

%%%%%%%%%%%%%%%%%%%%%%%%%%%%%%%%%%%%%%%%%%%%%%%%%%%%%%%%%%%%%%%%%%%%%%%%%%%%%%%%%

\subsubsection{Compressible homogeneous isotropic turbulence}

We select this 3D case as a more challenging problem to further stress the stability properties of our discretizations.
This problem is regarded as one of the cornerstones to elucidate the effects of compressibility for compressible turbulence \cite{lele1994compressibility}.
Based on the previous numerical experiments and theoretical analyses, isotropic compressible turbulence is divided into four main dynamical regimes \cite{sagaut2008homogeneous}, i.e., the low-Mach number quasi-isentropic regime, the low-Mach number thermal regime, the nonlinear subsonic regime, and the supersonic regime.
Most of the numerical schemes developed in the last decade utilized in the simulation of isotropic compressible turbulence with moderate turbulent Mach number to supersonic regimes fail to capture shocklets robustly and accurately resolve smooth regions.
For isotropic compressible turbulence in these two regimes, the stronger random shocklets and higher spatial-temporal gradients pose greater difficulties for numerical analyses than other regimes.
Both forced isotropic compressible turbulence with solenoidal and dilatational external force \cite{kida1990energy,jagannathan2016reynolds} and decaying isotropic compressible turbulence \cite{samtaney2001direct,pirozzoli2004direct,boukharfane2018contribution} are studied in the literature.
In this paper, we choose to test with the decaying isotropic compressible turbulence without external force.
The flow domain of numerical simulation is a cube box defined as $\Omega=[0,2\pi]^3$, with periodic boundary conditions in all directions.
The evolution of this artificial system is determined by initial thermodynamic quantities and two dimensionless parameters, i.e., the initial Taylor microscale Reynolds number $\Rey=\langle\rho\rangle \mathcal{U}_{\mathrm{RMS}}\lambda/\langle\mu\rangle$, and turbulent Mach number $\Ma_t=\mathcal{U}_{\mathrm{RMS}}/\langle c_s\rangle$, where $\langle .\rangle$ is the ensemble over the whole computational domain, $\rho$ is the density, $\mu$ is the initial dynamic viscosity, $c_s$ is the sound speed, $\mathcal{U}_{\mathrm{RMS}}$ is the root mean square of initial turbulent velocity field $\mathcal{U}_{\mathrm{RMS}}=\langle \mathcal{\bm{U}}\cdot\mathcal{\bm{U}}/3\rangle^{1/2}$, and the normalized Taylor micro-scale $\lambda$ is defined by
\begin{equation}
\lambda=\sqrt{\frac{\mathcal{U}_{\mathrm{RMS}}}{\left\langle\left(\frac{\partial\mathcal{U}_1}{\partial x_1}\right)^2+\left(\frac{\partial\mathcal{U}_2}{\partial x_2}\right)^2+\left(\frac{\partial\mathcal{U}_3}{\partial x_3}\right)^2\right\rangle}}.
\end{equation} 

A 3D solenoidal random initial velocity field $\mathcal{\bm{U}}$ can be generated by a specified spectrum, which is given by \cite{passot1987numerical}
\begin{equation}\label{eq_passot_pouquet}
E(k)=A_0k^4\exp\left(-2k^2/k_0^2\right),
\end{equation}
where $A_0$ is a constant to get a specified initial kinetic energy, $k$ is the wave number, $k_0$ is the wave number at which the spectrum peaks.
In this paper, fixed $A_0$ and $k_0$ in Equation \eqref{eq_passot_pouquet} are chosen for all cases, which are initialized by $A_0=0.00013$ and $k_0=8$.
Initial strategies play an important role in isotropic compressible turbulence simulation \cite{samtaney2001direct}, especially for the starting fast transient period during which the divergence of the velocity increases rapidly and negative temperature or pressure often appear.
In the computation, the initial pressure $\mathcal{P}_0$, density $\rho_0$ and temperature $\mathcal{T}_0$ are set following the procedure of \citet{ristorcelli1997consistent}.
For higher polynomial degree, the hydrodynamic and thermodynamic fields are computed from polynomial order $p=1$ and then mapped onto the particular grid.
Simulations are run for a relatively short time interval which is the three order of magnitude of a characteristic eddy turnover time $\tau=\mathcal{L}_1/\mathcal{U}_{\mathrm{RMS}}$ where $\mathcal{L}_1$ being the integral length scale defined as
\begin{equation}
\mathcal{L}_1=\frac{3\pi}{4}\frac{\int_0^\infty E(k)\mathrm{d}k/k}{\int_0^\infty E(k)\mathrm{d}k}.
\end{equation}
To assess the relative robustness of the discretizations, the turbulent Mach number is set to $\Ma_t=0.62$, where shocklets appear and thus the flow belong to the nonlinear subsonic regime in which compressible effects are important.
The number of wavenumber available in the computational box is intentionally limited to values that do not allow a correct representation of the smallest scales of the initially imposed turbulent field.
In all flow simulations the turbulent Reynolds number is $\Rey_\lambda=194$.

Several test runs at multiple resolutions and polynomial orders are performed to assess the numerical robustness of the conventional DC, SF-KG and ES-C schemes.

\begin{table}[!ht]
\centering
\scriptsize
\begin{tabular}{lc c|c|c|c}
\hline\hline
&Degree& 1 & 2 & 3 & 4  \\
\hline\hline
ES-C & $4^3$ &
\textcolor{mycolor2}{\Checkmark}&
\textcolor{mycolor2}{\Checkmark}&
\textcolor{mycolor2}{\Checkmark}&
\textcolor{mycolor2}{\Checkmark} \\
&    $8^3$ &
\textcolor{mycolor2}{\Checkmark}&
\textcolor{mycolor2}{\Checkmark}&
\textcolor{mycolor2}{\Checkmark}&
\textcolor{mycolor2}{\Checkmark} \\
&    $16^3$ &
\textcolor{mycolor2}{\Checkmark}&
\textcolor{mycolor2}{\Checkmark}&
\textcolor{mycolor2}{\Checkmark}&
\textcolor{mycolor2}{\Checkmark} \\
&    $32^3$ &
\textcolor{mycolor2}{\Checkmark}&
\textcolor{mycolor2}{\Checkmark}&
\textcolor{mycolor2}{\Checkmark}&
\textcolor{mycolor2}{\Checkmark} \\
&    $64^3$ &
\textcolor{mycolor2}{\Checkmark}&
\textcolor{mycolor2}{\Checkmark}&
\textcolor{mycolor2}{\Checkmark}&
\textcolor{mycolor2}{\Checkmark} \\
\hline\hline
SF-KG & $4^3$ &
\textcolor{mycolor1}{$\bm{\times}$}&
\textcolor{mycolor2}{\Checkmark}&
\textcolor{mycolor2}{\Checkmark}&
\textcolor{mycolor2}{\Checkmark} \\
&    $8^3$ &
\textcolor{mycolor1}{$\bm{\times}$}&
\textcolor{mycolor2}{\Checkmark}&
\textcolor{mycolor2}{\Checkmark}&
\textcolor{mycolor2}{\Checkmark} \\
&    $16^3$ &
\textcolor{mycolor1}{$\bm{\times}$}&
\textcolor{mycolor2}{\Checkmark}&
\textcolor{mycolor2}{\Checkmark}&
\textcolor{mycolor2}{\Checkmark} \\
&    $32^3$ &
\textcolor{mycolor1}{$\bm{\times}$}&
\textcolor{mycolor2}{\Checkmark}&
\textcolor{mycolor2}{\Checkmark}&
\textcolor{mycolor2}{\Checkmark} \\
&    $64^3$ &
\textcolor{mycolor2}{\Checkmark}&
\textcolor{mycolor2}{\Checkmark}&
\textcolor{mycolor2}{\Checkmark}&
\textcolor{mycolor2}{\Checkmark} \\
\hline\hline
DC & $4^3$ &
\textcolor{mycolor1}{$\bm{\times}$}&
\textcolor{mycolor1}{$\bm{\times}$}&
\textcolor{mycolor1}{$\bm{\times}$}&
\textcolor{mycolor1}{$\bm{\times}$} \\
&    $8^3$ &
\textcolor{mycolor1}{$\bm{\times}$}&
\textcolor{mycolor1}{$\bm{\times}$}&
\textcolor{mycolor1}{$\bm{\times}$}&
\textcolor{mycolor1}{$\bm{\times}$} \\
&    $16^3$ &
\textcolor{mycolor1}{$\bm{\times}$}&
\textcolor{mycolor1}{$\bm{\times}$}&
\textcolor{mycolor1}{$\bm{\times}$}&
\textcolor{mycolor2}{\Checkmark} \\
&    $32^3$ &
\textcolor{mycolor1}{$\bm{\times}$}&
\textcolor{mycolor1}{$\bm{\times}$}&
\textcolor{mycolor2}{\Checkmark}&
\textcolor{mycolor2}{\Checkmark} \\
&    $64^3$ &
\textcolor{mycolor2}{\Checkmark}&
\textcolor{mycolor2}{\Checkmark}&
\textcolor{mycolor2}{\Checkmark}&
\textcolor{mycolor2}{\Checkmark} \\
\end{tabular}
\caption{Numerical stability for the compressible homogeneous isotropic turbulence at $\Rey_\lambda=194$ and $\Ma_t=0.62$. \textcolor{mycolor2}{\Checkmark}= success, \textcolor{mycolor1}{$\bm{\times}$} = failure.}
\label{table_check_chit}
\end{table}

As summarized in \Cref{table_check_chit}, only the ES-C scheme shows the numerical robustness to compute this flow for all the selected meshes and solution polynomial degree and produce results past the initial start-up phase.
The conventional DC algorithm is numerically unstable at lower resolution at all polynomial orders, while SF-KG crashes if a polynomial order of $p=1$ is used for a mesh resolution up to $32^3$.

\begin{figure}[!ht]
\centering
\begin{subfigure}{0.47\textwidth}\centering 
\begin{tikzpicture}[scale=1.00]
\begin{axis}[grid=both,minor tick num=3, every major grid/.style={black, opacity=1.0}, minor grid style={gray!40, ultra thin}, legend entries={\scriptsize ES-C,\scriptsize SF-KG,\scriptsize DC},legend style={at={(0.57,0.98)},anchor=north west}, xlabel={$t/\tau$},ylabel={$\mathcal{E}_K/\mathcal{E}_K(t=0)$},y tick label style={/pgf/number format/.cd,fixed,fixed zerofill,precision=1,/tikz/.cd},x tick label style={/pgf/number format/.cd,fixed,fixed zerofill,precision=1,/tikz/.cd},legend columns=1,legend cell align={left},width=7cm, height=6cm,xmin=0,xmax=3.0,ymax=1,ymin=0.1]
\addplot[draw=mycolor1,line width=2.0pt] table[x expr=+\thisrowno{1}/0.030,y expr=+\thisrowno{3}/1.65161750e+00]{data/CHIT/logEC.dat};
\addplot[draw=mycolor2,line width=2.0pt] table[x expr=+\thisrowno{1}/0.030,y expr=+\thisrowno{3}/1.65161750e+00]{data/CHIT/logKG.dat};
\addplot[draw=mycolor11,line width=2.0pt] table[x expr=+\thisrowno{1}/0.030,y expr=+\thisrowno{3}/1.65161750e+00]{data/CHIT/logConv.dat};
\coordinate [](A1) at (0.17625409866666666,0);
\coordinate [](A2) at (0.17625409866666666,1);
\draw [thick,dashed,mycolor2] (A1)--(A2);
\coordinate [](B1) at (0.29340646266666665,0);
\coordinate [](B2) at (0.29340646266666665,1);
\draw [thick,dashed,mycolor11] (B1)--(B2);
\end{axis}
\end{tikzpicture}
\end{subfigure}
\begin{subfigure}{0.47\textwidth}\centering 
\begin{tikzpicture}[scale=1.00]
\begin{axis}[grid=both,minor tick num=3, every major grid/.style={black, opacity=1.0}, minor grid style={gray!40, ultra thin}, legend entries={\scriptsize ES-C,\scriptsize SF-KG,\scriptsize DC},legend style={at={(0.57,0.98)},anchor=north west}, xlabel={$t/\tau$},ylabel={$\Rey_\lambda$},y tick label style={/pgf/number format/.cd,fixed,fixed zerofill,precision=0,/tikz/.cd},x tick label style={/pgf/number format/.cd,fixed,fixed zerofill,precision=1,/tikz/.cd},legend columns=1,legend cell align={left},width=7cm, height=6cm,xmin=0,xmax=3.0,ymin=110,ymax=200]
\addplot[draw=mycolor1,line width=2.0pt] table[x expr=+\thisrowno{1}/0.030,y expr=+\thisrowno{5}/0.558792546875]{data/CHIT/logEC.dat};
\addplot[draw=mycolor2,line width=2.0pt] table[x expr=+\thisrowno{1}/0.030,y expr=+\thisrowno{5}/0.558792546875]{data/CHIT/logKG.dat};
\addplot[draw=mycolor11,line width=2.0pt] table[x expr=+\thisrowno{1}/0.030,y expr=+\thisrowno{5}/0.558792546875]{data/CHIT/logConv.dat};
\coordinate [](A1) at (0.17625409866666666,100);
\coordinate [](A2) at (0.17625409866666666,200);
\draw [thick,dashed,mycolor2] (A1)--(A2);
\coordinate [](B1) at (0.29340646266666665,100);
\coordinate [](B2) at (0.29340646266666665,200);
\draw [thick,dashed,mycolor11] (B1)--(B2);
\end{axis}
\end{tikzpicture}
\end{subfigure}	
\caption{Compressible isotropic turbulence at $\Rey_\lambda=194$ and $\Ma_t=0.62$. The points of termination are indicated by dashed lines.}
\label{fig_dhit_comp_mat}
\end{figure}
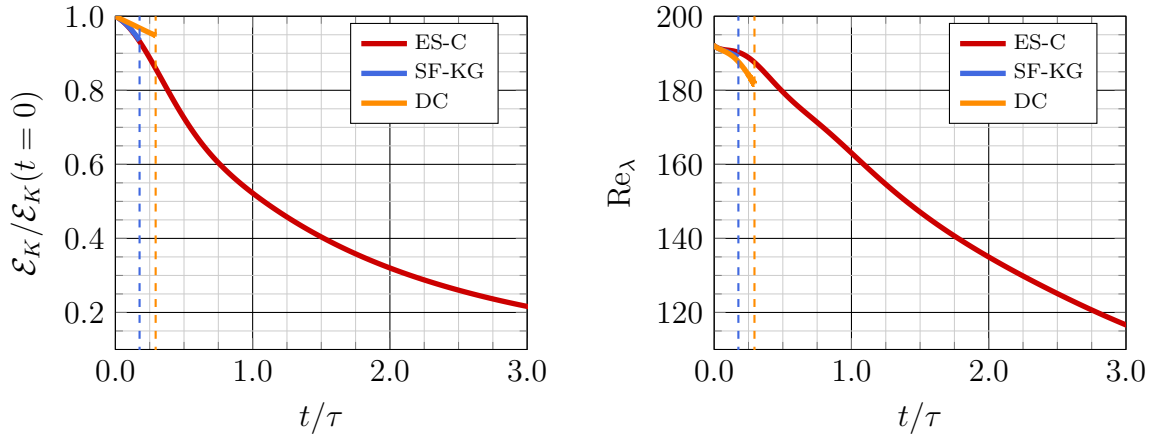

\Cref{fig_dhit_comp_mat} shows both the decaying history of the resolved turbulent kinetic energy $\mathcal{E}_K=\langle\rho\mathcal{U}_i\mathcal{U}_i\rangle/2$ versus $t/\tau$ and Reynolds number based on Taylor micro-scale versus $t/\tau$ for the different schemes and polynomial order $p=1$, where the dashed lines indicate the points of termination (crash).

\begin{figure}[!ht]
\centering
\begin{subfigure}{0.45\textwidth}\centering 
\includegraphics[width=0.80\columnwidth,trim = 0 0 -5 0,clip]{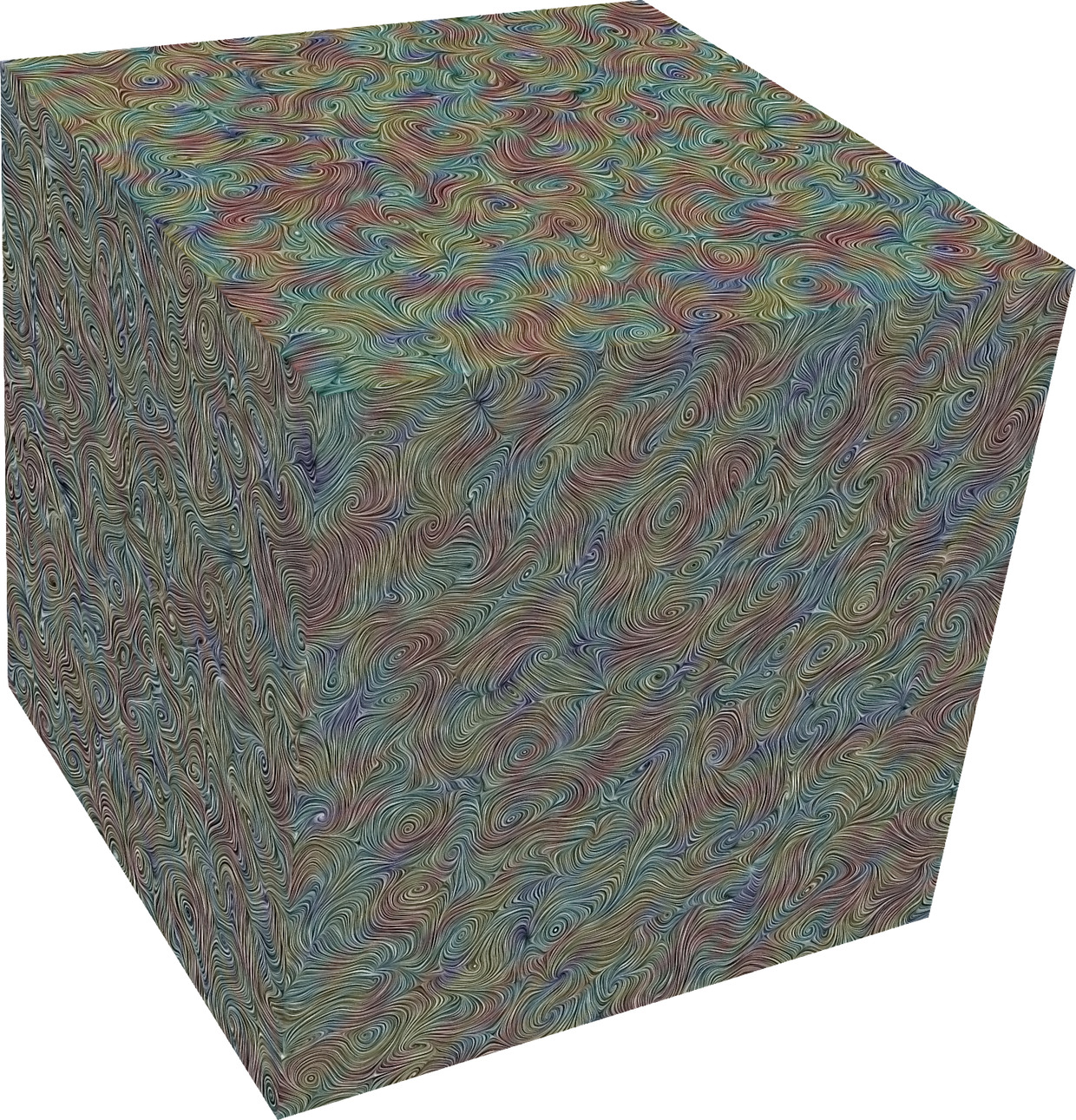}
\caption{$t/\tau=0$.}
\end{subfigure}
\begin{subfigure}{0.45\textwidth}\centering 
\includegraphics[width=0.80\columnwidth,trim = 0 0 -5 0,clip]{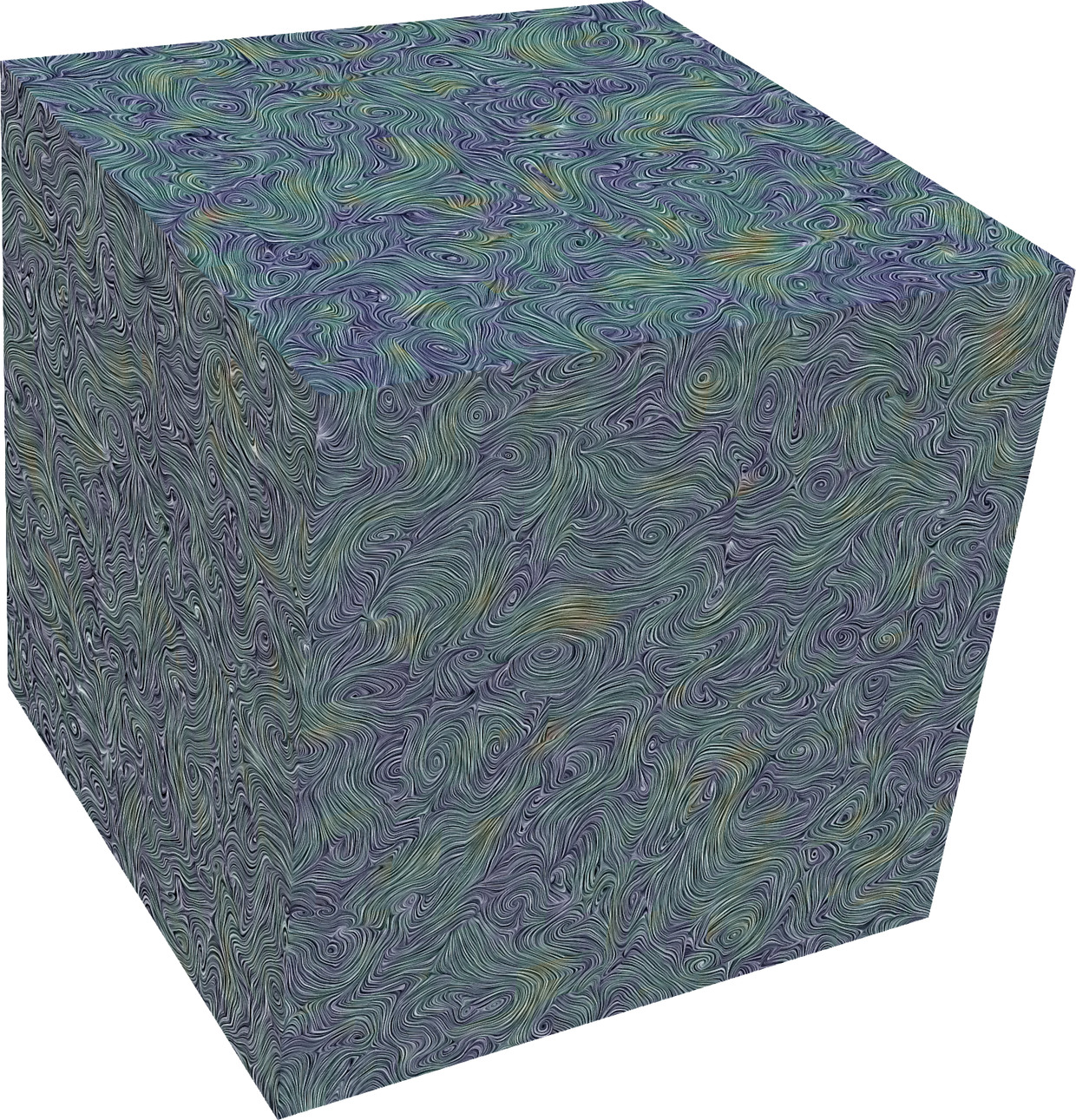}
\caption{$t/\tau=3$.}
\end{subfigure}
\caption{Path lines on the surface of velocity field from ES-C with $p=9$ using $32^3$ for the compressible isotropic turbulence at $\Rey_\lambda=194$ and $\Ma_t=0.62$.}
\label{fig-dhit_mt03_pathlines}
\end{figure}

In \Cref{fig-dhit_mt03_pathlines}, the path lines on the boundaries are plotted at two instants for the ES-C scheme using the line integral convolution colored by the magnitude of the velocity.
As expected, the velocity field at $t/\tau=3$ clearly contains more small scales than the initial field $t/\tau=2$, which shows the need to use a numerically robust scheme to capture them without filtering and also to properly handle shocklets that quickly form for such a higher turbulent Mach number. 
It is concluded from this test case that numerical schemes which are entropy stable
are more desirable for such a highly strong Mach-number regime. 
In fact, schemes that satisfy entropy condition are found to lead to stable density fluctuations in compressible isotropic turbulence simulations, while schemes that do not have this property can be unstable with respect to density fluctuations \cite{honein2004higher,pirozzoli2004direct}.

%%%%%%%%%%%%%%%%%%%%%%%%%%%%%%%%%%%%%%%%%%%%%%%%%%%%%%%%%%%%%%%%%%%%%%%%%%%%%%%%%%%%%%%
\subsubsection{Supersonic turbulent flows past a bar}

The last test we present is the 3D supersonic flow around a square bar, a problem that has been studied in depth in the literature related to supersonic flow around bluff bodies \cite{nakagawa1988effects,birch2003experimental,parsani_entropy_stability_solid_wall_2015}.
The similarity parameters are $\Rey_{\infty}=10{,}000$ and $\Ma_{\infty}=1.5$. 
In this regime the flow exhibits a complicated pattern characterized by shocks, expansion zones and separation \cite{parsani_entropy_stability_solid_wall_2015}.
This test is chosen as a challenging problem to stress the robustness of our discretizations in the context of supersonic turbulent flow.

\begin{figure}[!ht]
\centering
\includegraphics[width=0.45\textwidth]{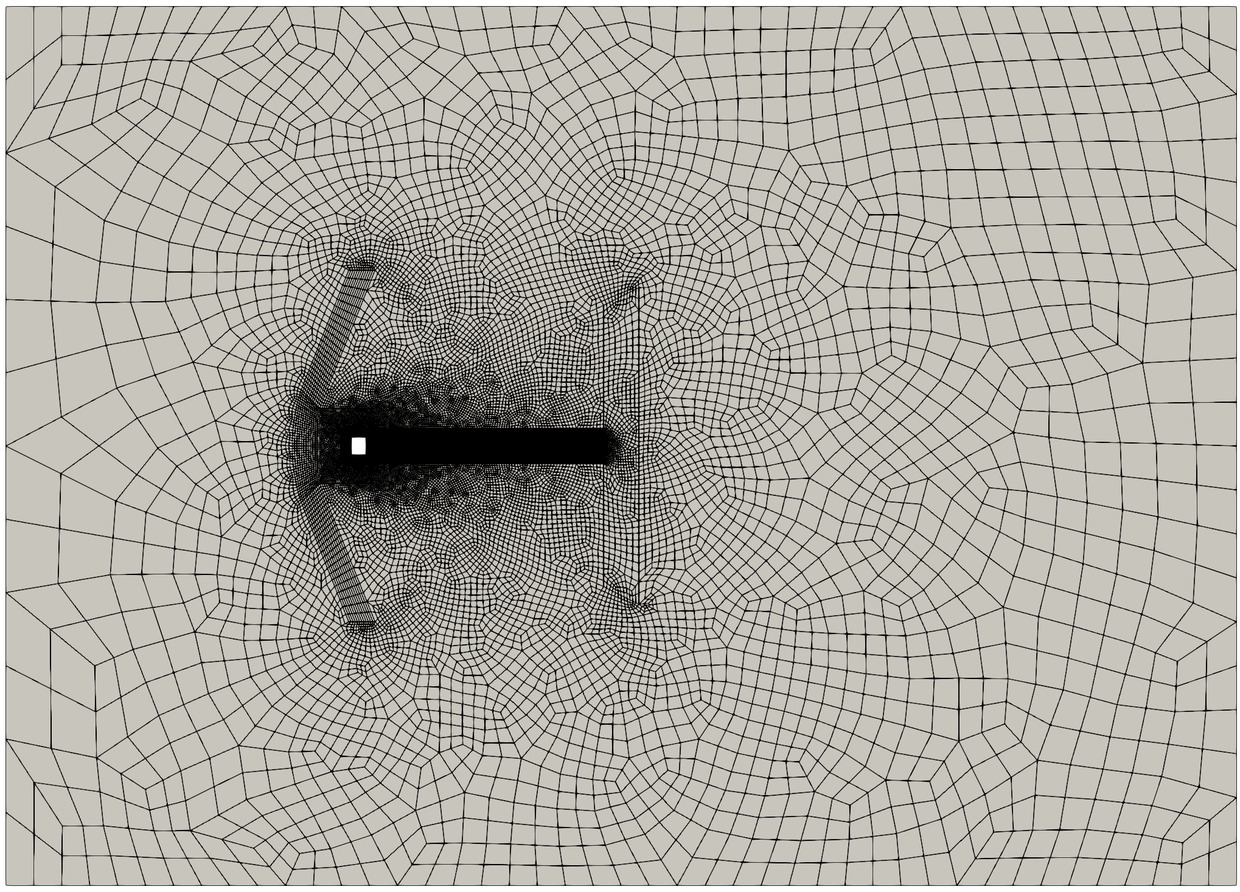}
\includegraphics[width=0.40\textwidth]{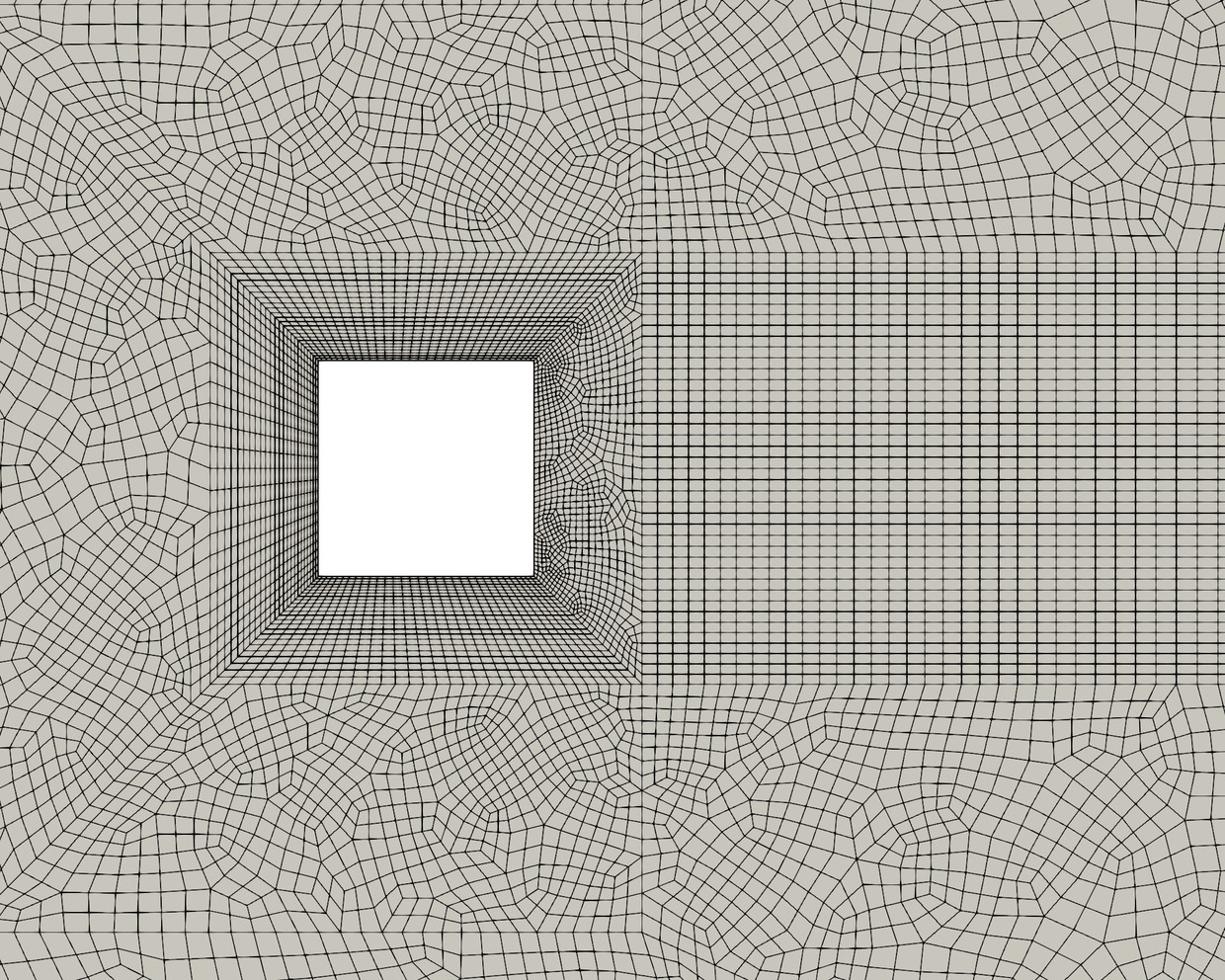}

\caption{Bi-dimensional view of (left) the 3D mesh used for the supersonic flow around a square bar and (right) zoom at the square bar.}
\label{fig_ssbar_mesh}
\end{figure}

A 3D square of unit side is placed at the origin of the system of reference and an unstructured quad mesh with refinement is extended out to a box $20$ and $50$ length units in the upstream and downstream directions, respectively, 25 length units in both normal directions to the bar, and then extruded one length unit in the spanwise direction.
Figure \ref{fig_ssbar_mesh} shows the resulting mesh with refinements around the shocks and the turbulent wake; we choose a polynomial degree of $p=3$ as a reasonable accuracy for this mesh.
The boundary conditions used are adiabatic wall on the surface of the bar \cite{parsani_entropy_stability_solid_wall_2015,dalcin2019conservative}, periodic in the spanwise direction and farfield on the remaining boundaries.
The initial condition is set to be a uniform supersonic flow. 
The simulation run until the main bow shock stabilizes.

\begin{figure}[!ht]
\centering
\begin{subfigure}{0.47\textwidth}\centering
\begin{tikzpicture}[scale=1.00]
\begin{axis}[grid=both,minor tick num=3, every major grid/.style={black, opacity=1.0}, minor grid style={gray!40, ultra thin}, legend entries={\scriptsize ES-C,\scriptsize SF-KG,\scriptsize DC},legend style={at={(0.57,0.98)},anchor=north west}, xlabel={$t$},ylabel={$C_d^{(\mathcal{P})}$},y tick label style={/pgf/number format/.cd,fixed,fixed zerofill,precision=1,/tikz/.cd},x tick label style={/pgf/number format/.cd,fixed,fixed zerofill,precision=1,/tikz/.cd},legend columns=1,legend cell align={left},width=7cm, height=6cm,xmin=0,xmax=1.0,ymax=7,ymin=1.]
\addplot[draw=mycolor1,line width=1.0pt] table[x expr=+\thisrowno{1},y expr=+\thisrowno{3}*2]{data/SUPCYL/es-forces.txt};
\addplot[draw=mycolor2,line width=1.0pt] table[x expr=+\thisrowno{1},y expr=+\thisrowno{3}*2]{data/SUPCYL/kg-forces.txt};
\coordinate [](A1) at (0.36196890593,0);
\coordinate [](A2) at (0.36196890593,10);
\draw [thick,dashed,mycolor2] (A1)--(A2);
\end{axis}
\end{tikzpicture}
\end{subfigure}
\begin{subfigure}{0.47\textwidth}\centering
\begin{tikzpicture}[scale=1.00]
\begin{axis}[grid=both,minor tick num=3, every major grid/.style={black, opacity=1.0}, minor grid style={gray!40, ultra thin}, legend entries={\scriptsize ES-C,\scriptsize SF-KG,\scriptsize DC },legend style={at={(0.57,0.98)},anchor=north west}, xlabel={$t$},ylabel={$C_d^{(\mathrm{V})}$},y tick label style={/pgf/number format/.cd,fixed,fixed zerofill,precision=1,/tikz/.cd},x tick label style={/pgf/number format/.cd,fixed,fixed zerofill,precision=1,/tikz/.cd},legend columns=1,legend cell align={left},width=7cm, height=6cm,xmin=0,xmax=1.0,ymin=0.,ymax=0.3]
\addplot[draw=mycolor1,line width=2.0pt] table[x expr=+\thisrowno{1},y expr=+\thisrowno{5}*2]{data/SUPCYL/es-forces.txt};
\addplot[draw=mycolor2,line width=2.0pt] table[x expr=+\thisrowno{1},y expr=+\thisrowno{5}*2]{data/SUPCYL/kg-forces.txt};
\coordinate [](A1) at (0.36196890593,0);
\coordinate [](A2) at (0.36196890593,1);
\draw [thick,dashed,mycolor2] (A1)--(A2);
\end{axis}
\end{tikzpicture}
\end{subfigure}
\caption{Drag coefficient for the supersonic bar at $\Rey_{\infty}=10{,}000$ and $\Ma_{\infty}=1.5$.}
\label{fig:ssbar_c_d}
\end{figure}
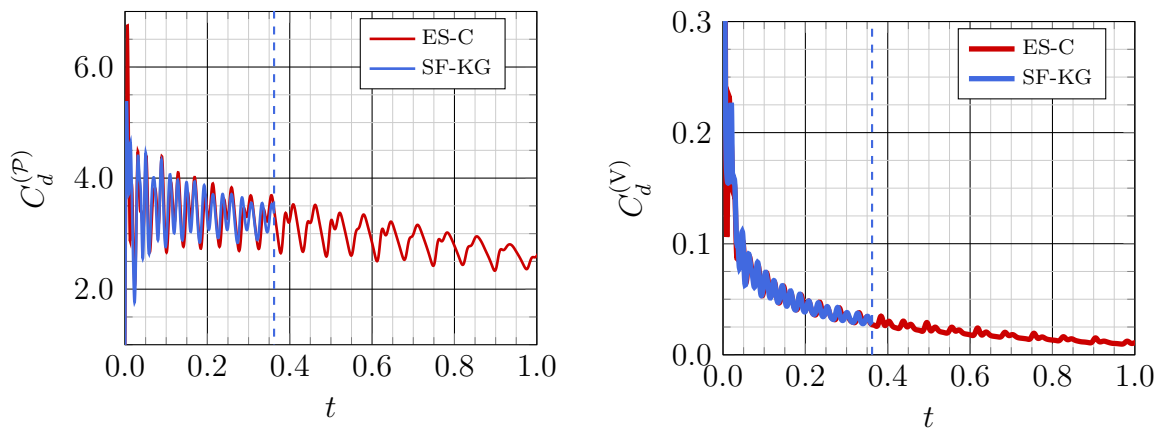

\Cref{fig:ssbar_c_d} shows the results for the pressure drag coefficient ($C_d^{(\mathcal{P})}$) and the viscous drag coefficient ($C_d^{(\mathrm{V})}$) for the initial time steps for the ES-C and the SF-KG discretizations. 
The dashed line indicates the termination time due to the solver crash.
Here the the DC discretization crashed very early in the simulation ($t\approx1.3\times10^{-3}$) yielding no useful results.

\begin{figure}
\centering
\begin{subfigure}{0.7\textwidth}\centering
\includegraphics[width=0.9\columnwidth]{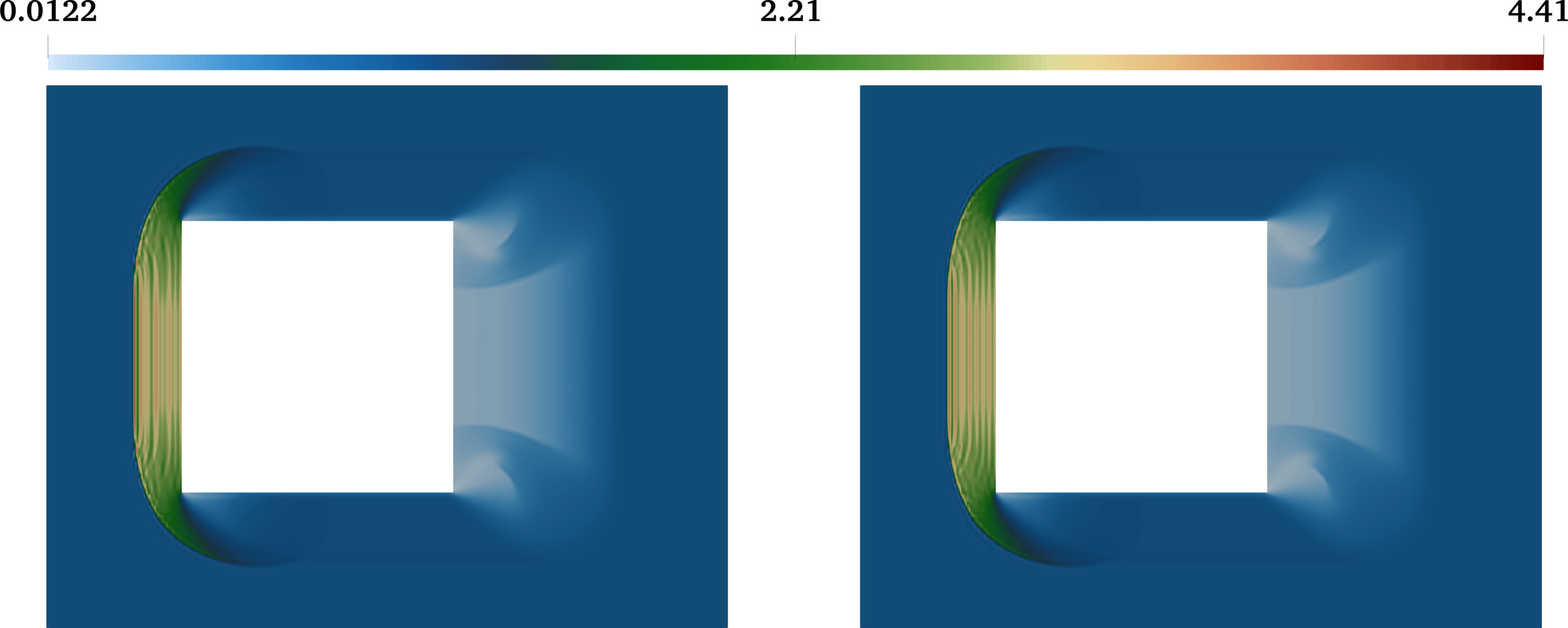}
\caption{Density $\rho$.}
\end{subfigure}
\begin{subfigure}{0.7\textwidth}\centering
\includegraphics[width=0.9\columnwidth]{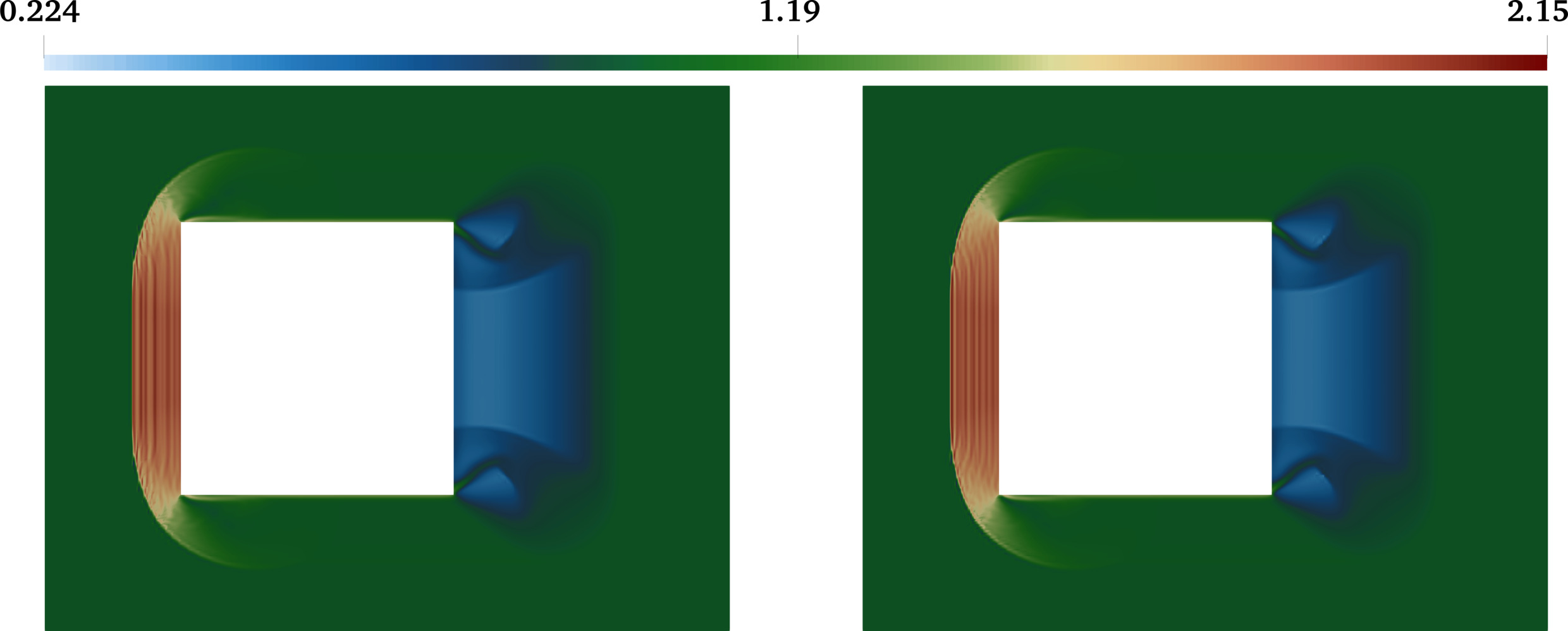}
\caption{Temperature $\mathcal{T}$.}
\end{subfigure}
\begin{subfigure}{0.7\textwidth}\centering
\includegraphics[width=0.9\columnwidth]{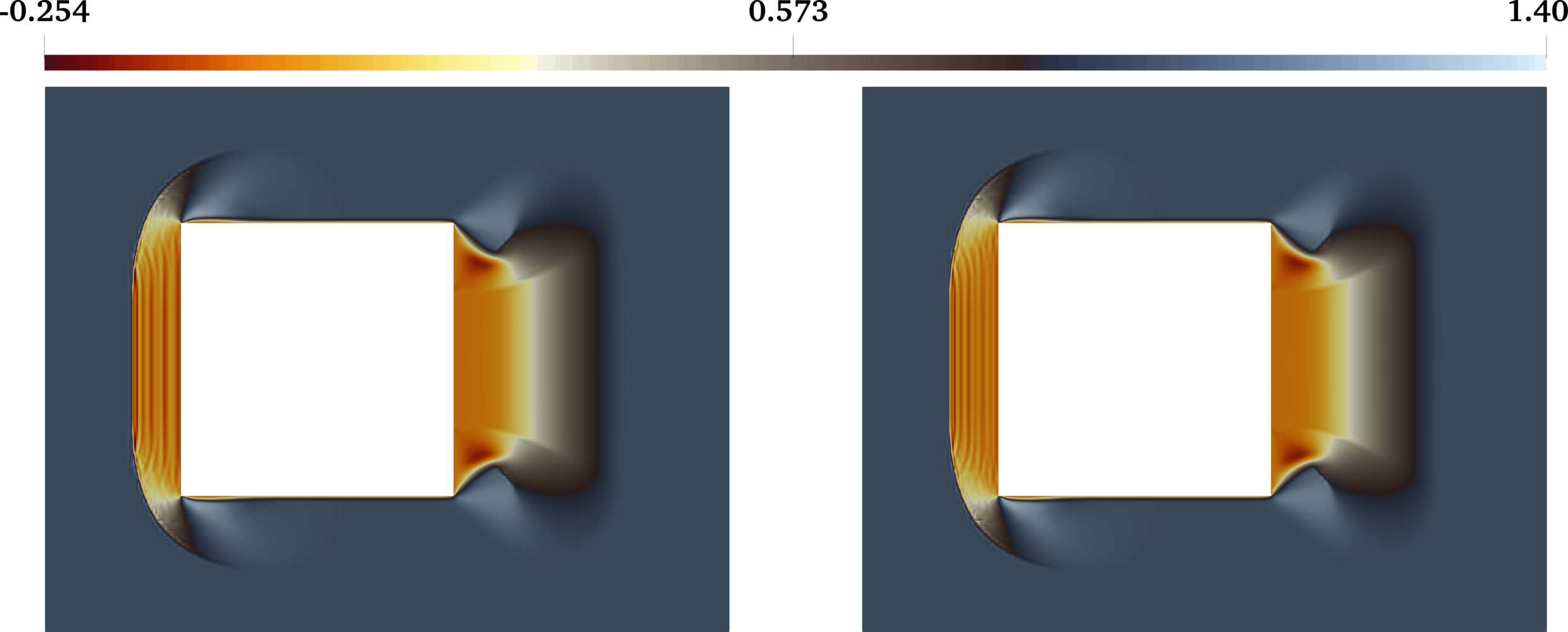}
\caption{Velocity component $\mathcal{U}_1$.}
\end{subfigure}
\begin{subfigure}{0.7\textwidth}\centering
\includegraphics[width=0.9\columnwidth]{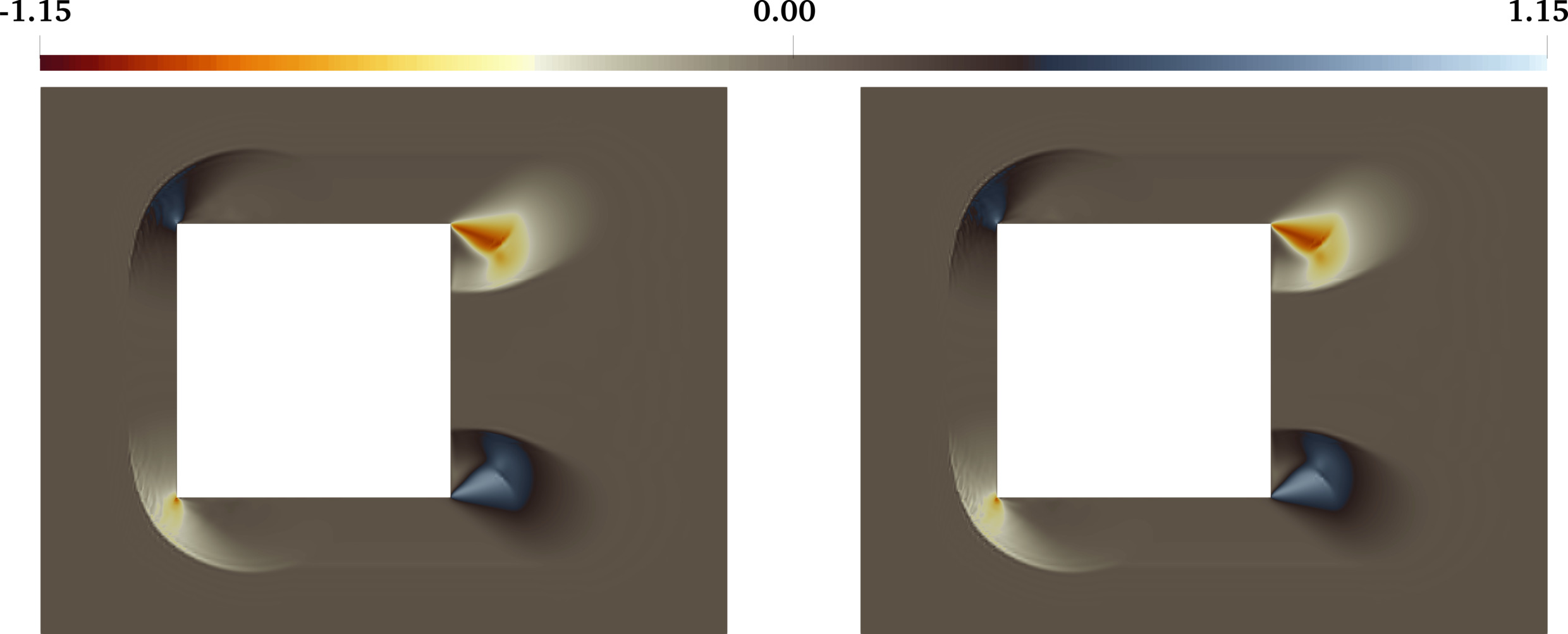}
\caption{Velocity component $\mathcal{U}_2$.}
\end{subfigure}
\caption{Comparison of the last recorded step for SF-KG (right) and ES-C (left) at $t=0.36$ for the supersonic bar at $\Rey_{\infty}=10{,}000$ and $\Ma_{\infty}=1.5$.}
\label{fig_comp_ssbar}
\end{figure}

\Cref{fig_comp_ssbar} shows a visualization of the last recorded step of the solution given by SF-KG (right side), compared with ES-C (left) at the same time.
This result is very early in the simulation when the shock propagating from the leading edge is still close to the wall.
Here we choose color schemes intended to highlight the differences between the two simulations, particularly in the shock region where the (expected) oscillations appear to be more regular for the ES-C case.
We remark that this visualization has no evidence for the onset of numerical instability; even though SF-KG immediately crashes at the next time step.

\begin{figure}[!ht]
\centering
\begin{subfigure}{0.45\textwidth}\centering
\includegraphics[width=0.9\columnwidth,trim = 0 0 -5 0,clip]{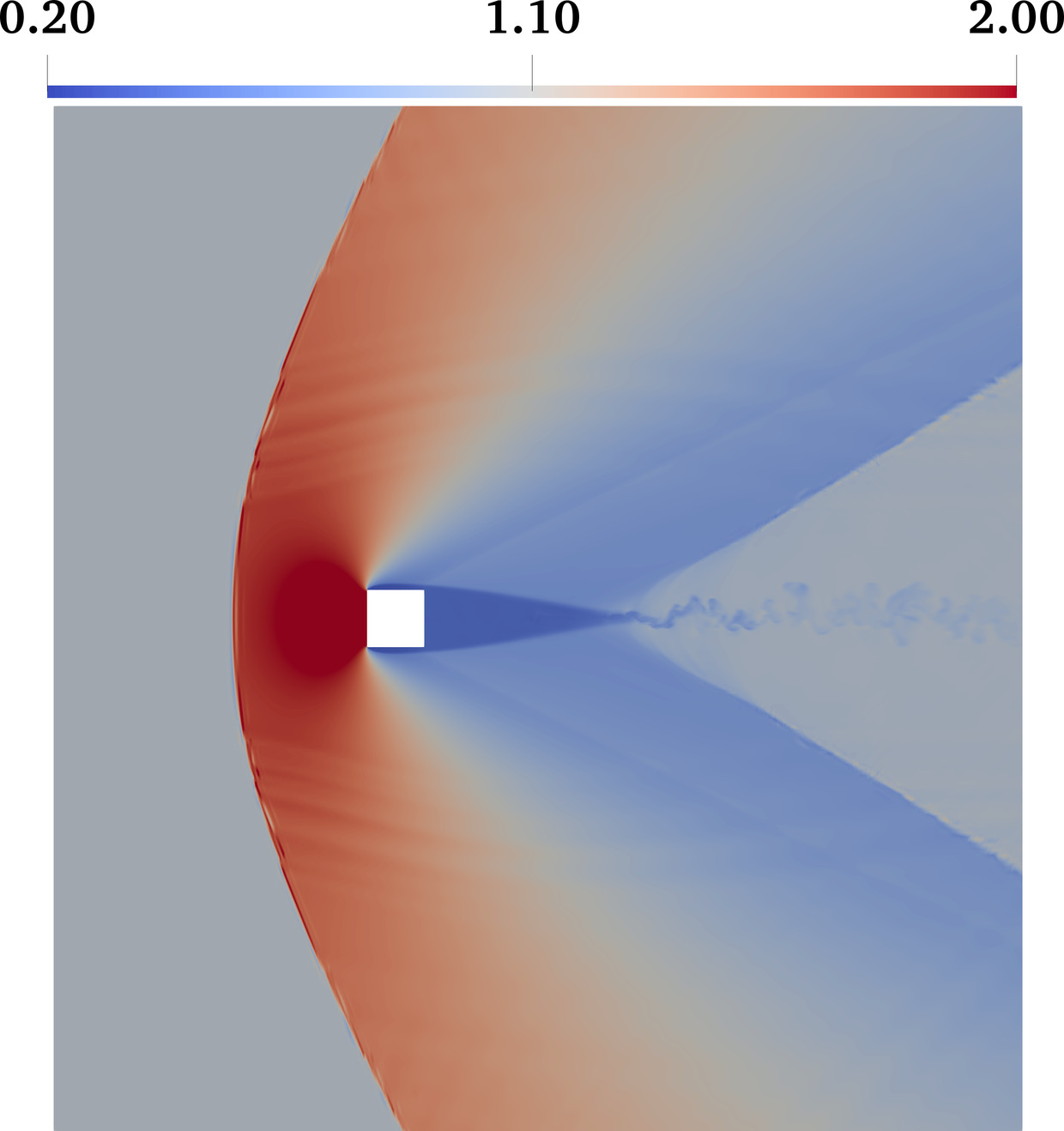}
\caption{Density $\rho$.}
\end{subfigure}
\begin{subfigure}{0.45\textwidth}\centering
\includegraphics[width=0.9\columnwidth,trim = 0 0 -5 0,clip]{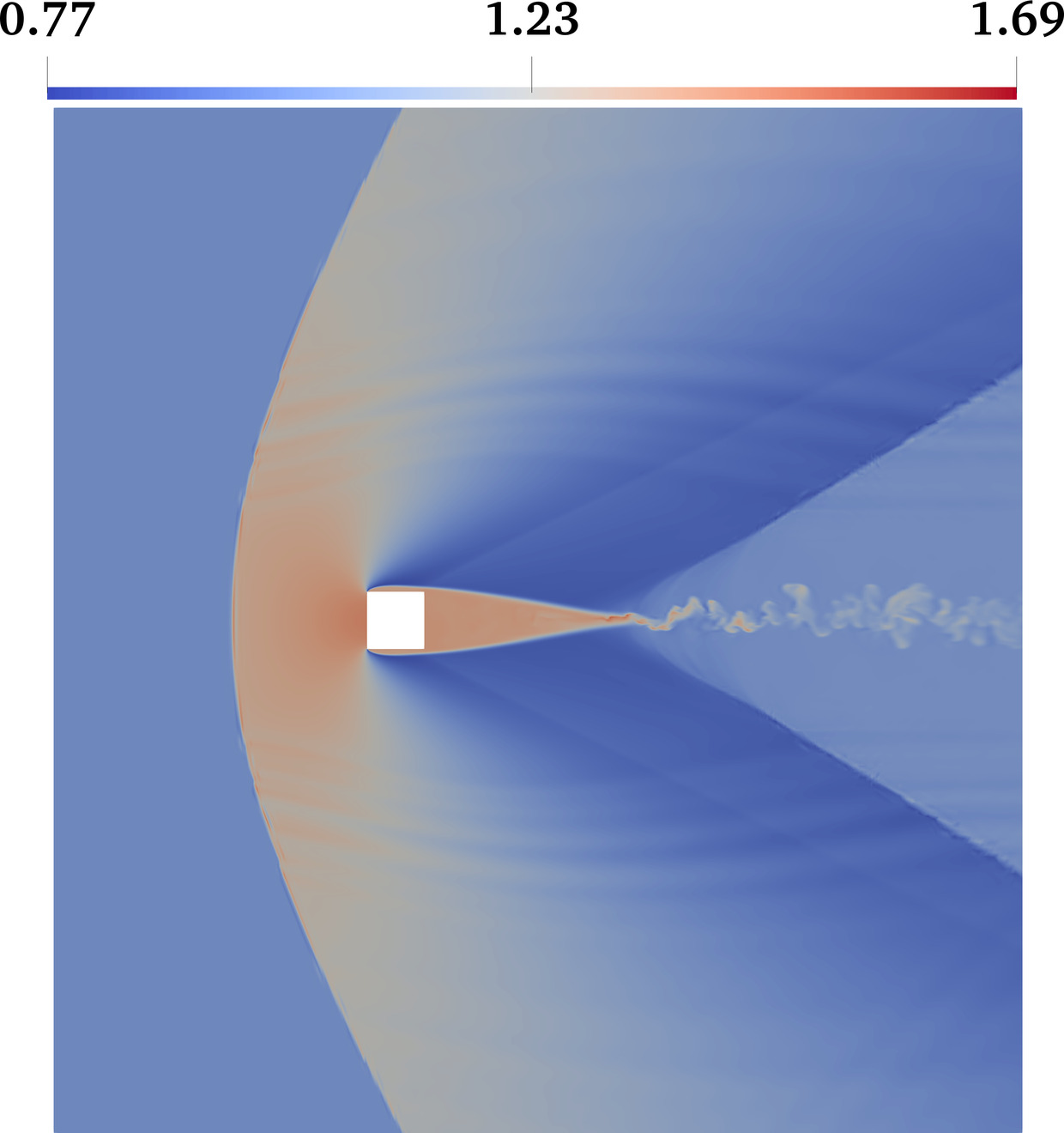}
\caption{Temperature $\mathcal{T}$.}
\end{subfigure}\\
\begin{subfigure}{0.45\textwidth}\centering
\includegraphics[width=0.9\columnwidth,trim = 0 0 -5 0,clip]{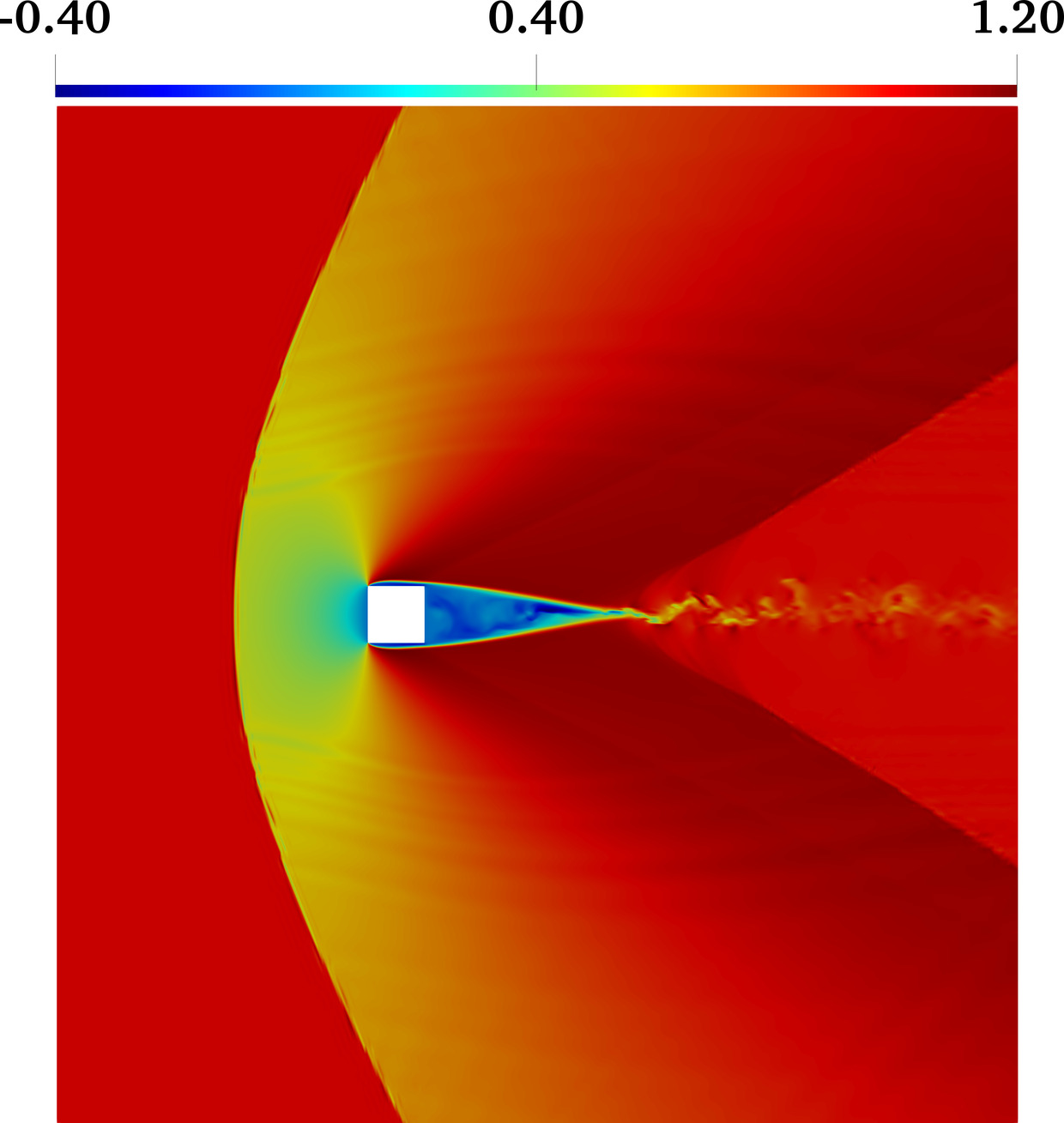}
\caption{Velocity component $\mathcal{U}_1$.}
\end{subfigure}
\begin{subfigure}{0.45\textwidth}\centering
\includegraphics[width=0.9\columnwidth,trim = 0 0 -5 0,clip]{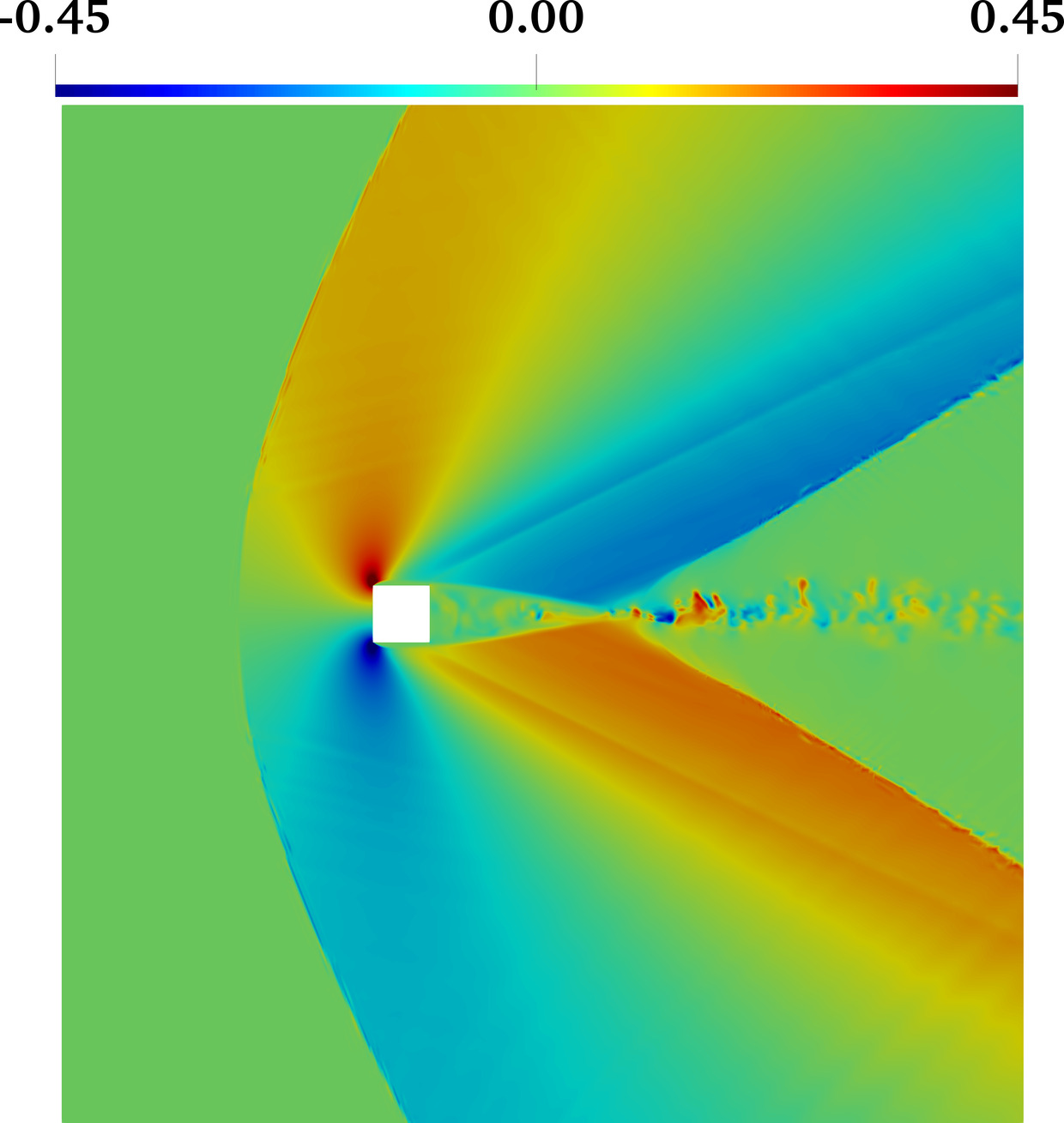}
\caption{Velocity component $\mathcal{U}_2$.}
\end{subfigure}
\caption{Instantaneous snapshots of fields at $t=130$ for the supersonic bar at $\Rey_{\infty}=10{,}000$ and $\Ma_{\infty}=1.5$.}
\label{fig_ssbar}
\end{figure}

\Cref{fig_ssbar} shows the resulting flow pattern for the ES-C discretization at $t=100$. 

%%%%%%%%%%%%%%%%%%%%%%%%%%%%%%%%%%%%%%%%%%%%%%%%%%%%%%%%%%%%%%%%%%%%%%%%%%%%%%%%%%%%%%%
\section{Conclusion}\label{sec:conclusion}

This work presents a three-way comparison between one entropy stable, and two non-entropy stable discretizations, in the context of split form and discontinuous collocation discretizations for the simulation of the compressible Navier--Stokes equations.
All three algorithms are constructed using SBP-SAT operators. 
The entropy stable discretization is based on the two-point entropy conservative flux of \citet{chandrashekar2013kinetic}, the kinetic energy preserving but not entropy stable discretization is built using the splitting of \citet{kennedy2008reduced}, and the conventional discontinuous collocation method is based on the standard inviscid flux arising from the compressible Euler equations. 

We have started our discussion with a note on the potential superiority of high order discretizations, not only measured in terms of accuracy but also in wall-clock time for the propagation of a 3D isentropic vortex and a 3D supersonic flow constructed with the method of manufactured solutions. 
Subsequently, by using the Taylor--Green vortex problem, we have shown the impact of the solution polynomial degree on the computational cost of the three discretizations.
In DC-type discretizations, the computational cost depends on the number of flux evaluations (multiplied by its computational complexity), and the cost of computing the surface terms (SATs).
These costs have opposite behaviors as the polynomial degree increases and therefore, an optimal wall-clock time per time step can be found.
The overall observation is that higher degrees ($p=3$ or higher) yield the lowest wall-clock time per time step.
This superiority is of paramount importance in the implementation of such simulation tools for current and future computing hardware. Hence, the emphasis is placed on higher-order discretizations throughout the entire paper.

We then moved on to the analysis of the reliability and robustness of the three discretizations under three stressful experimental conditions: i) turbulent flows with under-resolved features, ii) shocklets in compressible turbulent flow, and iii) shocks in supersonic flow.
As expected, only the entropy stable discretization is able to simulate all test cases successfully.
This robustness comes at a cost, which we consider to be affordable for reasonably high-order discretizations ($p$ up to 7), especially in large-scale complex simulations where any instability is unacceptable.
We interpret these results as significant evidence in favor of the applicability of entropy stable discretizations in future CFD solvers as reported in \cite{nasa_2030_vision}. 

We expect to see future extensions to this work addressing more applied aspects of the implementation of entropy stable discretizations such as performance at scale and $hp$-refinement, to name a few.

%%%%%%%%%%%%%%%%%%%%%%%%%%%%%%%%%%%%%%%%%%%%%%%%%%%%%%%%%%%%%%%%%%%%%%%%%%%%%%%%
\section*{Acknowledgments}
The research reported in this paper was funded by King Abdullah University of Science and Technology. 
We are thankful for the computing resources of the Supercomputing Laboratory and the Extreme Computing Research Center at King Abdullah University of Science and Technology.

%%%%%%%%%%%%%%%%%%%%%%%%%%%%%%%%%%%%%%%%%%%%%%%%%%%%%%%%%%%%%%%%%%%%%%%%%%%%%%%%
\section*{References}
\bibliography{fixed-comparison}

\begin{thebibliography}{}

\bibitem[Abhyankar et~al., 2018]{abhyankar2018petsc}
Abhyankar, S., Brown, J., Constantinescu, E.~M., Ghosh, D., Smith, B.~F., and
  Zhang, H. (2018).
\newblock {PETSc/TS}: {A} modern scalable {ODE}/{DAE} solver library.

\bibitem[Balay et~al., 2018]{petsc-user-ref}
Balay, S., Abhyankar, S., Adams, M.~F., Brown, J., Brune, P., Buschelman, K.,
  Dalcin, L., Dener, A., Eijkhout, V., Gropp, W.~D., Kaushik, D., Knepley,
  M.~G., May, D.~A., McInnes, L.~C., Mills, R.~T., Munson, T., Rupp, K., Sanan,
  P., Smith, B.~F., Zampini, S., Zhang, H., and Zhang, H. (2018).
\newblock {PETS}c users manual.
\newblock Technical Report ANL-95/11 - Revision 3.10, Argonne National
  Laboratory.

\bibitem[Birch et~al., 2003]{birch2003experimental}
Birch, T.~J., Prince, S.~A., and Simpson, G.~M. (2003).
\newblock An experimental and computational study of the aerodynamics of a
  square cross-section body at supersonic speeds.
\newblock Technical report, Defence Evaluation and Research Agency.

\bibitem[Boukharfane, 2018]{boukharfane2018contribution}
Boukharfane, R. (2018).
\newblock {\em Contribution {\`a} la simulation num{\'e}rique d'{\'e}coulements
  turbulents compressibles canoniques}.
\newblock PhD thesis, Chasseneuil-du-Poitou, {\'E}cole Nationale Sup{\'e}rieure
  de M{\'e}canique et d'A{\'e}rotechnique.

\bibitem[Carpenter et~al., 2014]{carpenter_ssdc_2014}
Carpenter, M.~H., Fisher, T., Nielsen, E., and Frankel, S. (2014).
\newblock Entropy stable spectral collocation schemes for the
  {N}avier--{S}tokes equations: Discontinuous interfaces.
\newblock {\em SIAM Journal on Scientific Computing}, 36(5):B835--B867.

\bibitem[Carpenter et~al., 2015]{carpenter_entropy_stable_staggered_2015}
Carpenter, M.~H., Parsani, M., Fisher, T.~C., and Nielsen, E.~J. (2015).
\newblock Entropy stable staggered grid spectral collocation for the {B}urgers'
  and the compressible {N}avier--{S}tokes equations.
\newblock {\em NASA TM-2015-218990}.

\bibitem[Carpenter et~al., 2016]{carpenter_entropy_stability_ssdc_2016}
Carpenter, M.~H., Parsani, M., Nielsen, E.~J., and Fisher, T.~C. (2016).
\newblock Towards an entropy stable spectral element framework for
  computational fluid dynamics.
\newblock In {\em 54th AIAA Aerospace Sciences Meeting}, AIAA 2016-1058.
  American Institute of Aeronautics and Astronautics.

\bibitem[Chan, 2018]{Chan2018}
Chan, J. (2018).
\newblock On discretely entropy conservative and entropy stable discontinuous
  {G}alerkin methods.
\newblock {\em Journal of Computational Physics}, 362:346 -- 374.

\bibitem[Chandrashekar, 2013]{chandrashekar2013kinetic}
Chandrashekar, P. (2013).
\newblock Kinetic energy preserving and entropy stable finite volume schemes
  for compressible euler and navier-stokes equations.
\newblock {\em Communications in Computational Physics}, 14(5):1252--1286.

\bibitem[Chen and Shu, 2017]{Chen2017}
Chen, T. and Shu, C.-W. (2017).
\newblock Entropy stable high order discontinuous {G}alerkin methods with
  suitable quadrature rules for hyperbolic conservation laws.
\newblock {\em Journal of Computational Physics}, 345:427 -- 461.

\bibitem[Crean et~al., 2018]{crean_entropy_stable_sbp_curvilinear_euler}
Crean, J., Hicken, J.~E., Del Rey~Fern\'andez, D.~C., Zingg, D.~W., and
  Carpenter, M.~H. (2018).
\newblock Entropy-stable summation-by-parts discretization of the {E}uler
  equations on general curved elements.
\newblock {\em Journal of Computational Physics}, 356:410 -- 438.

\bibitem[Dafermos, 2010]{dafermos_book_2010}
Dafermos, C.~M. (2010).
\newblock {\em Hyperbolic conservation laws in continuum physics}.
\newblock Springer-Verlag, Berlin.

\bibitem[Dalcin et~al., 2019]{dalcin2019conservative}
Dalcin, L., Rojas, D., Zampini, S., Del Rey~Fern\'andez, D.~C., Carpenter,
  M.~H., and Parsani, M. (2019).
\newblock Conservative and entropy stable solid wall boundary conditions for
  the compressible {Navier--Stokes} equations: {A}diabatic wall and heat
  entropy transfer.
\newblock {\em Journal of Computational Physics}, 397:108775.

\bibitem[DeBonis, 2013]{debonis2013solutions}
DeBonis, J. (2013).
\newblock Solutions of the taylor-green vortex problem using high-resolution
  explicit finite difference methods.
\newblock In {\em 51st AIAA Aerospace Sciences Meeting including the New
  Horizons Forum and Aerospace Exposition}, page 382.

\bibitem[Del Rey~Fern\'andez et~al., 2014a]{DCDRF2014}
Del Rey~Fern\'andez, D.~C., Boom, P.~D., and Zingg, D.~W. (2014a).
\newblock A generalized framework for nodal first derivative summation-by-parts
  operators.
\newblock {\em Journal of Computational Physics}, 266(1):214--239.

\bibitem[Del Rey~Fern\'andez et~al.,
  2019b]{fernandez_entropy_stable_p_euler_2019}
Del Rey~Fern\'andez, D.~C., Carpenter, M.~H., Dalcin, L., Fredrich, L.,
  Winters, A.~R., Gassner, G.~J., Zampini, S., and Parsani, M. (2019b).
\newblock Entropy stable $p-$nonconforming discretizations with the
  summation-by-parts property for the compressible {E}uler equations.
\newblock {\em Submitted to SIAM Journal of Scientific Computing}.

\bibitem[Del Rey~Fern\'andez et~al., 2019a]{fernandez_entropy_stable_p_ns_2019}
Del Rey~Fern\'andez, D.~C., Carpenter, M.~H., Dalcin, L., Fredrich, L.,
  Winters, A.~R., Gassner, G.~J., Zampini, S., and Parsani, M. (2019a).
\newblock Entropy stable $p-$nonconforming discretizations with the
  summation-by-parts property for the compressible {N}avier--{S}tokes
  equations.
\newblock {\em Submitted to Computer \& Fluids}.

\bibitem[Del Rey~Fern\'andez et~al.,
  2019c]{fernandez_entropy_stable_hp_ref_snpdea_2019}
Del Rey~Fern\'andez, D.~C., Carpenter, M.~H., Dalcin, L., Zampini, S., and
  Parsani, M. (2019c).
\newblock Entropy stable $h/p$ non-conforming discretization with the
  summation-by-parts property for the compressible {E}uler and {Navier--Stokes}
  equations.
\newblock Submitted to SN Partial Differential Equations and Applications.

\bibitem[Del Rey~Fern\'andez et~al., 2019d]{Fernandez2019_staggered}
Del Rey~Fern\'andez, D.~C., Crean, J., Carpenter, M.~H., and Hicken, J.~E.
  (2019d).
\newblock Staggered entropy-stable summation-by-parts discretization of the
  {E}uler equations on general curved elements.
\newblock {\em Journal of Computational Physics}, 392:161--186.

\bibitem[Del Rey~Fern\'andez et~al., 2014b]{Fernandez2014}
Del Rey~Fern\'andez, D.~C., Hicken, J.~E., and Zingg, D.~W. (2014b).
\newblock Review of summation-by-parts operators with simultaneous
  approximation terms for the numerical solution of partial differential
  equations.
\newblock {\em Computers \& Fluids}, 95(22):171--196.

\bibitem[Dormand and Prince, 1980]{DORMAND198019}
Dormand, J. and Prince, P. (1980).
\newblock A family of embedded runge-kutta formulae.
\newblock {\em Journal of Computational and Applied Mathematics}, 6(1):19 --
  26.

\bibitem[Fern\'{a}ndez et~al., 2019]{fernandez_entropy_stable_p_ref_nasa_2019}
Fern\'{a}ndez, D.~C., Carpenter, M.~H., Dalcin, L., Fredrich, L., Rojas, D.,
  Winters, A.~R., Gassner, G.~J., Zampini, S., and Parsani, M. (2019).
\newblock Entropy stable non-conforming discretizations with the
  summation-by-parts property for curvilinear coordinates.
\newblock {\em NASA TM-2019-}.

\bibitem[Fernandez et~al., 2019]{fernandez2019entropy}
Fernandez, D.~C., Carpenter, M.~H., Dalcin, L., Fredrich, L., Rojas, D.,
  Winters, A.~R., Gassner, G.~J., Zampini, S., and Parsani, M. (2019).
\newblock Entropy stable {$p$}-nonconforming discretizations with the
  summation-by-parts property for the compressible {E}uler equations.

\bibitem[Fisher, 2012]{fisher_phd_2012}
Fisher, T.~C. (2012).
\newblock {\em High-order {$L^{2}$} stable multi-domain finite difference
  method for compressible flows}.
\newblock PhD thesis, Purdue University.

\bibitem[Fisher and Carpenter, 2013]{FisherCarpenter2013JCPb}
Fisher, T.~C. and Carpenter, M.~H. (2013).
\newblock High-order entropy stable finite difference schemes for nonlinear
  conservation laws: Finite domains.
\newblock {\em Journal of Computational Physics}, 252:518--557.

\bibitem[Fisher et~al., 2013]{Fisher2013}
Fisher, T.~C., Carpenter, M.~H., Nordstr\"om, J., and Yamaleev, N.~K. (2013).
\newblock Discretely conservative finite-difference formulations for nonlinear
  conservation laws in split form: Theory and boundary conditions.
\newblock {\em Journal of Computational Physics}, 234(1):353--375.

\bibitem[Fjordholm et~al., 2012]{Fjordholm2012}
Fjordholm, U.~S., Mishra, S., and Tadmor, E. (2012).
\newblock Arbitrarily high-order accurate entropy stable essentially
  nonoscillatory schemes for systems of conservation laws.
\newblock {\em Communications in Computational Physics}, 50(2):554--573.

\bibitem[Flad and Gassner, 2017]{flad2017use}
Flad, D. and Gassner, G.~J. (2017).
\newblock On the use of kinetic energy preserving {DG}-schemes for large eddy
  simulation.
\newblock {\em Journal of Computational Physics}, 350:782--795.

\bibitem[Friedrich et~al., 2019]{Friedrich2019}
Friedrich, L., Shn\"ucke, G., Winters, A.~R., Del Rey~Fern\'andez, D.~C.,
  Gassner, G.~J., and Carpenter, M.~H. (2019).
\newblock Entropy stable space-time discontinuous {G}alerkin schemes with
  summation-by-parts property for hyperbolic conservation laws.
\newblock {\em Journal of Scientific Computing}, 80(1):175--222.

\bibitem[Friedrich et~al., 2018]{friedrich_hp_entropy_stability_2018}
Friedrich, L., Winters, A.~R., Del Rey~Fern\'{a}ndez, D.~C., Gassner, G.~J.,
  Parsani, M., and Carpenter, M.~H. (2018).
\newblock An entropy stable h/p non-conforming discontinuous {G}alerkin method
  with the summation-by-parts property.
\newblock {\em Journal of Scientific Computing}, 77(2).

\bibitem[Gassner, 2013]{gassner2013skew}
Gassner, G.~J. (2013).
\newblock A skew-symmetric discontinuous {G}alerkin spectral element
  discretization and its relation to {SBP}-{SAT} finite difference methods.
\newblock {\em SIAM Journal on Scientific Computing}, 35(3):A1233--A1253.

\bibitem[Gassner and Beck, 2013]{gassner_underresolved_turbulence_2013}
Gassner, G.~J. and Beck, A.~D. (2013).
\newblock On the accuracy of high-order discretizations for underresolved
  turbulence simulations.
\newblock {\em Theoretical and Computational Fluid Dynamics}, 27(3):221--237.

\bibitem[Gassner et~al., 2016a]{gassner_split_form_sbp_2016}
Gassner, G.~J., Winters, A.~R., and Kopriva, D.~A. (2016a).
\newblock Split form nodal discontinuous {G}alerkin schemes with
  summation-by-parts property for the compressible {Euler} equations.
\newblock {\em Journal of Computational Physics}, 327:39 -- 66.

\bibitem[Gassner et~al., 2016b]{Gassner2016}
Gassner, G.~J., Winters, A.~R., and Kopriva, D.~A. (2016b).
\newblock Split form nodal discontinuous {G}alerkin schemes with
  summation-by-parts property for the compressible {E}uler equations.
\newblock {\em Journal of Computational Physics}, 327(C):39--66.

\bibitem[Gassner et~al., 2016c]{gassner_entropy_shallow_water_2016}
Gassner, G.~J., Winters, A.~R., and Kopriva, D.~A. (2016c).
\newblock A well balanced and entropy conservative discontinuous {G}alerkin
  spectral element method for the shallow water equations.
\newblock {\em Applied Mathematics and Computation}, 272:291 -- 308.
\newblock Recent Advances in Numerical Methods for Hyperbolic Partial
  Differential Equations.

\bibitem[Gerritsen and Olsson, 1996]{Margot1996}
Gerritsen, M. and Olsson, P. (1996).
\newblock Designing an efficient solution strategy for fluid flows 1. {A}
  stable high order finite difference scheme and sharp shock resolution for the
  {E}uler equations.
\newblock {\em Journal of Computational Physics}, 129(2):245--262.

\bibitem[Hadri et~al., 2015]{shaheen_2}
Hadri, B., Kortas, S., Feki, S., Khurram, R., and Newby, G. (2015).
\newblock Overview of the {KAUST}'s {C}ray {X40} {S}ystem -- {S}haheen {II}.
\newblock {\em Proceedings of the Cray User Group Meeting}.

\bibitem[Hadri et~al., ]{hadri_ccpe_2019}
Hadri, B., Parsani, M., Hutchinson, M., Heinecke, A., Dalcin, L., and Keyes, D.
\newblock Performance study of sustained petascale direct numerical simulation
  on {C}ray {XC40} systems ({T}rinity, {S}haheen2 and {C}ori).
\newblock {\em Concurrency and Computation: Practice and Experience}.

\bibitem[Hesthaven and Warburton, 2008]{hesthaven_2008_nodal_dg}
Hesthaven, J.~S. and Warburton, T. (2008).
\newblock {\em Nodal discontinuous {Galerkin} methods: Algorithms, analysis,
  and applications}.
\newblock Texts in Applied Mathematics. Springer.

\bibitem[Honein and Moin, 2004]{honein2004higher}
Honein, A.~E. and Moin, P. (2004).
\newblock Higher entropy conservation and numerical stability of compressible
  turbulence simulations.
\newblock {\em Journal of Computational Physics}, 201(2):531--545.

\bibitem[Hughes et~al., 1986]{Hughes1986}
Hughes, T. J.~R., Franca, L.~P., and Mallet, M. (1986).
\newblock A new finite element formulation for computational fluid dynamics,
  {I}: symmetric forms of the compressible {N}avier--{S}tokes equations and the
  second law of thermodynamics.
\newblock {\em Computer Methods in Applied Mechanics and Engineering},
  54(2):223 -- 234.

\bibitem[Jagannathan and Donzis, 2016]{jagannathan2016reynolds}
Jagannathan, S. and Donzis, D.~A. (2016).
\newblock Reynolds and mach number scaling in solenoidally-forced compressible
  turbulence using high-resolution direct numerical simulations.
\newblock {\em Journal of Fluid Mechanics}, 789:669--707.

\bibitem[Katz and Sankaran, 2011]{katz2011}
Katz, A. and Sankaran, V. (2011).
\newblock Mesh quality effects on the accuracy of cfd solutions on unstructured
  meshes.
\newblock {\em Journal of Computational Physics}, 230(20):7670--7686.

\bibitem[Kennedy and Gruber, 2008]{kennedy2008reduced}
Kennedy, C.~A. and Gruber, A. (2008).
\newblock Reduced aliasing formulations of the convective terms within the
  {N}avier--{S}tokes equations for a compressible fluid.
\newblock {\em Journal of Computational Physics}, 227(3):1676--1700.

\bibitem[Kida and Orszag, 1990]{kida1990energy}
Kida, S. and Orszag, S.~A. (1990).
\newblock Energy and spectral dynamics in forced compressible turbulence.
\newblock {\em Journal of Scientific Computing}, 5(2):85--125.

\bibitem[Klose et~al., 2019]{klose2019robustness}
Klose, B.~F., Jacobs, G.~B., and Kopriva, D.~A. (2019).
\newblock On the robustness and accuracy of marginally resolved discontinuous
  {G}alerkin schemes for two dimensional {Navier--Stokes} flows.
\newblock In {\em AIAA Scitech 2019 Forum}, page 0780. American Institute of
  Aeronautics and Astronautics.

\bibitem[Knepley and Karpeev, 2009]{KnepleyKarpeev09}
Knepley, M.~G. and Karpeev, D.~A. (2009).
\newblock Mesh algorithms for {PDE} with {Sieve} {I}: {Mesh} distribution.
\newblock {\em Scientific Programming}, 17(3):215--230.

\bibitem[Kreiss and Scherer, 1974]{kreiss1974finite}
Kreiss, H.-O. and Scherer, G. (1974).
\newblock Finite element and finite difference methods for hyperbolic partial
  differential equations.
\newblock In de~Boor, C., editor, {\em Mathematical Aspects of Finite Elements
  in Partial Differential Equations}, pages 195--212, New York. Academic Press.

\bibitem[Lele, 1994]{lele1994compressibility}
Lele, S.~K. (1994).
\newblock Compressibility effects on turbulence.
\newblock {\em Annual Review of Fluid Mechanics}, 26(1):211--254.

\bibitem[Mengaldo et~al., 2015]{mengaldo_dealiasing_2015}
Mengaldo, G., De~Grazia, D., Moxey, D., Vincent, P., and Sherwin, S. (2015).
\newblock Dealiasing techniques for high-order spectral element methods on
  regular and irregular grids.
\newblock {\em Journal of Computational Physics}, 299:56--81.

\bibitem[Moura et~al., 2015]{moura2015dg}
Moura, R., Sherwin, S., and Peiro, J. (2015).
\newblock On dg-based iles approaches at very high reynolds numbers.
\newblock {\em Report, Research Gate}.

\bibitem[Nakagawa, 1988]{nakagawa1988effects}
Nakagawa, T. (1988).
\newblock Effects of an airfoil and shock waves on vortex shedding process
  behind a square cylinder.
\newblock {\em Acta mechanica}, 72(1-2):131--146.

\bibitem[Nolasco et~al., 2019]{nolasco_optim_metrics_2019}
Nolasco, I.~R., Dalcin, L., {Del Rey Fern\'andez}, D.~C., Zampini, S., and
  Parsani, M. (2019).
\newblock Optimized geometrical metrics satisfying free-stream preservation.
\newblock {\em Submitted to Computer \& Fluids}.

\bibitem[Nordstr{\"o}m and Bj{\"o}rck, 2001]{nordstrom2001finite}
Nordstr{\"o}m, J. and Bj{\"o}rck, M. (2001).
\newblock Finite volume approximations and strict stability for hyperbolic
  problems.
\newblock {\em Applied Numerical Mathematics}, 38(3):237--255.

\bibitem[Olsson and Oliger, 1994]{Olsson1994}
Olsson, P. and Oliger, J. (1994).
\newblock Energy and maximum norm estimates for nonlinear conservation laws.
\newblock Technical Report 94--01, The Research Institute of Advanced Computer
  Science.

\bibitem[Parsani et~al., 2016a]{parsani_ssdc_staggered_2016}
Parsani, M., Carpenter, M.~H., Fisher, T., and Nielsen, E. (2016a).
\newblock Entropy stable staggered grid discontinuous spectral collocation
  methods of any order for the compressible {N}avier--{S}tokes equations.
\newblock {\em SIAM Journal on Scientific Computing}, 38(5):A3129--A3162.

\bibitem[Parsani et~al., 2016b]{Parsani2016}
Parsani, M., Carpenter, M.~H., Fisher, T.~C., and Nielsen, E.~J. (2016b).
\newblock Entropy stable staggered grid discontinuous spectral collocation
  methods of any order for the compressible {N}avier--{S}tokes equations.
\newblock {\em SIAM Journal on Scientific Computing}, 38(5):A3129--A3162.

\bibitem[Parsani et~al., 2015a]{parsani_entropy_stable_interfaces_2015}
Parsani, M., Carpenter, M.~H., and Nielsen, E.~J. (2015a).
\newblock Entropy stable discontinuous interfaces coupling for the
  three-dimensional compressible {N}avier--{S}tokes equations.
\newblock {\em Journal of Computational Physics}, 290:132--138.

\bibitem[Parsani et~al., 2015b]{parsani_entropy_stability_solid_wall_2015}
Parsani, M., Carpenter, M.~H., and Nielsen, E.~J. (2015b).
\newblock Entropy stable wall boundary conditions for the three-dimensional
  compressible {N}avier--{S}tokes equations.
\newblock {\em Journal of Computational Physics}, 292:88--113.

\bibitem[Parsani et~al., 2015c]{parsani_wall_bc_entropy_2015}
Parsani, M., Carpenter, M.~H., and Nielsen, E.~J. (2015c).
\newblock Entropy stable wall boundary conditions for the three-dimensional
  compressible {N}avier--{S}tokes equations.
\newblock {\em Journal of Computational Physics}, 292(1):88--113.

\bibitem[Passot and Pouquet, 1987]{passot1987numerical}
Passot, T. and Pouquet, A. (1987).
\newblock Numerical simulation of compressible homogeneous flows in the
  turbulent regime.
\newblock {\em Journal of Fluid Mechanics}, 181:441--466.

\bibitem[Pazner and Persson, 2019]{pazner_es_line_dg_2019}
Pazner, W. and Persson, P.-O. (2019).
\newblock Analysis and entropy stability of the line-based discontinuous
  {Galerkin} method.
\newblock {\em Journal of Scientific Computing}, 80(1):376--402.

\bibitem[Pirozzoli and Grasso, 2004]{pirozzoli2004direct}
Pirozzoli, S. and Grasso, F. (2004).
\newblock Direct numerical simulations of isotropic compressible turbulence:
  influence of compressibility on dynamics and structures.
\newblock {\em Physics of Fluids}, 16(12):4386--4407.

\bibitem[Ranocha, 2019]{ranocha2019mimetic}
Ranocha, H. (2019).
\newblock Mimetic properties of difference operators: Product and chain rules
  as for functions of bounded variation and entropy stability of second
  derivatives.
\newblock {\em BIT Numerical Mathematics}, 59(2):547--563.

\bibitem[Ranocha et~al., 2016]{ranocha2016summation}
Ranocha, H., {\"O}ffner, P., and Sonar, T. (2016).
\newblock Summation-by-parts operators for correction procedure via
  reconstruction.
\newblock {\em Journal of Computational Physics}, 311:299--328.

\bibitem[Ranocha et~al., 2019]{ranocha2019relaxation}
Ranocha, H., Sayyari, M., Dalcin, L., Parsani, M., and Ketcheson, D.~I. (2019).
\newblock Relaxation {R}unge--{K}utta methods: Fully-discrete explicit
  entropy-stable schemes for the {E}uler and {N}avier--{S}tokes equations.
\newblock Accepted in SIAM Journal on Scientific Computing.

\bibitem[Ray et~al., 2016]{Ray2016}
Ray, D., Chandrashekar, P., Fjordhom, U.~S., and Mishra, S. (2016).
\newblock Entropy stable scheme on two-dimensional unstructured grids for
  {E}uler equations.
\newblock {\em Communications in Computational Physics}, 19(5):1111--1140.

\bibitem[Ristorcelli and Blaisdell, 1997]{ristorcelli1997consistent}
Ristorcelli, J.~R. and Blaisdell, G.~A. (1997).
\newblock Consistent initial conditions for the dns of compressible turbulence.
\newblock {\em Physics of Fluids}, 9(1):4--6.

\bibitem[Roache, 2001]{10.1115/1.1436090}
Roache, P.~J. (2001).
\newblock {Code Verification by the Method of Manufactured Solutions }.
\newblock {\em Journal of Fluids Engineering}, 124(1):4--10.

\bibitem[Roy, 2005]{roy2005}
Roy, C.~J. (2005).
\newblock Review of code and solution verification procedures for computational
  simulation.
\newblock {\em Journal of Computational Physics}, 205(1):131--156.

\bibitem[Sagaut and Cambon, 2008]{sagaut2008homogeneous}
Sagaut, P. and Cambon, C. (2008).
\newblock {\em Homogeneous turbulence dynamics}, volume~10.
\newblock Springer.

\bibitem[Samtaney et~al., 2001]{samtaney2001direct}
Samtaney, R., Pullin, D.~I., and Kosovi{\'c}, B. (2001).
\newblock Direct numerical simulation of decaying compressible turbulence and
  shocklet statistics.
\newblock {\em Physics of Fluids}, 13(5):1415--1430.

\bibitem[Sandham et~al., 2002]{Sandham2002}
Sandham, N.~D., Li, Q., and Yee, H.~C. (2002).
\newblock Entropy splitting for high-order numerical simulation of compressible
  turbulence.
\newblock {\em Journal of Computational Physics}, 178(2):307--322.

\bibitem[Shu, 1998]{shu1998essentially}
Shu, C.-W. (1998).
\newblock Essentially non-oscillatory and weighted essentially non-oscillatory
  schemes for hyperbolic conservation laws.
\newblock In {\em Advanced Numerical Approximation of Nonlinear Hyperbolic
  Equations}, pages 325--432. Springer.

\bibitem[Sj\"orn and Yee, 2018]{Bjorn2018}
Sj\"orn, B. and Yee, H.~C. (2018).
\newblock High order entropy conservative central schemes for wide ranges of
  compressible gas dynamics and {MHD} flows.
\newblock {\em Journal of Computational Physics}, 364:153--185.

\bibitem[Slotnick et~al., 2014]{nasa_2030_vision}
Slotnick, J., Khodadoust, A., Alonso, J., Darmofal, D., Gropp, W., Lurie, E.,
  and Mavriplis, D. (2014).
\newblock Cfd vision 2030 study: A path to revolutionary computational
  aerosciences.
\newblock {\em NASA-CR-2014-218178}.

\bibitem[S\"{o}derlind, 2003]{Soderlind2003}
S\"{o}derlind, G. (2003).
\newblock Digital filters in adaptive time-stepping.
\newblock {\em ACM Transactions on Mathematical Software}, 29(1):1--26.

\bibitem[S\"{o}derlind and Wang, 2006]{Soderlind2006}
S\"{o}derlind, G. and Wang, L. (2006).
\newblock Adaptive time-stepping and computational stability.
\newblock {\em Journal of Computational and Applied Mathematics},
  185(2):225--243.

\bibitem[Sv\"{a}rd et~al., 2018]{svard_entropy_stable_solid_wall_2018}
Sv\"{a}rd, M., Carpenter, M.~H., and Parsani, M. (2018).
\newblock Entropy stability and the no-slip wall boundary condition.
\newblock {\em SIAM Journal on Numerical Analysis}, 56(1):256--273.

\bibitem[Sv\"ard and Nordstr\"om, 2014]{Svard2014}
Sv\"ard, M. and Nordstr\"om, J. (2014).
\newblock Review of summation-by-parts schemes for
  initial-boundary-value-problems.
\newblock {\em Journal of Computational Physics}, 268(1):17--38.

\bibitem[Sv\"{a}rd and \"{O}zcan, 2014]{svard_entropy_stable_euler_wall_2014}
Sv\"{a}rd, M. and \"{O}zcan, H. (2014).
\newblock Entropy-stable schemes for the {E}uler equations with far-field and
  wall boundary conditions.
\newblock {\em Journal of Scientific Computing}, 58(1):61--89.

\bibitem[Tadmor, 2003]{Tadmor2003}
Tadmor, E. (2003).
\newblock Entropy stability theory for difference approximations of nonlinear
  conservation laws and related time-dependent problems.
\newblock {\em Acta Numerica}, 12:451--512.

\bibitem[Wang et~al., 2013]{wang_high_order_workshop_2013}
Wang, Z., Fidkowski, K., Abgrall, R., Bassi, F., Caraeni, D., Cary, A.,
  Deconinck, H., Hartmann, R., Hillewaert, K., Huynh, H., Kroll, N., May, G.,
  Persson, P.-O., Leer, B., and Visbal, M. (2013).
\newblock High-order cfd methods: current status and perspective.
\newblock {\em International Journal for Numerical Methods in Fluids},
  72(8):811--845.

\bibitem[Winters et~al., 2017]{Winters2017}
Winters, A.~R., Derigs, D., Gassner, G.~J., and Walch, S. (2017).
\newblock Uniquely defined entropy stable matrix dissipation operator for high
  {M}ach number ideal {MHD} and compressible {E}uler simulations.
\newblock {\em Journal of Computational Physics}, 332(1):274--289.

\bibitem[Winters and Gassner, 2015]{Winters2015}
Winters, A.~R. and Gassner, G.~J. (2015).
\newblock A comparison of two entropy stable discontinuous {G}alerkin spectral
  element approximations to the shallow water equations with non-constant
  topography.
\newblock {\em Journal of Computational Physics}, 301(1):357--376.

\bibitem[Winters et~al., 2018]{winters_split_form_2018}
Winters, A.~R., Moura, R.~C., Mengaldo, G., Gassner, G.~J., Walch, S., Peiro,
  J., and Sherwin, S.~J. (2018).
\newblock A comparative study on polynomial dealiasing and split form
  discontinuous {Galerkin} schemes for under-resolved turbulence computations.
\newblock {\em Journal of Computational Physics}, 372:1 -- 21.

\bibitem[Yamaleev and Carpenter, 2017]{Yamaleev2017}
Yamaleev, N.~K. and Carpenter, M.~H. (2017).
\newblock A family of fourth-order entropy stable non-oscillatory spectral
  collocation schemes for the 1-d {N}avier-{S}tokes equations.
\newblock {\em Journal of Computational Physics}, 331:90--107.

\bibitem[Yee et~al., 2000]{Yee2000}
Yee, H.~C., Vinokur, M., and Djomehri, M.~J. (2000).
\newblock Entropy splitting and numerical dissipation.
\newblock {\em Journal of Computational Physics}, 162(1):33--81.

\bibitem[Zhu et~al., 2013]{zhu_dg_weno_2013}
Zhu, J., Zhong, X., Shu, C.-W., and Qiu, J. (2013).
\newblock {Runge--Kutta} discontinuous {G}alerkin method using a new type of
  {WENO} limiters on unstructured meshes.
\newblock {\em Journal of Computational Physics}, 248:200 -- 220.

\end{thebibliography}

\end{document}